\documentclass[smallcondensed]{svjour3}        
\smartqed  

\usepackage{geometry}
\usepackage{xcolor}

\definecolor{matlab_blue}  {rgb}{0.0,    0.4470 ,   0.7410}
\definecolor{matlab_red}   {rgb}{1.0,    0.0,    0.0}
\definecolor{matlab_yellow}{rgb}{0.9290, 0.6940 ,   0.1250}
\definecolor{matlab_green} {rgb}{0.4660, 0.6740,    0.1880}
\usepackage{graphicx}
\usepackage{subfigure}
\usepackage{booktabs}

\usepackage{amsmath,amssymb}

\usepackage{accents}
\usepackage{stmaryrd}

\usepackage[numbers]{natbib}

\newcommand{\aver}[1]{\ensuremath{\{\!\{#1\}\!\}}}
\newcommand{\jump}[1]{\ensuremath{\left\llbracket #1 \right\rrbracket}}

\newcommand{\Sec}{Sec.~}
\newcommand{\R}{\mathbb R}

\newcommand*{\diff}{\mathop{}\!\mathrm{d}}

\DeclareMathAccent{\svec}{\mathord}{letters}{126}

\newcommand{\stvec}[1]{\mathbf #1}
\newcommand{\stvecg}[1]{\boldsymbol #1}
\newcommand{\ssvec}[1]{\overset{\leftrightarrow}{\stvec{#1}}}
\newcommand{\ssvecg}[1]{\overset{\leftrightarrow}{\stvecg{#1}}}
\newcommand{\cssvec}[1]{\overset{\leftrightarrow}{\tilde{\stvec{#1}}}}
\newcommand{\cssvecg}[1]{\overset{\leftrightarrow}{\tilde{\stvecg{#1}}}}
\newcommand{\smat}[1]{\underline{{#1}}}
\newcommand{\bmat}[1]{\mathcal{#1}}
\newcommand{\snabla}{\svec{\nabla}}
\newcommand{\interiorfaces}{\genfrac{}{}{0pt}{}{\mathrm{interior}}{\mathrm{faces}}}
\newcommand{\boundaryfaces}{\genfrac{}{}{0pt}{}{\mathrm{boundary}}{\mathrm{faces}}}

\newcommand{\GP}{\mathrm{GP}}

\newcommand{\SVV}{\mathrm{SVV}}
\newcommand{\diag}{\mathrm{diag}}
\newcommand{\rL}{\ensuremath{\mathrm{L}}}
\newcommand{\rF}{\ensuremath{\mathrm{F}}}
\newcommand{\rE}{\ensuremath{\mathrm{E}}}

\usepackage{scalerel,stackengine}
\stackMath
\newcommand\reallywidehat[1]{%
\savestack{\tmpbox}{\stretchto{%
  \scaleto{%
    \scalerel*[\widthof{\ensuremath{#1}}]{\kern.1pt\mathchar"0362\kern.1pt}%
    {\rule{0ex}{\textheight}}
  }{\textheight}%
}{2.4ex}}%
\stackon[-6.9pt]{#1}{\tmpbox}%
}
\usepackage{hyperref}
\newtheorem{prop}{Property}

\newcommand{\revone}[1]{#1}
\newcommand{\revtwo}[1]{#1}

\begin{document}

\title{An entropy stable spectral vanishing viscosity for discontinuous Galerkin schemes:
application to shock capturing and LES models}

\titlerunning{Entropy stable SVV for DG schemes: application to shock capturing/LES models}

\author{Andr\'es Mateo-Gab\'in \and Juan Manzanero \and Eusebio Valero}

\institute{Andr\'es Mateo Gab\'in (\email{andres.mgabin@upm.es}) \and
    Juan Manzanero \and Eusebio Valero \at ETSIAE-UPM - School of Aeronautics,
        Universidad Polit\'ecnica de Madrid. Plaza Cardenal Cisneros 3, 28040 Madrid, Spain
    \and
    Andr\'es Mateo-Gab\'in \and Juan Manzanero \and Eusebio Valero \at Center for Computational
        Simulation, Universidad Polit\'ecnica de Madrid, Campus de Montegancedo, Boadilla del Monte,
        28660, Madrid, Spain.
}
\date{Received: date / Accepted: date}

\maketitle

\begin{abstract}

We present a stable spectral vanishing viscosity for discontinuous Galerkin schemes, with
applications to turbulent and supersonic flows. The idea behind the SVV is to spatially filter the
dissipative fluxes, such that it concentrates in higher wavenumbers, where the flow is typically
under--resolved, leaving low wavenumbers dissipation--free. Moreover, we derive a stable
approximation of the Guermond--Popov fluxes with the Bassi--Rebay 1 scheme, used to introduce
density regularization in shock capturing simulations. This filtering uses a Cholesky decomposition
of the fluxes that ensures the entropy stability of the scheme, which also includes a stable
approximation of boundary conditions for adiabatic walls. For turbulent flows, we test the method
with the three--dimensional Taylor--Green vortex and show that energy is correctly dissipated, and
the scheme is stable when a kinetic energy preserving split--form is used in combination with a low
dissipation Riemann solver. Finally, we test the shock capturing capabilities of our method with
the Shu--Osher and the supersonic forward facing step cases, obtaining good results without spurious
oscillations even with coarse meshes.


\end{abstract}

\section{Introduction}

The central role that differential equations play in science and engineering explain the great
interest and the efforts made to develop ever more efficient and accurate methods. In the field of
fluid mechanics, high--order discontinuous Galerkin methods~\cite{Reed1977,Cockburn2000} seem to
lay in a sweet spot, being suitable for complex geometries and easily allowing the introduction of
h/p mesh adaptation to spectrally increase the
accuracy~\cite{Woopen2014,Kompenhans2016,Friedrich2018,Ntoukas2021}. Nevertheless, the numerical
dissipation they can introduce is limited and, even with the addition of inter--element numerical
fluxes~\cite{Toro2009}, it is not sufficient for challenging cases such as supersonic flows, where
the Navier--Stokes equations allow discontinuous solutions that need higher rates of dissipation
when approximated by high--order polynomials. In this situation it is desirable that the solver is
able to maintain an ever decreasing entropy, according to the second law of thermodynamics.

The mathematical entropy is a well established framework to represent conservation laws. This
mathematical entropy is arbitrarily defined such that is ever
decreasing~\cite{Harten1983,Tadmor1986,Tadmor2003}. Hence, in kinetic--energy preserving schemes,
the mathematical entropy is the kinetic energy, whereas in thermodynamic entropy preserving schemes,
it is minus the thermodynamic entropy (which is ever increasing). In the case of the Discontinuous
Galerkin Spectral Element Method (DGSEM)~\cite{Kopriva2009} with Gauss--Lobatto (GL) points, the
polynomial derivative operator satisfies the discrete Summation--By--Parts
Simultaneous--Approximation--Term (SBP--SAT) property~\cite{Strand1994,Abgrall2020}, which has led
the development of split--form schemes that preserve the discrete entropy~\cite{Kravchenko1997}.
However, in the context of supersonic flows and shock capturing, where the flow variables can
experience large oscillations, the use of an entropy preserving method does not guarantee its
stability. This is because the entropy stability condition rests on a positive density, which is
usually not the case as a result of the violent density oscillations at the vicinity of shocks. One
way to alleviate this problem is the use of shock capturing methods that introduce some sort of
density regularization, while keeping discrete entropy preservation. This point is a focus of active
research and the literature contains several strategies that can be divided in two main branches.
The first group makes use of the discretization to introduce the required additional dissipation. By
lowering the order of the approximation, the dissipation of the numerical fluxes is introduced at
lower wavenumbers, reducing the slope of the wave and smoothing the solution~\cite{Hennemann2021}.
The other approach consists in adding dissipation directly in the driving equations, usually through
a second order elliptic term~\cite{Persson2006,Klockner2011}. Nevertheless, the discontinuity limits
the attainable accuracy and elements containing them can never reach spectral
convergence~\cite{Gottlieb2001}. In this work we focus on the second method and we derive two
different approaches; first, a thermodynamic entropy stable scheme that uses a spectral vanishing
viscosity (SVV) based on the fluxes proposed by Guermond and Popov~\cite{Guermond2014}, and second,
a kinetic energy preserving scheme with its fluxes based on the physical Navier--Stokes viscosity.
The use of the SVV makes it possible to modulate the dissipation in the wavenumber domain, and
allows us to continuously shift from a low (second order) to a high--order dissipation.

Moreover, the SVV approach is also suitable for under--resolved turbulent flows, where the mesh
cannot capture all the scales involved and the high--order DGSEM might not dissipate the correct
amount of energy~\cite{Karamanos2000,Kirby2006,Moura2016}. Several strategies have also been
proposed for these types of flows, being the Reynolds Averaged Navier--Stokes (RANS)
equations~\cite{Lodares2021} and Large Eddy Simulations (LES) among the most extended. Implicit LES
(iLES) methods benefit from the scheme's numerical dissipation~\cite{Manzanero2020}, while explicit
LES methods model the effect of the small scales that cannot be accurately resolved into aditional
dissipation. In this work, we follow~\cite{Kirby2002,Pasquetti2008} and couple the SVV and LES
aproaches so that the intensity is computed with the LES--Smagorinsky model, and it is then filtered
by the SVV. However, we modify the SVV discretization such that the resulting scheme is kinetic
energy preserving.

The structure of the rest of the paper is as follows: we briefly describe the compressible
Navier--Stokes equations in \Sec\ref{sec:NSE}, where we also highlight the entropy variables and
introduce the artificial viscosity fluxes that will be used later. We then show the foundation of
the DGSEM in \Sec\ref{sec:DGSEM} and, in \Sec\ref{sec:DG:SVV}, we detail the novel filtering
technique that we use in the SVV formulation. \Sec\ref{sec:stability} contains a detailed analysis
of the stability of the entropy formulation with artificial viscosity, and in
\Sec\ref{sec:vonNeumann} we perform von Neumann analyses of the SVV applied to a one--dimensional
advection--diffusion equation to preview its dissipation properties. Finally, \Sec\ref{sec:results}
shows some numerical results with test cases in one to three dimensions: the Shu--Osher shock tube,
a supersonic forward facing step and the inviscid Taylor--Green vortex.

\subsection{Notation}\label{subsec:intro:notation}

The mathematical formulation of partial differential equations introduces the concepts of space
vectors in~$\R^3$ like~$\svec{x}=(x,y,z)^T$, and state vectors in~$\R^5$, such
as~$\stvec{q}=(\rho,\rho u,\rho v,\rho w,\rho e)^T$, to gather the field variables. We adopt the
notation in~\cite{Gassner2018}, where a block vector represents an entity contained in the
space~$\R^3\times\R^5$, e.g., fluxes, as three state vectors stacked on top of each other,
\begin{equation}
    \ssvec{f} = \left(\begin{array}{c}
            \stvec{f}_1 \\
            \stvec{f}_2 \\
            \stvec{f}_3
        \end{array}\right) = \left(\begin{array}{c}
            \stvec{f} \\
            \stvec{g} \\
            \stvec{h}
        \end{array}\right).
\end{equation}
We can now extend the usual multiplication operators of calculus to vectors of different nature by
acting on each subspace separately,
\begin{equation}
    \ssvec{f}\cdot\ssvec{g} = \sum_{i=1}^{3}\stvec{f}_{i}^T\stvec{g}_{i},~~
    \svec{g}\cdot\ssvec{f} = \sum_{i=1}^{3}g_i \stvec{f}_{i},~~
    \svec{g}\stvec{f} = \left(\begin{array}{c}
        g_1\stvec{f} \\
        g_2\stvec{f} \\
        g_3\stvec{f}
    \end{array}\right),
\label{eq:governing:space-state-block-ops}
\end{equation}
and, by setting~$\svec{g}=\svec{\nabla}$, we obtain the expressions for the divergence and gradient
operators,
\begin{equation}
    \svec{\nabla}\cdot\ssvec{f} = \sum_{i=1}^{3}\frac{\partial\stvec{f}_{i}}{\partial x_i},~~
    \svec{\nabla}\stvec{q} = \left(\begin{array}{c}
        \stvec{q}_{x} \\
        \stvec{q}_{y} \\
        \stvec{q}_z
    \end{array}\right).
\end{equation}
Likewise, we define state matrices,~$\smat{A}$, of size~$5\times5$ to operate on the field
variables, and block matrices as a combination of state matrices,
\begin{equation}
    \bmat{A} = \left(\begin{array}{ccc}
        \smat{A}_{11} & \smat{A}_{12} & \smat{A}_{13} \\
        \smat{A}_{21} & \smat{A}_{22} & \smat{A}_{23} \\
        \smat{A}_{31} & \smat{A}_{32} & \smat{A}_{33}
    \end{array}\right).
\label{eq:notation:block-matrix}
\end{equation}
and thus, the product of matrices and vectors of the same type is well defined. We can also
construct block versions of space matrices as follows. For example, consider the product
\begin{equation}
    \svec{g} = \tens{M}\svec{f}.
\label{eq:notation:rotation}
\end{equation}
that involves a space matrix~$\tens{M} \in \R^{3 \times 3}$. Its associated block
matrix~$\bmat{M} \in \R^{15 \times 15}$ is created so that it is equivalent to
apply~\eqref{eq:notation:rotation} to each of the variables of the state components of the block
vector,~$\ssvec{f}$,
\begin{equation}
    \ssvec{g} = \mathcal M \ssvec{f},\quad
    \bmat{M} = \left(\begin{array}{ccc}
        \tens{M}_{11}\smat{I}_{5} & \tens{M}_{12}\smat{I}_{5} & \tens{M}_{13}\smat{I}_{5} \\
        \tens{M}_{21}\smat{I}_{5} & \tens{M}_{22}\smat{I}_{5} & \tens{M}_{23}\smat{I}_{5} \\
        \tens{M}_{31}\smat{I}_{5} & \tens{M}_{32}\smat{I}_{5} & \tens{M}_{33}\smat{I}_{5}
    \end{array}\right),
\label{eq:notation:block-matrix-from-space-matrix}
\end{equation}
For more details, see \cite{Gassner2018}. We finally represent integrals in the domain~$\Omega$
in bracket notation,~$\langle\cdot\rangle$, being~$\langle\cdot,\cdot\rangle$ the inner product in
the same region,
\begin{equation}
    \left\langle f\right\rangle = \int_{\Omega} f\diff \svec{x},\quad
    \left\langle f,g\right\rangle = \int_{\Omega}fg\diff\svec{x}.
\end{equation}

\section{The compressible Navier--Stokes equations with artificial viscosity}\label{sec:NSE}

In this section we describe the compressible Navier--Stokes equations. For the purpose of this
work, we will complement the set of equations with additional dissipative terms, of which we
consider two: locally increasing the molecular viscosity (Boussinesq's approximation), and the
entropy stable dissipative flux derived by Guermond and Popov in~\cite{Guermond2014}. These
fluxes are detailed in \Sec\ref{sec:NSE:AD}. Finally, in \Sec\ref{sec:NSE:ES} we study the kinetic
energy and entropy stability of the schemes derived.

The compressible Navier--Stokes equations with artificial viscosity are a set of non--linear
advection--diffusion equations,
\begin{equation}
    \stvec{q}_{t} + \svec{\nabla}\cdot\ssvec{f}_{e} =
        \svec{\nabla}\cdot\ssvec{f}_{v} + \svec{\nabla}\cdot\ssvec{f}_{a}.
\label{eq:NSE:adv-diff}
\end{equation}
The state vector is~$\stvec{q}=\left(\rho,\rho\svec{u},\rho e\right)^T$ where~$\rho$ is the
density,~$\svec{u}=\left(u,v,w\right)\revone{^T}$ is the velocity, and~$\rho e$ is the total energy,
\begin{equation}
    \rho e = \rho e_{i} + \frac{1}{2}\rho |\svec{u}|^{2},~~
    e_{i} = \frac{p/\rho}{\gamma-1},~~
    p = \rho R T,
\end{equation}
being~$e_{i}$ the internal energy,~$p$ the pressure,~$T$ the temperature,~$R$ the ideal gas
constant, and~$\gamma$ the specific heat ratio. The inviscid fluxes~$\ssvec{f}_{e}$ depend on the
state vector,
\begin{equation}
    \stvec{f}_{e} = \left(\begin{array}{c}
        \rho u     \\
        \rho u^2+p \\
        \rho u v   \\
        \rho u w   \\
        \rho h u
    \end{array}\right),~~
    \stvec{g}_{e} = \left(\begin{array}{c}
        \rho v       \\
        \rho uv      \\
        \rho v^2 + p \\
        \rho vw      \\
        \rho hv
    \end{array}\right),~~
    \stvec{h}_{e} = \left(\begin{array}{c}
        \rho w     \\
        \rho uw    \\
        \rho vw    \\
        \rho w^2+p \\
        \rho hw
    \end{array}\right),
\end{equation}
where~$h = e + p/\rho$ is the total enthalpy. The viscous fluxes~$\ssvec{f}_{v}$ are
\begin{equation}
    \stvec{f}_{v} = \left(\begin{array}{c} 0 \\
        \tau_{11} \\
        \tau_{21} \\
        \tau_{31} \\
        \svec{\tau}_{1}\cdot\svec{u} + q_{1}
    \end{array}\right),~~
    \stvec{g}_{v} = \left(\begin{array}{c}
        0         \\
        \tau_{12} \\
        \tau_{22} \\
        \tau_{32} \\
        \svec{\tau}_{2}\cdot\svec{u} + q_{2}
    \end{array}\right),~~
    \stvec{h}_{v} = \left(\begin{array}{c}
        0         \\
        \tau_{13} \\
        \tau_{23} \\
        \tau_{33} \\
        \svec{\tau}_{3}\cdot\svec{u} + q_{3}
    \end{array}\right),
\end{equation}
with~$\tau_{ij} = \mu\left(\partial u_i /\partial x_j+\partial u_j/\partial x_i -
\frac{2}{3}\svec{\nabla}\cdot\svec{u}~\delta_{ij} \right)$ the stress
tensor,~$\svec{\tau}_{i} = (\tau_{1i},\tau_{2i},\tau_{3i})$ the stress in the three spatial
directions, and~$\svec{q}$ the heat flux,
\begin{equation}
    \svec{q} = \kappa \snabla T.
\end{equation}
The coefficients~$\mu$ and~$\kappa$ are the molecular viscosity and thermal conductivity, related
through the Prandtl number,
\begin{equation}
    \kappa = \theta \mu R,~~\theta = \frac{\gamma}{(\gamma-1)\Pr}
\label{eq:NSE:kappatheta}
\end{equation}
Finally, two different approaches to include additional dissipation are considered:
\begin{enumerate}
\item Increasing the molecular viscosity,~$\mu_a, \ssvec{f}_a = \ssvec{f}_v(\mu_a)$,
\item The Guermond--Popov flux,~$\ssvec{f}_a = \ssvec{f}_{\GP}(\mu_a, \alpha_a)$, developed
in~\cite{Guermond2014} that introduces two additional viscosity parameters,~$\alpha_a$ and~$\mu_a$,
which will be explained in more detail in~\Sec\ref{sec:NSE:GP},
\begin{equation}
    \ssvec{f}_{\mathrm{GP}} = \alpha_{a}\left(\begin{array}{c}
        \svec{\nabla}\rho \\
        \svec{\nabla}\rho\otimes\svec{u} \\
        \svec{\nabla}\left(\rho e_{i}\right) + \frac{1}{2}|\svec{u}|^2\svec{\nabla}\rho
    \end{array}\right) + \mu_{a} \left(\begin{array}{c}
        0 \\
        \rho \svec{\nabla}^{s}\svec{u}\\
        \rho \svec{u}\cdot\svec{\nabla}^s\svec{u}
    \end{array}\right),\quad
    \svec{\nabla}^{s}\svec{u} = \frac{1}{2} \left(\svec{\nabla}\svec{u} +
    \svec{\nabla}\svec{u}^{T}\right).
\label{eq:NSE:GP-fluxes}
\end{equation}
\end{enumerate}

\subsection{Entropy pairs and entropy variables}\label{subsec:entropy}

Now, in order to further analyze the entropy stability properties of our formulation in later
sections, we introduce the concept of mathematical entropy. Starting from a generic conservation law
in three dimensions,
\begin{equation}
    \stvec{q}_t + \snabla \cdot \ssvec{f} = 0,
\label{eq:NSE:conservation}
\end{equation}
we define the entropy~\cite{Harten1983,Tadmor1986,Tadmor2003,Friedrichs1971},~$\mathcal{E}$, as a
convex function on the variables~$\stvec{q}$ that also satisfies
\begin{equation}
    \left(\mathcal{E}_{\stvec{q}}\right)^T \ssvec{f}_{\stvec{q}} =
        \svec{f}^{\mathcal E}_{\stvec{q}}, \quad
    \ssvec{f}_\stvec{q} = \left(\stvec{f}_\stvec{q}, \stvec{g}_\stvec{q},
        \stvec{h}_\stvec{q}\right)^T,
\label{eq:NSE:entropypair}
\end{equation}
for the entropy--entropy flux pair~$(\mathcal E, \svec{f}^{\mathcal E})$. The derivative
of~$\mathcal E$ is given the name of entropy variables,
\begin{equation}
    \stvec{w} = \mathcal{E}_{\stvec{q}},
\label{eq:NSE:entropyvars}
\end{equation}
which also introduces the mapping~$\stvec{q} \to \stvec{w}$. After
multiplying~\eqref{eq:NSE:conservation} from the left by~$\stvec{w}^T$ and introducing the
conditions of~\eqref{eq:NSE:entropypair}, the initial conservation law also expresses the
conservation of the mathematical entropy,~$\mathcal E$,
\begin{equation}
    \mathcal{E}_t + \snabla \cdot \svec{f}^{\mathcal E} = 0.
\label{eq:NSE:entropy1D}
\end{equation}
Going back to~\eqref{eq:NSE:adv-diff}, the terms to the left of the equal sign correspond to the
Euler equations which, being hyperbolic, add no dissipation and can be reduced to the
form~\eqref{eq:NSE:entropy1D} with the use of a certain family of entropy
functions~\cite{Tadmor1986}. The dissipative fluxes at the right hand side
of~\eqref{eq:NSE:adv-diff} are also modified by the introduction of these entropy variables. We will
show in the following sections that both,~$\ssvec{f}_v$ and~$\ssvec{f}_a$, are symmetrized by this
change of variables in the sense that,
\begin{equation}
    \ssvec{f}_v = \bmat{B}_v(\stvec{q}) \snabla\stvec{w}, \quad
    \ssvec{f}_a = \bmat{B}_a(\stvec{q}) \snabla\stvec{w},
\label{eq:NSE:symmetric_diss}
\end{equation}
with~$\bmat{B}_v$ and~$\bmat{B}_a$ symmetric matrices. We also define the dissipation introduced by
any flux of the form~\eqref{eq:NSE:symmetric_diss} as,
\begin{equation}
    D = \left(\svec{\nabla}\stvec{w}\right)^T\bmat{B}~\svec{\nabla}\stvec{w},
\label{eq:NSE:entropy_diss}
\end{equation}
which will satisfy~$D\geqslant 0$ if~$\bmat{B}$ is positive semi--definite. In this case, after
multiplying~\eqref{eq:NSE:adv-diff} from the left by~$\stvec{w}^T$ and integrating over~$\Omega$, we
get a new conservation law for the integrated value of the entropy,~$\langle\mathcal E\rangle$,
\begin{equation}
    \langle\mathcal E_t\rangle +
        \int_{\Omega} \snabla\cdot\svec{f}^{\mathcal E} \diff V =
        \int_{\partial\Omega} \stvec{w}^T\left(\ssvec{f}_v+\ssvec{f}_a\right)\cdot\svec{n}\diff S -
        \left\langle\snabla\stvec{w}, \ssvec{f}_v\right\rangle -
        \left\langle\snabla\stvec{w}, \ssvec{f}_a\right\rangle,
\end{equation}
that can be further simplified when introducing~\eqref{eq:NSE:symmetric_diss}
and~\eqref{eq:NSE:entropy_diss} into the volume flux terms,
\begin{equation}
    \langle\mathcal E_t\rangle + \int_{\partial\Omega}
        \left(\svec{f}^{\mathcal E} - \stvec{w}^T\ssvec{f}_v - \stvec{w}^T\ssvec{f}_a\right)
        \cdot\svec{n}\diff S = -\langle D_v\rangle - \langle D_a\rangle
        \leqslant 0.
\label{eq:NSE:entropy_inequality}
\end{equation}
The last inequality comes from the fact that the dissipation introduced by the viscous and
artificial fluxes considered in this work is always positive, as it is explained in
\Sec\ref{sec:NSE:AD}. Therefore, the entropy remains bounded during the flow evolution at a
continuous level. Since the Navier--Stokes equations verify~\eqref{eq:NSE:entropy_inequality}, it is
desirable to have it satisfied also at a discrete level because entropy stable schemes can be shown
to be non--linearly stable~\cite{Merriam1989}, and therefore, very robust.

In this work we specifically use two different mathematical entropies:
\begin{enumerate}
\item The kinetic energy,~$\mathcal K = \rho |\svec{u}|^2 / 2$, with entropy variables,
\begin{equation}
    \stvec{w}^{\mathcal K} = \frac{\partial \mathcal K}{\partial \stvec{q}} =
        \left(-\frac{|\svec{u}|^2}{2}, u, v, w, 0\right)^T.
\label{eq:NSE:kinetic-energy-variables}
\end{equation}
\item The thermodynamic entropy,~$\mathcal S = -\rho s / (\gamma-1)$,
with~$s = \ln p - \gamma \ln \rho$ and,
\begin{equation}
    \stvec{w}^{\mathcal S} = \frac{\partial \mathcal S}{\partial \stvec{q}} =
        \left(\frac{\gamma-s}{\gamma-1}-\frac{\rho |\svec{u}|^{2}}{2p},
        \frac{\rho u}{p}, \frac{\rho v}{p}, \frac{\rho w}{p}, -\frac{\rho}{p}\right)^T.
\label{eq:NSE:entropy-variables}
\end{equation}
\end{enumerate}

\subsection{Entropy stable dissipation}\label{sec:NSE:AD}

From~\eqref{eq:NSE:entropy_diss} and~\eqref{eq:NSE:entropy_inequality} we see that the viscous
fluxes are a natural way to introduce dissipation to the entropy equation and thus, we only need to
prove that these fluxes are actually adding a positive dissipative term by showing that the
matrices~$\bmat{B}$ are positive semi--definite.

This has been already shown for the Navier--Stokes viscous fluxes in~\cite{Fisher2013}. In any case,
in this section we complete the previous analysis for the particular cases summarized in
Table~\ref{tab:NSE:AD:recap}.
\begin{table}
\centering
\begin{tabular}{lcccc}
    \toprule
    Flux type & Kinetic energy & Thermodynamic entropy & $\mu_{\alpha}$ & $\alpha_{\alpha}$ \\
    \midrule
    Navier--Stokes  & \checkmark & \checkmark & constant / LES & --       \\
    Guermond--Popov & --         & \checkmark & constant       & constant \\
    \bottomrule
\end{tabular}
\caption{Combinations of artificial viscous fluxes and entropy variables used throughout this work.}
\label{tab:NSE:AD:recap}
\end{table}

\subsubsection{The kinetic energy variables viscous flux}\label{sec:NSE:kin-visc}

The Navier--Stokes viscous flux can be expressed in matrix--vector form with~$\bmat{B}_v^{\mathcal
K} = \mu\bmat{C}_v^{\mathcal K}$ (see Appendix~\ref{sec:NS_kinetic_matrices} for the actual
definition of the matrix) so, in addition to the physical dissipation, more can be added by
augmenting the molecular viscosity,~$\bmat{B}_a^{\mathcal K} = \mu_a\bmat{C}_v^{\mathcal K}$. Its
associated dissipation is
\begin{equation}
\begin{split}
    D_a^{\mathcal K} &= \mu_a\left(\svec{\nabla}\stvec{w}^{\mathcal K}\right)^T
        \bmat{C}_{v}^{\mathcal K}\svec{\nabla}\stvec{w}^{\mathcal K} =
        \mu_a\svec{\nabla}\svec{u}\ : \left[\left(\svec{\nabla}\svec{u}\right)^{T}
        + \svec{\nabla}\svec{u} - \frac{2}{3}\left(\svec{\nabla} \cdot
        \svec{u}\right)I_{3}\right]\\
        &= \mu_a\left(2|\tens{S}|^{2} - \frac{2}{3}\left(\svec{\nabla} \cdot
        \svec{u}\right)^2\right) \geqslant 0.
\end{split}
\end{equation}

\subsubsection{The thermodynamic entropy variables viscous flux}\label{sec:NSE:ther-visc}

Since the matrix corresponding to the Navier--Stokes viscous flux is also linearly dependent
on~$\rho$ with this set of variables,~$\bmat{B}_v^{\mathcal S} = \frac{\mu p}{\rho}
\tilde{\bmat{B}}_v^{\mathcal S}$, artificial viscosity can also be added by augmenting the
molecular viscosity,~$\bmat{B}_a^{\mathcal S} = \frac{\mu_a p}{\rho}
\tilde{\bmat{B}}_v^{\mathcal S}$. However in this case, the proof for non--linear matrices rests on
a Cholesky decomposition of the matrix~$\tilde{\bmat B}_{a}^{\mathcal S} =
\bmat{L}_{v}^{\mathcal S,T}\bmat{D}_{v}^{\mathcal S}\bmat{L}_{v}^{\mathcal S}$, hence,
\begin{equation}
    D_a^{\mathcal S} = \frac{\mu_a p}{\rho}\svec{\nabla}\stvec{w}^{\mathcal S,T}
        \tilde{\bmat{B}}_{v}^{\mathcal S}\svec{\nabla}\stvec{w}^{\mathcal S} =
        \frac{\mu_a p}{\rho}\left(\bmat{L}_{v}^{\mathcal S}
        \svec{\nabla}\stvec{w}^{\mathcal S}\right)^{T}{\bmat{D}}_{v}^{\mathcal S}
        \left(\bmat{L}_{v}^{\mathcal S} \svec{\nabla}\stvec{w}^{\mathcal S}\right) \geqslant 0,
\end{equation}
if all the terms in the diagonal matrix, the artificial viscosity~$\mu_a$, and the
temperature~$T=p/\rho R$ are positive. The diagonal matrix was given in~\cite{Fisher2013},
\begin{equation}
    \bmat{D}_{v}^{\mathcal S} = \mathrm{diag}\left( 0, \frac{4}{3}, 1, 1, \frac{\theta p}{\rho},
        0, 0, 1, 1, \frac{\theta p}{\rho}, 0, 0, 0, 0, \frac{\theta p}{\rho}\right),
\end{equation}
where~$\theta$ takes the value given in~\eqref{eq:NSE:kappatheta}. For further information, we
include the precise form of the matrices~$\tilde{\bmat{B}}_v^{\mathcal S}$
and~$\bmat{L}_v^{\mathcal S}$ in Appendix~\ref{sec:NS_entropy_matrices}.

\subsubsection{The Guermond--Popov flux}\label{sec:NSE:GP}

The last flux studied in this work was derived in~\cite{Guermond2014}, and introduces two
parameters~\eqref{eq:NSE:GP-fluxes}. The first one,~$\alpha_{a}$, introduces dissipation that
originates in the continuity equation, proportional to the gradient of the density. This dissipation
is carried through the rest of the equations to get an entropy stable dissipation. The second
one,~$\mu_{a}$, further introduces dissipation that originates in the momentum equation and whose
work is carried to the energy equation, in a similar fashion to the physical dissipation of the
Navier--Stokes equations.

To study the dissipation introduced by this flux, it is written in the
form~$\ssvec{f}_{\GP}=\bmat{B}_{\GP}\svec{\nabla}\stvec{w}^{\mathcal S}$, using the thermodynamic
entropy variables~\eqref{eq:NSE:entropy-variables}. Therefore, the matrix
$\bmat{B}_{\GP}$ is written as the sum of the separate contributions of~$\alpha_{a}$ and~$\mu_{a}$
\begin{equation}
    \bmat{B}_{\GP} = \alpha_{a}\rho\bmat{B}_{\GP}^{\alpha} + \mu_{a}p\bmat{B}_{\GP}^{\mu},
\label{eq:NSE:ad:GP-B}
\end{equation}
where~$\bmat{B}_{\GP}^{\alpha}$ and~$\bmat{B}_{\GP}^{\mu}$ are given in
Appendix~\ref{sec:GP_entropy_matrices}. As in the Navier--Stokes fluxes, having an entropy stable
dissipation rests on the positive--definiteness of~\eqref{eq:NSE:ad:GP-B}. In a similar fashion, we
perform a Cholesky decomposition of~$\bmat{B}_{\GP}=\bmat{L}_{\GP}^{T}\bmat{D}_{\GP}\bmat{L}_{\GP}$,
with
\begin{equation}
    \bmat{D}_{\GP} = \mathrm{diag}\left(\alpha_{a}\rho, \mu_{a}p, \frac{1}{2}\mu_{a}p,
        \frac{1}{2}\mu_{a}p, \alpha_{a} \rho,\alpha_{a} \rho, 0, \mu_{a}p, \frac{1}{2}\mu_{a}p,
        \alpha_{a} \rho, \alpha_{a} \rho, 0, 0, \mu_{a}p, \alpha_{a} \rho\right),
\label{eq:NSE:GP-Dmat}
\end{equation}
thus confirming that
\begin{equation}
    D_a^{\GP} = \svec{\nabla}\stvec{w}^{\mathcal S,T}
        \tilde{\bmat{B}}_{\GP}\svec{\nabla}\stvec{w} =
        \left(\bmat{L}_{\GP}\svec{\nabla}\stvec{w}^{\mathcal S}\right)^{T}
        {\bmat{D}}_{\GP} \left(\bmat{L}_{\GP}\svec{\nabla}\stvec{w}^{\mathcal S}\right)\geqslant 0,
\end{equation}
as long as~$\rho,p \geqslant 0$. More details of the formulation are given in
Appendix~\ref{sec:GP_entropy_matrices}.

\subsection{Continuous entropy analysis}\label{sec:NSE:ES}

We now compute the final expression of the entropy equation~\eqref{eq:NSE:entropy_inequality} for
the two entropy variables that we consider. After multiplying~\eqref{eq:NSE:adv-diff} by the entropy
variables and integrating the result over~$\Omega$ we get,
\begin{equation}
    \left\langle \stvec{q}_{t},\stvec{w}\right\rangle +
        \left\langle\svec{\nabla}\cdot\ssvec{f}_{e},\stvec{w}\right\rangle =
        \left\langle\svec{\nabla}\cdot\left(\ssvec{f}_{v}+\ssvec{f}_{a}\right),
        \stvec{w}\right\rangle.
\label{eq:NSE:contentropy}
\end{equation}
The final form of the different terms depends on the choice of entropy variables. We know
from~\eqref{eq:NSE:entropyvars} that the first term gives the evolution of the entropy,
\begin{equation}
    \left\langle \stvec{q}_{t},\stvec{w}^{\mathcal K}\right\rangle =
        \left\langle\mathcal K_{t}\right\rangle,~~\text{or}~~
    \left\langle\stvec{q}_{t},\stvec{w}^{\mathcal S}\right\rangle =
        \left\langle\mathcal S_{t}\right\rangle,
\end{equation}
and also from~\eqref{eq:NSE:entropypair} we can derive that the second term is the integral of the
entropy flux over the boundary~$\partial \Omega$. However, the case of the kinetic energy needs a
special treatment since it does not meet all the requirements set in~\Sec\ref{subsec:entropy}.
Specifically, the entropy flux expression~\eqref{eq:NSE:entropypair} is not true and this second
term of~\eqref{eq:NSE:contentropy} represents the kinetic energy flux through the boundaries plus a
volume term,~$\mathcal{W}_p$, usually associated to the work of the pressure. A more detailed
discussion about this term and its role in kinetic energy preserving schemes can be found
in~\cite{Jameson2008,Gassner2014,gassner2016split,Lodares2021}. Using~\eqref{eq:NSE:entropypair} to
compute the entropy fluxes,
\begin{equation}
\begin{split}
    \left\langle \svec{\nabla}\cdot\ssvec{f}_{e},\stvec{w}^{\mathcal K}\right\rangle &=
        \int_{\partial\Omega}\frac{1}{2}\rho|\svec{u}|^2\svec{u}\cdot\svec{n} \diff S +
        \left\langle\svec{u},\svec{\nabla}p\right\rangle =
        \int_{\partial\Omega}\svec{f}^{\mathcal K}_{e}\cdot\svec{n}\diff S + \mathcal W_{p},\\
    \left\langle\svec{\nabla}\cdot\ssvec{f}_{e},\stvec{w}^{\mathcal S}\right\rangle &=
        \int_{\partial \Omega}\mathcal S\svec{u}\cdot\svec{n}\diff S =
        \int_{\partial\Omega}\svec{f}^{\mathcal S}_{e}\cdot\svec{n}\diff S,
\end{split}
\end{equation}
Finally, the last term gives the dissipative terms after introducing~\eqref{eq:NSE:entropy_diss} for
the physical and artificial viscous fluxes,
\begin{equation}
    \left\langle\svec{\nabla}\cdot\left(\ssvec{f}_{v}+\ssvec{f}_{a}\right),\stvec{w}\right\rangle =
        \int_{\partial\Omega}\stvec{w}^{T}\left( \ssvec{f}_{v}+\ssvec{f}_{a}\right) \cdot
        \svec{n}\,\mathrm{d}S - \langle D_{v}\rangle - \langle D_{a}\rangle.
\end{equation}
Thus, for both the kinetic energy and thermodynamic entropy we get
\begin{equation}
    \langle\mathcal E_t\rangle +
        \int_{\partial\Omega}\svec{f}^{\mathcal E}\cdot\svec{n}\mathrm{d}S +
        \Theta \mathcal W_{p} + \langle D_{v}\rangle + \langle D_{a}\rangle = 0,\quad
    \svec{f}^{\mathcal E} = \svec{f}^{\mathcal E}_{e} - \stvec{w}^{T}\left(\ssvec{f}_{v} +
        \ssvec{f}_{a}\right),
\end{equation}
with~$\Theta=1$ for~$\mathcal E=\mathcal K$, and~$\Theta=0$ for~$\mathcal E=\mathcal S$, which,
provided that~$D_{v}$ and~$D_{a}$ are positive, is equivalent to~\eqref{eq:NSE:entropy_inequality}
except for the pressure work with kinetic energy entropy variables. We will consider the effect
of~$\mathcal W_p$ in \Sec\ref{sec:stability}, where we perform a semi--discrete entropy stability
analysis.

\section{Numerical approximation}\label{sec:DGSEM}

\subsection{Geometrical transformations}

The domain~$\Omega$ is tessellated into non--overlapping hexahedral elements~$e$ that are
geometrically transformed from a cube~$\rE=[-1,1]^{3}$ called the~\textit{reference element}. To do
that, we define a transfinite mapping from the local
($\svec{\xi}=\left(\xi,\eta,\zeta\right)^T\in \rE$) to the physical ($\svec{x}=(x,y,z)^T\in e$)
coordinates, $\svec{x}=\svec{X}\left(\svec{\xi}\right)$. This mapping also relates the differential
operators from the physical to the computational space. We also compute the covariant and
contravariant basis, and the transformation Jacobian,
\begin{equation}
    \svec{a}_{i} = \frac{\partial\svec{X}}{\partial\xi^{i}},~~
    J\svec{a}^{i}=\svec{a}_{j}\times\svec{a}_{k},~~
    J=\svec{a}_{1}\cdot\left(\svec{a}_{2}\times\svec{a}_{3}\right),~~
    \left(i,j,k\right)\text{cyclical},
\end{equation}
that build the matrix~$\tens{M}=\left(J\svec{a}^{\xi},J\svec{a}^{\eta},J\svec{a}^{\zeta}\right)$.
This matrix is divergence--free,~$\svec{\nabla}_{\xi}\cdot\tens{M}=0$, which is known as the
\textit{metric identities}~\cite{Kopriva2006}. With the matrix~$\tens{M}$, we transform the
gradient and divergence operators~\cite{Kopriva2009},
\begin{equation}
    J\svec{\nabla}u = \tens{M}\svec{\nabla}_{\xi}u,~~J\svec{\nabla}\cdot\svec{f} =
        \svec{\nabla}_{\xi}\cdot\left(\tens{M}^{T}\svec{f}\right) =
        \svec{\nabla}_{\xi}\cdot\svec{\tilde f}.
\end{equation}
Next, we construct a block matrix $\bmat{M}$ from~$\tens{M}$
using~\eqref{eq:notation:block-matrix-from-space-matrix}. By doing so, we extend the gradient and
divergence operators to state and block vectors,
\begin{equation}
    J\svec{\nabla}\stvec{u} = \bmat{M}\svec{\nabla}_{\xi}\stvec{u},~~
    J\svec{\nabla}\cdot\ssvec{f} = \svec{\nabla}_{\xi}\cdot\left(\bmat{M}^{T}\ssvec{f}\right) =
        \svec{\nabla}_{\xi}\cdot\cssvec{f},
\label{eq:dg:geom:transform-operators}
\end{equation}
and we also introduce the matrices~$\tens{M}^T$ and~$\bmat{M}^T$ as transformation matrices
to project covariant space and block vectors into the contravariant frame. We will apply this result
throughout this work to express the fluxes in the reference space.

To write the advection--diffusion system~\eqref{eq:NSE:adv-diff} in the reference space, it is
first casted to a first order system with the definition of the auxiliary variables
$\ssvec{g}=\svec{\nabla}\stvec{w}$ (where the gradient is computed from the entropy
variables~$\stvec{w}$~\eqref{eq:NSE:entropy-variables}). Then these two equations are transformed to
the reference space using~\eqref{eq:dg:geom:transform-operators},
\begin{equation}
\begin{split}
    J\stvec{q}_{t} &+ \svec{\nabla}_{\xi}\cdot\cssvec{f}_{e}\left(\stvec{q}\right) =
        \svec{\nabla}_{\xi}\cdot\left(\cssvec{f}_{v}\left(\stvec{q},\ssvec{g}\right) +
        \cssvec{f}_{a}\left(\stvec{q},\ssvec{g}\right)\right), \\
    J\ssvec{g} &= \bmat M\svec{\nabla}_{\xi}\stvec{w}.
\end{split}
\label{eq:dg:geom:1st-order-ref}
\end{equation}

Finally, we construct two weak forms from~\eqref{eq:dg:geom:1st-order-ref}. We multiply the two
equations by two test functions~$\stvecg{\phi}$~and $\ssvecg{\varphi}$, we integrate the result over
the reference element~$\rE$, and we use the Gauss law on the divergence and gradients,
\begin{equation}
\begin{split}
    \left\langle J\stvec{q}_{t},\stvecg{\phi}\right\rangle_{\rE} &+
        \int_{\partial \rE}\stvecg{\phi}^{T}\cssvec{f}_{e}\cdot\hat{n}\diff\hat{S} -
        \left\langle\cssvec{f}_{e},\svec{\nabla}_{\xi}\stvecg{\phi}\right\rangle_{\rE} \\
    =& \int_{\partial \rE}\stvecg{\phi}^{T}\left(\cssvec{f}_{v} +
        \cssvec{f}_{a}\right)\cdot\hat{n}\diff\hat{S} - \left\langle\cssvec{f}_{v} +
        \cssvec{f}_{a},\svec{\nabla}_{\xi}\stvecg{\phi}\right\rangle_{\rE},\\
    \left\langle J\ssvec{g},\ssvecg{\varphi}\right\rangle_{\rE} =&
        \int_{\partial \rE}\stvec{w}^{T}\cssvecg{\varphi}\cdot\hat{n}\diff\hat{S} -
        \left\langle\stvec{w},\svec{\nabla}_{\xi}\cdot\cssvecg{\varphi}\right\rangle_{\rE}.
\end{split}
\label{eq:dg:geom:weak-forms}
\end{equation}
In~\eqref{eq:dg:geom:weak-forms},~$\hat{n}$ and~$\mathrm{d}\hat{S}$ are the unit normal vector and
the surface differential of the six planar faces of~$\rE$ (e.g. for the
faces~$\eta=\pm 1$,~$\hat{n}=(0,\pm 1,0)^T$ and~$\diff\hat{S}=\diff\xi\diff\zeta$).

\subsection{Polynomial approximation}

In this work we use a high–order Discontinuous Galerkin Spectral Element Method
(DGSEM)~\cite{Kopriva2009,Winters2021}. In particular, this work uses the Gauss-Lobatto version of
the DGSEM, which makes it possible to construct entropy-stable schemes using the summation-by-parts
simultaneous-approximation-term (SBP-SAT) property and two-point entropy--conservative fluxes.
Moreover, it handles arbitrary three dimensional curvilinear hexahedral meshes while maintaining
high order, spectral accuracy and entropy--stability.
\begin{equation}
    \stvec{q}\bigr|_{\rE} \approx \stvec{Q}\left(\svec{\xi},t\right) =
        \sum_{i,j,k=0}^{N}\stvec{Q}_{ijk}(t)l_i(\xi)l_j(\eta)l_k(\zeta).
\label{eq:dg:pol:Q}
\end{equation}
In~\eqref{eq:dg:pol:Q}, and following the notation of~\cite{Kopriva2017},~$l_i\left(\xi\right)$ are
the Lagrange polynomials associated to the Gauss--Lobatto
points,~$\left\{\xi_i,\eta_i,\zeta_i\right\}_{i,j,k=0}^{N}$. The coefficients
$\stvec{Q}_{ijk}\left(t\right)$ are called the nodal degrees of freedom, whose values coincide with
the interpolant evaluation at the GL points,~$\stvec{Q}_{ijk} =
\stvec{Q}\left(\xi_i,\eta_j,\zeta_k,t\right)$. Similarly, the fluxes are also represented by their
polynomial approximation~\cite{Kopriva2009},
\begin{equation}
    \ssvec{f}\bigr|_{\rE}\approx \ssvec{F}\left(\svec{\xi},t\right) =
        \sum_{i,j,k=0}^{N}\ssvec{F}_{ijk}(t)l_i(\xi)l_j(\eta)l_k(\zeta),~~
    \ssvec{F}_{ijk} = \ssvec{f}\left(\stvec{Q}_{ijk}\right).
\label{eq:dg:pol:F}
\end{equation}

The rest of the quantities involved are approximated following~\eqref{eq:dg:pol:Q}
and~\eqref{eq:dg:pol:F}, with the exception of the contravariant basis~$J\svec{a}^{i}$ that build
the matrix~$\mathcal{M}$. In the continuous setting, the metric terms satisfy the metric identities,
$\svec{\nabla}_{\xi}\cdot\tens{M}=0$. In its discrete counterpart, the discrete metric identities
represent the ability of the scheme to be free--stream preserving~\cite{Kopriva2006}, and it is also
a requirement for discrete stability~\cite{Gassner2018}. Thus, we follow~\cite{Kopriva2006} and use
a conservative form of the metrics that ensure the preservation of the metric identities discretely.

We introduce the polynomial approximation in~\eqref{eq:dg:geom:weak-forms},
\begin{subequations}\label{eq:dg:pol:weak-forms}
\begin{align}
\begin{split}
    \left\langle\mathcal J\stvec{Q}_{t},\stvecg{\phi}\right\rangle_{\rE} &+
        \int_{\partial \rE}\stvecg{\phi}^{T}\cssvec{F}_{e}\cdot\hat{n}\diff\hat{S} -
        \left\langle\cssvec{F}_{e},\svec{\nabla}_{\xi}\stvecg{\phi}\right\rangle_{\rE}\\
    =& \int_{\partial \rE}\stvecg{\phi}^{T}\left(\cssvec{F}_{v} +
        \cssvec{F}_{a}\right)\cdot\hat{n}\diff\hat{S} -
        \left\langle\left( \cssvec{F}_{v}+\cssvec{F}_{a}\right),
        \svec{\nabla}_{\xi}\stvecg{\phi}\right\rangle_{\rE},
        \label{eq:dg:pol:weak-forms:Q}
\end{split}\\
    \left\langle \mathcal J\ssvec{G},\ssvecg{\varphi}\right\rangle_{\rE} =&
        \int_{\partial \rE}\stvec{W}^{T}\cssvecg{\varphi}\cdot\hat{n}\diff\hat{S} -
        \left\langle\stvec{W},\svec{\nabla}_{\xi}\cdot\cssvecg{\varphi}\right\rangle_{\rE},
        \label{eq:dg:pol:weak-forms:G}
\end{align}
\end{subequations}
where the test functions are also restricted to~$N$ order polynomials in the reference space.

Next, we approximate the exact integrals by quadrature rules. In the reference element
${\rE=[-1,1]^{3}}$, the integrals are computed as the product of the nodal degrees of freedom by the
quadrature weights $w_{i}=\int_{-1}^{1}l_i\left(\xi\right)\diff \xi$~\cite{Kopriva2009},
\begin{equation}
    \left\langle f \right\rangle_{\rE} \approx \left\langle f\right\rangle_{\rE,N} =
        \sum_{i,j,k=0}^{N}w_{i}w_{j}w_{k}F_{ijk} = \sum_{i,j,k=0}^{N}w_{ijk}F_{ijk}.
\end{equation}
The result is exact if~$f\in\mathcal P^{2N-1}$. Then, the surface integrals at the
boundary~$\partial \rE$ are also approximated with quadratures. We take into account that the normal
vectors at the surfaces $\xi=\pm 1$, $\eta=\pm 1$ and $\zeta=\pm 1$ that make the reference element
faces are $\hat{n}_1=\pm (1,0,0)$, $\hat{n}_2=\pm (0,1,0)$ and $\hat{n}_3=\pm (0,0,1)$ to simplify
the expressions,
\begin{equation}
    \int_{\partial \rE}\svec{\tilde{f}}\cdot\hat{n}\diff\hat{S} \approx
        \int_{\partial \rE,N}\svec{\tilde{f}}\cdot\hat{n}\diff\hat{S} =
        \sum_{i,j=0}^{N}w_iw_j\left(\tilde{F}^{\xi}_{Nij} - \tilde{F}^{\xi}_{0ij} +
        \tilde{F}^{\eta}_{iNj} - \tilde{F}^{\eta}_{i0j} + \tilde{F}^{\zeta}_{ijN} -
        \tilde{F}^{\zeta}_{ij0}\right).
\end{equation}
With Gauss--Lobatto points no additional projection is required and one only chooses the appropriate
nodal degree of freedom that corresponds to a given surface node
(e.g.~$\stvec{Q}(1,\eta_j,\zeta_k)=\stvec{Q}_{Njk}$). The latter, alongside the exactness of the
quadrature, makes the scheme satisfy the discrete Gauss law~\cite{gassner2013skew},
\begin{equation}
    \left\langle \svec{\tilde F},\svec{\nabla}_{\xi}G\right\rangle_{\rE,N} =
        \int_{\partial \rE,N}G\svec{\tilde{F}}\cdot\hat{n}\diff \hat{S} -
        \left\langle\svec{\nabla}_{\xi}\cdot\svec{\tilde{F}},G\right\rangle_{\rE,N}.
\label{eq:dg:pol:disc-Gauss-law}
\end{equation}
This is a direct consequence of the aforementioned SBP derivative
operator~\cite{Carpenter1999,Nordstrom2013,Fisher2013,DelRey2014} that, when expressed in matrix
form,~$D$, has the properties,
\begin{equation}
    PD = Q, \quad P = P^T, \quad Q+Q^T=B, \quad B = e_N e_N^T - e_0 e_0^T,
\label{eq:dg:pol:sbp}
\end{equation}
where~$P$ is a diagonal matrix with the quadrature weights and~$e_i$ are column vectors
filled with zeros except for position~$i$. Equation~\eqref{eq:dg:pol:sbp} is simply the matrix form
of~\eqref{eq:dg:pol:disc-Gauss-law} for SBP operators.

As a final remark, it is possible to write surface integrals both in computational or physical
coordinates. The matrix $\bmat{M}$ relates the normal unit vector and surface differentials,
$\svec{n}\diff S = \mathcal{M}\cdot\hat{n}\diff\hat{S}$. Hence, we can write,
\begin{equation}
    \int_{\partial \rE,N}\svec{\tilde{F}}\cdot\hat{n}\diff\hat{S} =
        \int_{\partial e,N}\svec{{F}}\cdot\svec{n}\diff{S}.
\label{eq:dg:pol:int-f-de-dE}
\end{equation}
The representation in physical coordinates is convenient because at the interface between two
elements,~$e^{+}$ and~$e^{-}$, the normal vectors satisfy
${\svec{n}^{+}\mathrm{d}S^{+}=-\svec{n}^{-}\mathrm{d}S^{-}}$.

We approximate the integrals in~\eqref{eq:dg:pol:weak-forms} with quadratures, and
use~\eqref{eq:dg:pol:int-f-de-dE} to get
\begin{subequations}\label{eq:dg:pol:weak-forms-quad}
\begin{align}
\begin{split}
    \left\langle\mathcal J\stvec{Q}_{t},\stvecg{\phi}\right\rangle_{\rE,N} &+
        \int_{\partial e,N}\stvecg{\phi}^{T}\ssvec{F}_{e}^{\star}\cdot\svec{n}\diff{S} -
        \left\langle\cssvec{F}_{e},\svec{\nabla}_{\xi}\stvecg{\phi}\right\rangle_{\rE,N}\\
    =& \int_{\partial e,N}\stvecg{\phi}^{T}\left(\ssvec{F}_{v}^{\star} +
        \ssvec{F}_{a}^{\star}\right)\cdot\svec{n}\diff{S} - \left\langle\cssvec{F}_{v} +
        \cssvec{F}_{a},\svec{\nabla}_{\xi}\stvecg{\phi}\right\rangle_{\rE,N},
        \label{eq:dg:pol:weak-forms-quad:Q}
\end{split}\\
    \left\langle \mathcal J\ssvec{G},\ssvecg{\varphi}\right\rangle_{\rE,N} =&
        \int_{\partial e,N}\stvec{W}^{\star,T}\ssvecg{\varphi}\cdot\svec{n}\diff{S} -
        \left\langle\stvec{W},\svec{\nabla}_{\xi}\cdot\cssvecg{\varphi}\right\rangle_{\rE,N}.
        \label{eq:dg:pol:weak-forms-quad:G}
\end{align}
\end{subequations}
In~\eqref{eq:dg:pol:weak-forms-quad}, we have replaced the interface fluxes by the uniquely defined
\textit{numerical fluxes} that depend on the two adjacent states~\cite{Kopriva2009}. These are
detailed in \Sec\ref{sec:dg:f-star}.

Finally, we use a split--form scheme for the inviscid fluxes. The split--form method allows us to
obtain schemes that are entropy stable discretely through the use of a two--point volume flux
function~\cite{Fisher2013}. Thus, we apply the discrete Gauss law on the inviscid fluxes
of~\eqref{eq:dg:pol:weak-forms-quad:Q},
\begin{equation}
\begin{split}
    \left\langle\mathcal J\stvec{Q}_{t},\stvecg{\phi}\right\rangle_{\rE,N} &+
        \int_{\partial e,N}\stvecg{\phi}^{T}\left(\ssvec{F}_{e}^{\star} -
        \ssvec{F}_{e}\right)\cdot\svec{n}\diff{S} +
        \left\langle\svec{\nabla}_{\xi}\cdot\cssvec{F}_{e},\stvecg{\phi}\right\rangle_{\rE,N} \\
    =& \int_{\partial e,N}\stvecg{\phi}^{T}\left(\ssvec{F}_{v}^{\star} +
        \ssvec{F}_{a}^{\star}\right)\cdot\svec{n}\diff{S} - \left\langle\cssvec{F}_{v} +
        \cssvec{F}_{a},\svec{\nabla}_{\xi}\stvecg{\phi}\right\rangle_{\rE,N},
\end{split}
\label{eq:dg:pol:weak-forms-sbp}
\end{equation}
and compute the divergence with a two--point volume flux function
$\ssvec{F}^{\#}\left(\stvec{Q}_{ijk},\stvec{Q}_{nml}\right)$,
\begin{equation}
    \left(\svec{\nabla}_{\xi}\cdot\cssvec{F}_{e}\right)_{ijk} \!\!\!\!\approx
        \mathbb D\left(\ssvec{F}^{\#}_{e}\right)_{ijk} \!\!\!\! = 2\sum_{n=0}^{N}
        \left(D_{in}\cssvec{F}^{\#,1}\!\!\!\!\!\!\!\!\left(\stvec{Q}_{ijk},\stvec{Q}_{njk}\right) +
        D_{jn}\cssvec{F}^{\#,2}\!\!\!\!\!\!\!\!\left(\stvec{Q}_{ijk},\stvec{Q}_{ink}\right) +
        D_{kn}\cssvec{F}^{\#,3}\!\!\!\!\!\!\!\!\left(\stvec{Q}_{ijk},\stvec{Q}_{ijn}\right)\right).
\label{eq:dg:pol:two-point-divergence}
\end{equation}
where~$D_{ij} = l'_j(\xi_i)$. The interested reader can find the details on~\cite{gassner2016split}.
Throughout this work, two different two--point volume flux functions will be used. For kinetic
energy preserving schemes, we use the two--point flux by Pirozzoli~\cite{Pirozzoli2010}, whereas for
thermodynamic entropy stable schemes we use the two--point flux by
Chandrashekar~\cite{Chandrashekar2013}. With the two--point divergence, the DG
scheme~\eqref{eq:dg:pol:weak-forms-sbp} is
\begin{subequations}\label{eq:dg:pol:NSE-split-form}
\begin{align}
\begin{split}
    \left\langle\mathcal J\stvec{Q}_{t},\stvecg{\phi}\right\rangle_{\rE,N} &+
        \int_{\partial e,N}\stvecg{\phi}^{T}\left(\ssvec{F}_{e}^{\star} -
        \ssvec{F}_{e}\right)\cdot\svec{n}\diff{S} +
        \left\langle\mathbb{D}\left(\cssvec{F}_{e}^{\#}\right),\stvecg{\phi}\right\rangle_{\rE,N} \\
    =& \int_{\partial e,N}\stvecg{\phi}^{T}\left(\ssvec{F}_{v}^{\star} +
        \ssvec{F}_{a}^{\star}\right)\cdot\svec{n}\diff{S} - \left\langle\cssvec{F}_{v} +
        \cssvec{F}_{a},\svec{\nabla}_{\xi}\stvecg{\phi}\right\rangle_{\rE,N},
        \label{eq:dg:pol:NSE-split-form:Q}
\end{split}\\
    \left\langle \mathcal J\ssvec{G},\ssvecg{\varphi}\right\rangle_{\rE,N} =&
        \int_{\partial e,N}\stvec{W}^{\star,T}\ssvecg{\varphi}\cdot\svec{n}\diff{S} -
        \left\langle\stvec{W},\svec{\nabla}_{\xi}\cdot\cssvecg{\varphi}\right\rangle_{\rE,N}.
        \label{eq:dg:pol:NSE-split-form:G}
\end{align}
\end{subequations}

Recall at this point that the flux~$\ssvec{F}_{a}$ represents the additional viscosity added by
either the molecular viscosity~$\ssvec{F}_{a}=\ssvec{F}_{v}(\mu_{a})$ or the Guermond--Popov
fluxes~$\ssvec{F}_{a}=\ssvec{F}_{\GP}(\alpha_{a},\mu_{a})$.

\subsection{Intercell numerical fluxes}\label{sec:dg:f-star}

One of the characteristics that make DG schemes attractive is that they can introduce dissipation at
the inter--element faces. The characteristics of this dissipation are dictated by the choice of
numerical fluxes~$\ssvec{F}^{\star}_e$,~$\ssvec{F}^{\star}_v$ and~$\ssvec{F}^{\star}_a$. We write
the one--dimensional numerical fluxes, which are later transformed to three--dimensions using the
rotational invariance property~\cite{Toro2009},
\begin{equation}
    \stvec{F}_{e}^{\star}\cdot\svec{n} =
    \smat{\mathrm T}^{T}\stvec{F}^{1,\star}_{e}\left(\stvec{Q}_{nL},\stvec{Q}_{nR}\right),\quad
    \stvec{Q}_{n} = \smat{\mathrm T}\stvec{Q},\quad
    \smat{\mathrm T}=\left(\begin{array}{ccccc}
        1 & 0 & 0 & 0 & 0 \\
        0 & n_{x} & n_{y} & n_{z} & 0 \\
        0 & t_{1,x} & t_{1,y} & t_{1,z} & 0 \\
        0 & t_{2,x} & t_{2,y} & t_{2,z} & 0 \\
        0 & 0 & 0 & 0 & 1 \end{array}\right),
\end{equation}
where~$\svec{n}$,~$\svec{t_1}$ and~$\svec{t_2}$ are the three vectors that define the local
reference frame at a point of the interface. For the inviscid fluxes, one can construct entropy
conserving schemes (i.e. without numerical dissipation) using the two--point flux as a numerical
flux,
\begin{equation}
    \ssvec{F}^{\star}_{e}(\stvec{Q}_L, \stvec{Q}_R)\cdot\svec{n} =
        \smat{\mathrm T}^{T}\ssvec{F}^{1,\#}_{e}(\stvec{Q}_{nL}, \stvec{Q}_{nR}),
\label{eq:dg:numerical-fluxes:EC}
\end{equation}
or construct an entropy stable (dissipative) approximation by adding extra terms to it. For the
latter, there are several choices available, which depend on whether one constructs a kinetic energy
preserving or a thermodynamic entropy stable scheme~\cite{Chandrashekar2013}. For kinetic energy
consistency (e.g. Pirozzoli's), the Lax--Friedrichs numerical flux is stable~\cite{gassner2016split}
\begin{equation}
    \stvec{F}^{\star}_{e}\left(\stvec{Q}_{nL},\stvec{Q}_{nR}\right)\cdot\svec{n} =
        \smat{\mathrm T}^{T}\left(\stvec{F}^{1,\#}_{e}\left(\stvec{Q}_{nL},\stvec{Q}_{nR}\right) -
        \frac{1}{2}|\lambda_{\max}|\jump{\stvec{Q}_n}\right),~~
    \jump{\stvec{Q}} = \stvec{Q}_{R}-\stvec{Q}_{L},
\label{eq:dg:f-star:LxF}
\end{equation}
and for the thermodynamic entropy (e.g. Chandrashekar's), a popular choice is to construct a matrix
dissipation flux based on Roe's Riemann solver~\cite{Chandrashekar2013},
\begin{equation}
    \stvec{F}^{\star}_{e}\left(\stvec{Q}_{nL},\stvec{Q}_{nR}\right)\cdot\svec{n} =
        \smat{\mathrm T}^{T}\left(\stvec{F}^{1,\#}_{e}\left(\stvec{Q}_{nL},\stvec{Q}_{nR}\right) -
        \frac{1}{2}\smat{\mathrm M}\jump{\stvec{W}_{n}}\right),
\label{eq:dg:f-star:mat-diss}
\end{equation}
where the positive definite matrix~$\smat{\mathrm M}$ depends on the two states and is constructed
from the eigenvectors and the absolute value of the eigenvalues of the Euler equations.

For the entropy variables, viscous fluxes, and artificial viscosity, we follow~\cite{Gassner2018}
and use the Bassi--Rebay~1 (BR1) scheme~\cite{BR1},
\begin{equation}
    \stvec{W}^{\star} = \aver{\stvec{W}},~~
    \ssvec{F}_{v}^{\star} = \aver{\ssvec{F}_{v}},~~
    \ssvec{F}_{a}^{\star} = \aver{\ssvec{F}_{a}}.
\label{eq:dg:BR1}
\end{equation}

\subsection{Physical boundaries numerical fluxes}\label{sec:dg:f-star-bound}

The boundary conditions are weakly enforced through numerical fluxes applied to the interior state
and to a ghost state constructed from the boundary data. In this work, we use free-- and no--slip
adiabatic walls, inflow, and outflow boundary conditions. The implementation of the adiabatic wall
boundary conditions is entropy stable, while inflow and outflow boundaries can increase and/or
decrease the entropy.

\subsubsection{Adiabatic wall}

For the inviscid fluxes, both free-- and no--slip walls cancel the normal velocity. This enforcement
is done through the numerical flux, after the construction of a mirrored ghost state with
negative normal velocity,
\begin{equation}
    \stvec{Q}_{n,\mathrm{ghost}} = \left(\rho, -\rho U, \rho V, \rho W, \rho E\right)^T,\quad
    \ssvec{F}_e^{\star}\cdot\svec{n} = \ssvec{F}_e^{\star}\left(\stvec{Q}_n,
        \stvec{Q}_{n,\mathrm{ghost}}\right)\cdot\svec{n}
\label{eq:dg:f-star-bound:Q-ghost}
\end{equation}

For the viscous fluxes, the free--slip adiabatic wall cancels the stress tensor and the heat flux at
the wall, thus we enforce Neumann boundary conditions in the state variables,
\begin{equation}
    \stvec{W}_{n}^{\star} = \stvec{W}_{n}, \quad \ssvec{F}_{v}^{\star}\cdot\svec{n}=0.
\label{eq:dg:f-bound:free-slip-visc}
\end{equation}

No--slip adiabatic walls, on the contrary, are of Dirichlet type for the three velocities and
Neumann for the energy,
\begin{equation}
    \stvec{W}_{n}^{\star} = \left(W_{1n},0,0,0,W_{5n}\right)^T, \quad
    \ssvec{F}_{v}^{\star}\cdot\svec{n} = \left(0,F_{v,2},F_{v,3},F_{v,4},0\right)^T.
\label{eq:dg:f-bound:no-slip-visc}
\end{equation}
For the Guermond--Popov fluxes, we cancel the flux at the physical boundary for the first and
fifth equations, and use either Dirichlet boundary conditions at no--slip walls,
($F^{\star}=F, W^{\star}=0$) or Neumann boundary conditions at free--slip walls
($F^{\star}=0, W^{\star}=W$) for the momentum equations,
\begin{equation}\label{eq:dg: physical-bc:art-fluxes}
    \stvec{W}_{n}^{\star} =\left(W_{1n},W^{\star}_{2n},W^{\star}_{3n},
        W^{\star}_{4n},W_{5n}\right)^T, \quad
    \ssvec{F}_{\GP}^{\star}\cdot\svec{n} =
        \left(0,F^{\star}_{\GP,2},F^{\star}_{\GP,3},F^{\star}_{\GP,4},0\right)^T.
\end{equation}

\subsubsection{Supersonic inflow}

Inflow boundary conditions at supersonic velocities only carry information from the exterior into
the domain and thus, they can be imposed in the inviscid term by setting a ghost state with the
boundary values and using inter--element fluxes to couple it with the interior solution,
\begin{equation}
    \stvec{Q}_{n,\mathrm{ghost}} =
        \left(\rho_0,\rho_0U_0,\rho_0V_0,\rho_0W_0,\rho_0E_0\right)^T.
\end{equation}
For the viscous fluxes we use Neumann boundary conditions and cancel the viscous stress as done for
free--slip walls.

\subsubsection{Outflow}

Since information travels upstream only in subsonic flows, the value of the
pressure at the boundary for the inviscid flux is imposed through ghost states only for this case,
whereas it takes the value of the interior if the local Mach number is greater than one,
\begin{equation}
    \begin{array}{ll}
        \stvec{Q}_{n,\mathrm{ghost}} =
            \left(\rho, \rho U, \rho V, \rho W, \rho E\right)^T, & \quad \mathrm{if}\,M > 1, \\
        \stvec{Q}_{n,\mathrm{ghost}} = \left(\rho_0, \rho_0 U_0, \rho_0 V_0, \rho_0 W_0,
            \rho_0 E_0\right)^T, & \quad \mathrm{if}\,M \leqslant 1,
    \end{array}
\end{equation}
where we compute the exterior values from the Riemann invariants,
\begin{equation}
\begin{gathered}
    \rho_0 = \rho \left( 1 + \frac{p_0/p-1}{\gamma} \right), \\
    \svec{U}_0 = \svec{U}_t + \svec{U}_{0,n}, \\
        \svec{U}_{0,n} = r^+ - \frac{2a_0}{\gamma-1},~~ r^+ = \svec{U}_n + \frac{2a}{\gamma-1},
\end{gathered}
\end{equation}
being~$\svec{U}_n$ and~$\svec{U}_t$ the normal and tangent velocity vectors. As for the viscous
fluxes, we again use Neumann boundary conditions with no viscous stress, and the entropy variables
at the outlet use the values of the interior points. For more details,
see~\cite{Carlson2011,Mengaldo2014}.

\section{Spectral vanishing viscosity: filtered artificial dissipation}\label{sec:DG:SVV}

In this section we introduce the construction of an entropy stable spectral vanishing viscosity that
helps us to modulate the dissipation by spatially filtering the artificial viscosity fluxes. The SVV
approach was first considered by Tadmor~\cite{Tadmor1989}, where he added a filtered viscous term to
the Burgers' equation,
\begin{equation}
    \frac{\partial U}{\partial t} + \frac{1}{2}\frac{\partial U^2}{\partial x} =
        \varepsilon\frac{\partial}{\partial x}\left[\mathcal{F} \star
        \frac{\partial U}{\partial x}\right], \quad
    \varepsilon \sim \frac{1}{N^{\beta}\log N}, \quad \beta < 1.
\end{equation}
\revtwo{In this work we use Legendre polynomials,~$L_k(x)$, as a modal basis, where~$k$ not only
indicates the mode number but also serves as a measure of the wavenumber. In this case, the
filtering operator,~$\star$, can be expressed as an element-wise product where, for each
mode,~$L_k$, there is a filter kernel coefficient,~$\mathcal{F}_k$, that modulates its intensity.}
Under this assumptions, \revtwo{Tadmor} proved that the solution converges to the weak solution that
is also physically plausible. However, there is not only one possible shape for this filter kernel
and different \revtwo{options} have been proposed. For instance, in~\cite{Tadmor1989} it is defined
as,
\begin{equation}
    \hat{\mathcal{F}}_k = \left\{\begin{array}{ll}
            0, & k \leqslant M \\
            1, & k > M, \\
        \end{array}\right. \quad \varepsilon M = 0.25,
\end{equation}
or, in~\cite{Maday1993,Kaber1996}, a~$C^{\infty}$ function is used to improve the resolution,
\begin{equation}
    \hat{\mathcal{F}}_k = \left\{\begin{array}{ll}
            0, & k \leqslant M \\
            \exp\left[-\dfrac{(k-N)^2}{(k-M)^2}\right], & k > M, \\
        \end{array}\right. \quad M \approx 5\sqrt{N}.
\end{equation}
Finally, taking advantage of the fact that the viscosity varies for the different frequencies
considered by the solution approximation, we can also use it to add more dissipation in turbulent,
unresolved regions, where the one from the underlying numerical scheme is not enough. This approach
was followed in~\cite{Karamanos2000,Kirby2006} and it can also be combined with other LES methods,
as in~\cite{Kirby2002,Pasquetti2008}.

We first establish the foundation in one dimension, and then extend the methodology to the
three--dimensional curvilinear setting. We also define the reference line~$\rL=[-1,1]$,
face~$\rF=[-1,1]^2$ and element~$\rE=[-1,1]^3$.

\subsection{Polynomial filtering: one dimension}\label{sec:DG:filtering:1D}

In the nodal DG method, the solution approximated by degree~$N$ polynomials is represented using the
Lagrange polynomials,
\begin{equation}
    Q(\xi) = \sum_{j=0}^{N} Q_j l_j(\xi),
\end{equation}
but other representations are also possible. Precisely, we are interested in a modal reconstruction
since it provides a natural approach to polynomial filtering, as we can directly manipulate the
intensity of each polynomial mode. We introduce the Legendre
polynomials,~$L_j(\xi)$, to represent the solution,
\begin{equation}
    Q(\xi) = \sum_{j=0}^{N} Q_j l_j(\xi)=\sum_{j=0}^{N}\hat{Q}_{j}L_{j}\left(\xi\right),
\label{eq:DG:filtering:Q-nodal-modal}
\end{equation}
where~$\hat{Q}_j$ are the modal coefficients of the solution (as opposed to the nodal
coefficients,~$Q_j$).

The Legendre polynomials are an orthogonal basis (both continuosly and discretely) since they
satisfy,
\begin{equation}
    \left\langle L_n,L_m\right\rangle_{\rL,N} = \left\Vert L_n\right\Vert_{N}^{2}\delta_{nm}.
\end{equation}
Note that all the norms $\left\Vert L_n\right\Vert_{N}$ are exactly computed with quadratures (i.e.
they involve a quadrature of degree less than $2N-1$) except for the one associated to the highest
mode, $\left\Vert L_{N}\right\Vert_{N}$. The exactness of the norm is not essential however, as long
as it is positive.

We define a \textit{forward} (F) operation that computes the modal coefficients from the nodal
quantities. To do that, we multiply~\eqref{eq:DG:filtering:Q-nodal-modal} by one of the Legendre
polynomials~$L_i\left(\xi\right)$ and compute its quadrature over~$\rL$,
\begin{equation}
    \left\langle Q,L_i\right\rangle_{\rL,N} =
        \sum_{j=0}^{N}Q_{j}\left\langle l_j,L_i\right\rangle_{\rL,N} =
        \sum_{j=0}^{N}\hat{Q}_{j}\left\langle L_j,L_i\right\rangle_{\rL,N} =
        \left\Vert L_{i}\right\Vert_N^{2} \hat{Q}_{i},
\end{equation}
therefore,
\begin{equation}
    \hat{Q}_i = \sum_{j=0}^{N}\frac{\left\langle l_j,L_i\right\rangle}
        {\left\Vert L_i\right\Vert^{2}}Q_{j} = \sum_{j=0}^{N}F_{ij}Q_{j},\quad
    F_{ij} = \frac{\left\langle l_j,L_i\right\rangle}{\left\Vert L_i\right\Vert^{2}}.
\end{equation}

Similarly, the \textit{backward} (B) operation recovers the nodal coefficients from the modal form.
In this case, we multiply~\eqref{eq:DG:filtering:Q-nodal-modal} by one of the Lagrange
polynomials,~$l_{i}\left(\xi\right)$ and compute its quadrature over $\rL$,
\begin{equation}
    \left\langle Q,l_i\right\rangle_{\rL,N} =
        \sum_{j=0}^{N}Q_{j}\left\langle l_j,l_i\right\rangle_{\rL,N} = w_iQ_i =
        \sum_{j=0}^N\hat{Q}_{j}\left\langle L_j,l_i\right\rangle_{\rL,N},
\end{equation}
and therefore,
\begin{equation}
    Q_{i} = \sum_{j=0}^{N}\frac{\left\langle L_j,l_i\right\rangle_{\rL,N}}{w_{i}}\hat{Q}_{j} =
        \sum_{j=0}^{N}B_{ij}\hat{Q}_{j},\quad
    B_{ij} = \frac{\left\langle L_j,l_i\right\rangle_{\rL,N}}{w_{i}}.
\end{equation}

This two transformations and its associated matrices,~$F$ and~$B$, satisfy two properties required
for the stability of the scheme.
\begin{prop}\label{prop:1}
They are inverse matrices, $FB=BF=I$, as the direct multiplication shows,
\begin{equation}
\begin{split}
    (FB)_{ij} &= \sum_{k=0}^{N}F_{ik}B_{kj} =
        \sum_{k=0}^{N}\frac{\left\langle L_i,l_k\right\rangle_{\rL,N}}
            {\left\Vert L_i\right\Vert_{N}^{2}}
            \frac{\left\langle l_k,L_j\right\rangle_{\rL,N}}{w_k}\\
    &= \sum_{k=0}^{N}\frac{1}{w_{k}\left\Vert L_i\right\Vert^2_{N}}
        \sum_{l=0}^{N} w_{l}^2 L_i(\xi_l)L_j(\xi_l)l_k^2(\xi_l)\\
    &= \frac{1}{\left\Vert L_i\right\Vert^2_{N}}
        \sum_{k=0}^{N}w_{k}L_i(\xi_{k})L_j(\xi_{k}) =
        \frac{\left\langle L_i,L_j\right\rangle_{\rL,N}}{\Vert L_i\Vert^2_{N}} = \delta_{ij}
\end{split}
\end{equation}
\end{prop}
\begin{prop}\label{prop:2}
Multiplied by the discrete nodal and modal weights, they are transposed matrices,
\begin{equation}
    w_iB_{ij} = \Vert L_j\Vert_{N}^{2}F_{ji}.
\label{eq:DG:filtering:prop2}
\end{equation}
\end{prop}

Properties~\ref{prop:1} and~\ref{prop:2} allow us to relate the discrete norms in both nodal and
modal spaces (Parseval's identity),
\begin{prop}\label{prop:3}
Discrete Parseval's identity: discrete inner products can be equivalently computed from nodal or
modal coefficients
\begin{equation}
    \left\langle U,V\right\rangle_{\rL,N} = \sum_{i=0}^{N}w_i U_i V_i =
        \sum_{j=0}^{N}\Vert L_j\Vert_{N}^{2} \hat{U}_{j}\hat{V}_{j}.
\end{equation}
The proof for this property uses Properties~\ref{prop:1} (at the second equality) and~\ref{prop:2}
(at the fourth equality),
\begin{equation}
\begin{split}
    \left\langle U,V\right\rangle_{\rL,N} &= \sum_{i=0}^{N}w_i U_i V_i =
        \sum_{i=0}^{N}w_i U_i \sum_{j,k=0}^{N}B_{ij}F_{jk}V_{k} =
        \sum_{i,j=0}^{N}U_i \left(w_i B_{ij}\right)\sum_{k=0}^{N}F_{jk}V_{k} \\
    &= \sum_{i,j=0}^{N}U_i \Vert L_j\Vert_{N}^{2}F_{ji}\hat{V}_{j} =
        \sum_{j=0}^{N}\Vert L_j\Vert_{N}^{2} \hat{V}_{j}\sum_{i=0}^{N}F_{ji}U_i =
        \sum_{j=0}^{N}\Vert L_j\Vert_{N}^{2} \hat{U}_{j}\hat{V}_{j}.
\end{split}
\label{eq:DG:filtering:discrete-norm}
\end{equation}
The ability to compute equivalent discrete L$^2$ norms either from the nodal or modal spaces is key
to show the discrete stability of filtered dissipative fluxes.
\end{prop}

In this work, we locally filter the polynomial functions by applying a given filter kernel according
to,
\begin{equation}
    \mathcal F\left(\xi\right) = \sum_{j=0}^{N}\hat{\mathcal F}_{j}L_{j}\left(\xi\right).
\end{equation}
In the case of the solution vector,~$Q$, the application of the filter kernel returns another
polynomial whose modal coefficients are the point-wise product of the modal coefficients
of~$\mathcal{F}$ and~$Q$,
\begin{equation}
    \bar{Q} = \mathcal F{\star}Q =
        \sum_{j=0}^{N} \hat{\mathcal F}_{j}\hat{Q}_{j}L_j\left(\xi\right).
\label{eq:DG:filtering:modal}
\end{equation}
Thus we get the filtered solution~$\bar{Q}$ in the modal space, which is then returned to a nodal
representation,~$\bar{Q}_{i}$, afterwards. For computational efficiency, we can define a single
matrix--vector product operation that encapsulates the forward, filtering, and backwards steps,
\begin{equation}
    \bar{Q}_{i} = \sum_{j=0}^{N}\mathcal H_{ij}Q_{j},\quad
    \mathcal H = B\mathrm{diag}\left(\hat{\mathcal{F}}\right)F.
\end{equation}

The next theorem, also used similarly in~\cite{Nordstrom2021}, will prove to be fundamental for
entropy stability,
\begin{theorem}\label{theorem:1}
The inner product of~$Q$ and~$\bar{Q}$ is positive (assuming positive filtering
coefficients~$\hat{\mathcal F}_{i}$),
\begin{equation}
    \left\langle Q,\bar{Q}\right\rangle_{\rL,N} =
        \left\langle Q, \mathcal F{\star}Q\right\rangle_{\rL,N}\geqslant 0.
\label{eq:DG:filtering:theorem1}
\end{equation}
\end{theorem}
\begin{proof}
Using Property~\ref{prop:3}, the inner product of~$Q$ and~$\bar{Q}$ can be expressed in modal space
where, according to~\eqref{eq:DG:filtering:modal},
\begin{equation}
    \bar{Q} = \sum_{i=0}^N \hat{\bar{Q}}_i L_i, \quad \hat{\bar{Q}}_i =
        \hat{\mathcal F}_i \hat{Q}_i.
\end{equation}
Introducing this definition into the inner product in~\eqref{eq:DG:filtering:theorem1},
\begin{equation}
    \left\langle Q,\bar{Q}\right\rangle_{\rL,N} = \sum_{j=0}^{N}\Vert L_{j}\Vert^2
        \hat{Q}_{j}\hat{\bar{Q}}_{j} = \sum_{j=0}^{N}\Vert L_{j}\Vert^2
        \hat{\mathcal F}_{j}\hat{Q}_{j}^2,
\end{equation}
which is positive as long as the filter kernel coefficients are positive. \qed
\end{proof}

\subsection{Polynomial filtering: three dimensions}\label{sec:DG:filtering:3D}

In three dimensions, the polynomial framework is constructed from a tensor product form of the
one--dimensional one. Therefore, it automatically inherits the properties described for the
one--dimensional space in \Sec\ref{sec:DG:filtering:1D}.

We get the modal representation of a solution,
\begin{equation}
    \stvec{Q} =
        \sum_{i,j,k=0}^{N}\stvec{Q}_{ijk}l_i(\xi)l_j(\eta)l_k(\zeta) =
        \sum_ {i,j,k}^{N}\hat{\stvec{Q}}_{ijk}L_{i}(\xi)L_{j}(\eta)L_{k}(\zeta),
\end{equation}
from a weak--form with the tensor product Legendre polynomials
$L_{ijk}=L_i(\xi)L_j(\eta)L_k(\zeta)$,
\begin{equation}
    \left\langle \stvec{Q}, L_{ijk}\right\rangle_{\rE,N} =
        \sum_{m,n,l=0}^{N}\stvec{Q}_{mnl}\left\langle l_{mnl},  L_{ijk}\right\rangle_{\rE,N} =
        \sum_{m,n,l=0}^{N}\hat{\stvec{Q}}_{mnl}\left\langle L_{mnl}, L_{ijk}\right\rangle_{\rE,N}.
\end{equation}

The tensor product inner products in the three--dimensional reference element~$\rE$ reduce to the
product of the inner products in the reference line~$\rL$,
\begin{equation}
    \left\langle l_{mnl}, L_{ijk}\right\rangle_{\rE,N} =
        \left\langle l_m,L_i\right\rangle_{\rL,N} \left\langle l_n,L_j\right\rangle_{\rL,N}
        \left\langle l_l,L_k\right\rangle_{\rL,N}.
\end{equation}
Then, we get the modal coefficients from three matrix--vector multiplications with the
one--dimensional \textit{forward} matrix~$F$ in each of the local coordinates directions,
\begin{equation}
    \hat{\stvec{Q}}_{ijk} = \sum_{m,n,l=0}^{N}\frac{\left\langle l_m,L_i\right\rangle_{\rL,N}}
        {\Vert L_i\Vert^2_{N}}\frac{\left\langle l_n,L_j\right\rangle_{\rL,N}}{\Vert L_j\Vert^2_{N}}
        \frac{ \left\langle l_l,L_k\right\rangle_{\rL,N}  }{\Vert L_k\Vert^2_{N}}\stvec{Q}_{mnl} =
        \sum_{m,n,l=0}^{N}F_{im}F_{jn}F_{kl}\stvec{Q}_{mnl}
\end{equation}

\begin{theorem}
Parseval's discrete identity: the three--dimensional framework also satisfies a discrete version of
Parseval's identity,
\begin{equation}
    \left\langle\stvec{U},\stvec{V}\right\rangle_{\rE,N} =
        \sum_{i,j,k=0}^{N}w_{ijk}\stvec{U}_{ijk}^{T}\stvec{V}_{ijk} =
        \sum_{i,j,k=0}^{N}\Vert L_{i}\Vert_{N}^{2}\Vert L_{j}\Vert_{N}^{2}\Vert L_{k}\Vert_{N}^{2}
        \hat{\stvec{U}}^{T}_{ijk}\hat{\stvec{V}}_{ijk}.
\label{eq:DG:filtering:theorem2}
\end{equation}
\end{theorem}
\begin{proof}
We start by defining a partial modal transformation as the result of performing the transformation
to the modal space only in one direction,
\begin{equation}
    \hat{\stvec{Q}}^{\xi}_{ijk} = \sum_{m=0}^{N}F_{im}\stvec{Q}_{mjk}, \quad
    \hat{\stvec{Q}}^{\eta}_{ijk} = \sum_{m=0}^{N}F_{jm}\stvec{Q}_{imk}, \quad
    \hat{\stvec{Q}}^{\zeta}_{ijk} = \sum_{m=0}^{N}F_{km}\stvec{Q}_{ijm},
\end{equation}
which can be applied recursively in two,
\begin{equation}
    \hat{\stvec{Q}}^{\xi^i\xi^j}_{ijk} =
        \reallywidehat{\hat{\stvec{Q}}^{\xi^i}_{ijk}}^{\xi^j},
\end{equation}
and three directions,
\begin{equation}
    \hat{\stvec{Q}}^{\xi\eta\zeta}_{ijk} = \hat{\stvec{Q}}_{ijk}.
\end{equation}
Expressing the inner product of~\eqref{eq:DG:filtering:theorem2} in nodal space, and applying
Property~\ref{prop:3} in each direction at a time,
\begin{equation}
\begin{split}
    \left\langle\stvec{U},\stvec{V}\right\rangle_{\rE,N} &=
        \sum_{i,j,k=0}^{N}w_{ijk}\stvec{U}_{ijk}^{T}\stvec{V}_{ijk} =
        \sum_{i,j,k=0}^{N}\Vert L_{i}\Vert_{N}^{2}
        w_{jk}\hat{\stvec{U}}^{\xi,T}_{ijk}\hat{\stvec{V}}^{\xi}_{ijk} \\
    &= \sum_{i,j,k=0}^{N}\Vert L_{i}\Vert_{N}^{2}\Vert L_{j}\Vert_{N}^{2}
        w_{k}\hat{\stvec{U}}^{\xi\eta,T}_{ijk}\hat{\stvec{V}}^{\xi\eta}_{ijk} \\
    &= \sum_{i,j,k=0}^{N}\Vert L_{i}\Vert_{N}^{2}\Vert L_{j}\Vert_{N}^{2}\Vert L_{k}\Vert_{N}^{2}
        \hat{\stvec{U}}^{T}_{ijk}\hat{\stvec{V}}_{ijk},
\end{split}
\end{equation}
which is the equivalent of Property~\ref{prop:3} in three dimensions. \qed
\end{proof}

In three dimensions, the polynomials are also locally filtered as,
\begin{equation}
    \bar{\stvec{Q}} = \mathcal F{\star}\stvec{Q} = \sum_{i,j,k=0}^{N}\hat{\mathcal F}_{ijk}
        \hat{\stvec{Q}}_{ijk}L_i(\xi)L_j(\eta)L_k(\zeta).
\end{equation}
For the most general filter kernel with coefficients~$\hat{\mathcal F}_{ijk}$,
\begin{equation}
    \bar{\stvec{Q}}_{ijk} = \sum_{m,n,l=0}^{N}\sum_{p,q,r=0}^{N}B_{ip}B_{jq}B_{kr}
        \hat{\mathcal{F}}_{pqr}F_{pm}F_{qn}F_{rl}\stvec{Q}_{mnl},
\end{equation}
we get a fully three--dimensional tensor product (i.e. with~$(N+1)^3$ operations). However,
computations are drastically reduced if we restrict ourselves to a tensor product version of the
filter,~${\hat{\mathcal F}_{ijk}=\hat{\mathcal F}_{i}\hat{\mathcal F}_{j}\hat{\mathcal F}_{k}}$,
which reduces the filtering to three one--dimensional matrix--vector multiplications (i.e.~$3(N+1)$
operations),
\begin{equation}
    \bar{\stvec{Q}}_{ijk} = \sum_{m,n,l=0}^{N}
        \mathcal{H}_{im}\mathcal{H}_{jn}\mathcal{H}_{kl}\stvec{Q}_{mnl} =
        \sum_{m=0}^{N}\mathcal{H}_{im}\sum_{n=0}^{N}\mathcal{H}_{jn}\sum_{l=0}^{N}\mathcal{H}_{kl}
        \stvec{Q}_{mnl},\quad
    \mathcal H = B\hat{\mathcal F}F.
\end{equation}

\begin{figure}
    \centering
    \subfigure[High--pass filter: $\hat{\mathcal F}_{ij} = \hat{\mathcal F}_{i}^{\mathrm{1D}}
        \hat{\mathcal F}_{j}^{\mathrm{1D}}$\label{fig:DG:filtering:3D:hihgpass}]{%
            \includegraphics[width=0.45\textwidth,trim=1cm 1cm 2cm 1cm,clip]
            {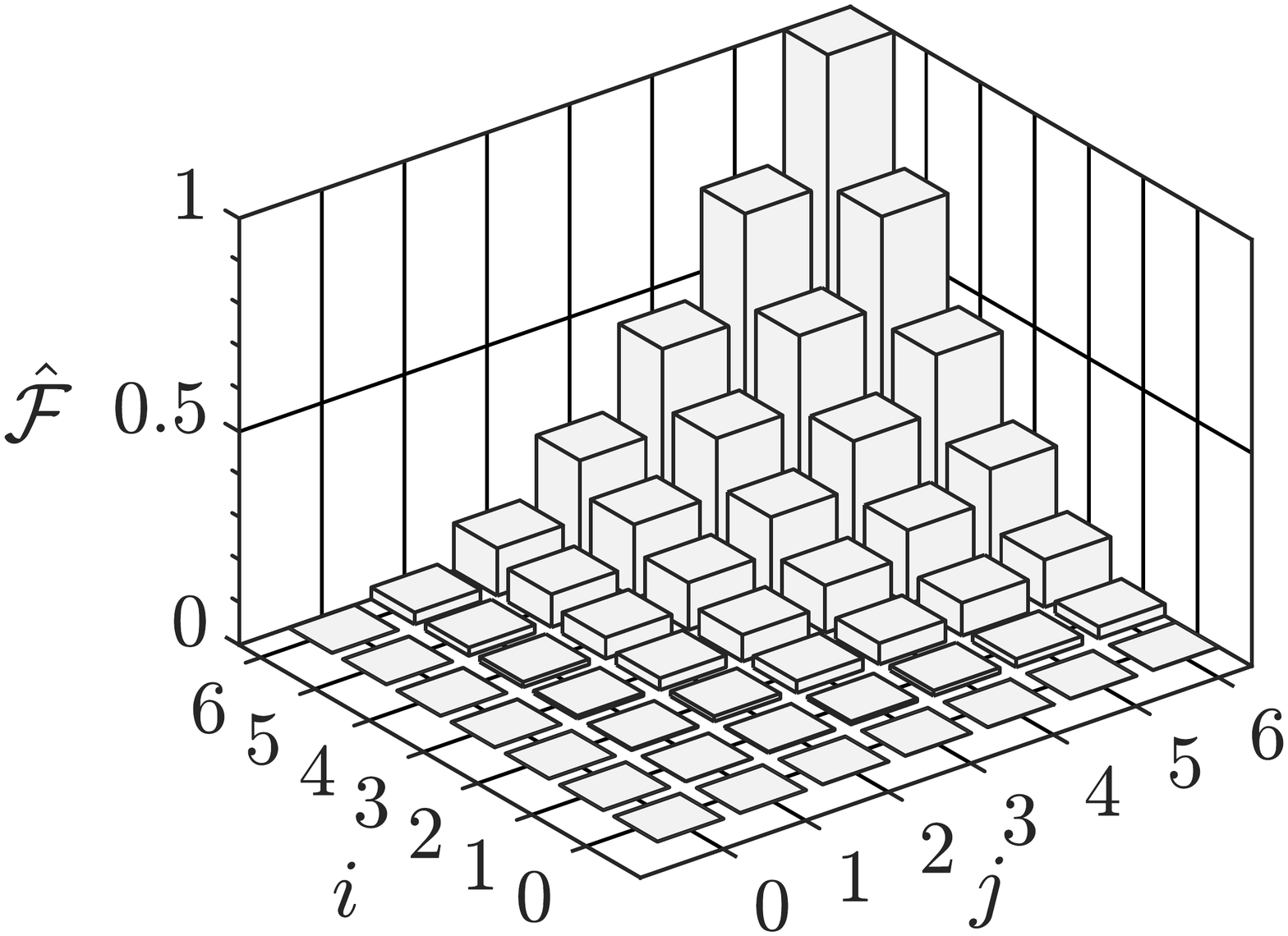}}
    \subfigure[No high--pass \mbox{filter:}~${\hat{\mathcal F}_{ij} =
        1-\left(1-\hat{\mathcal F}_{i}^{\mathrm{1D}}\right)
        \left(1-\hat{\mathcal F}_{j}^{\mathrm{1D}}\right)}$\label{fig:DG:filtering:3D:lowpass}]{%
            \includegraphics[width=0.45\textwidth,trim=1cm 1cm 2cm 1cm,clip]
            {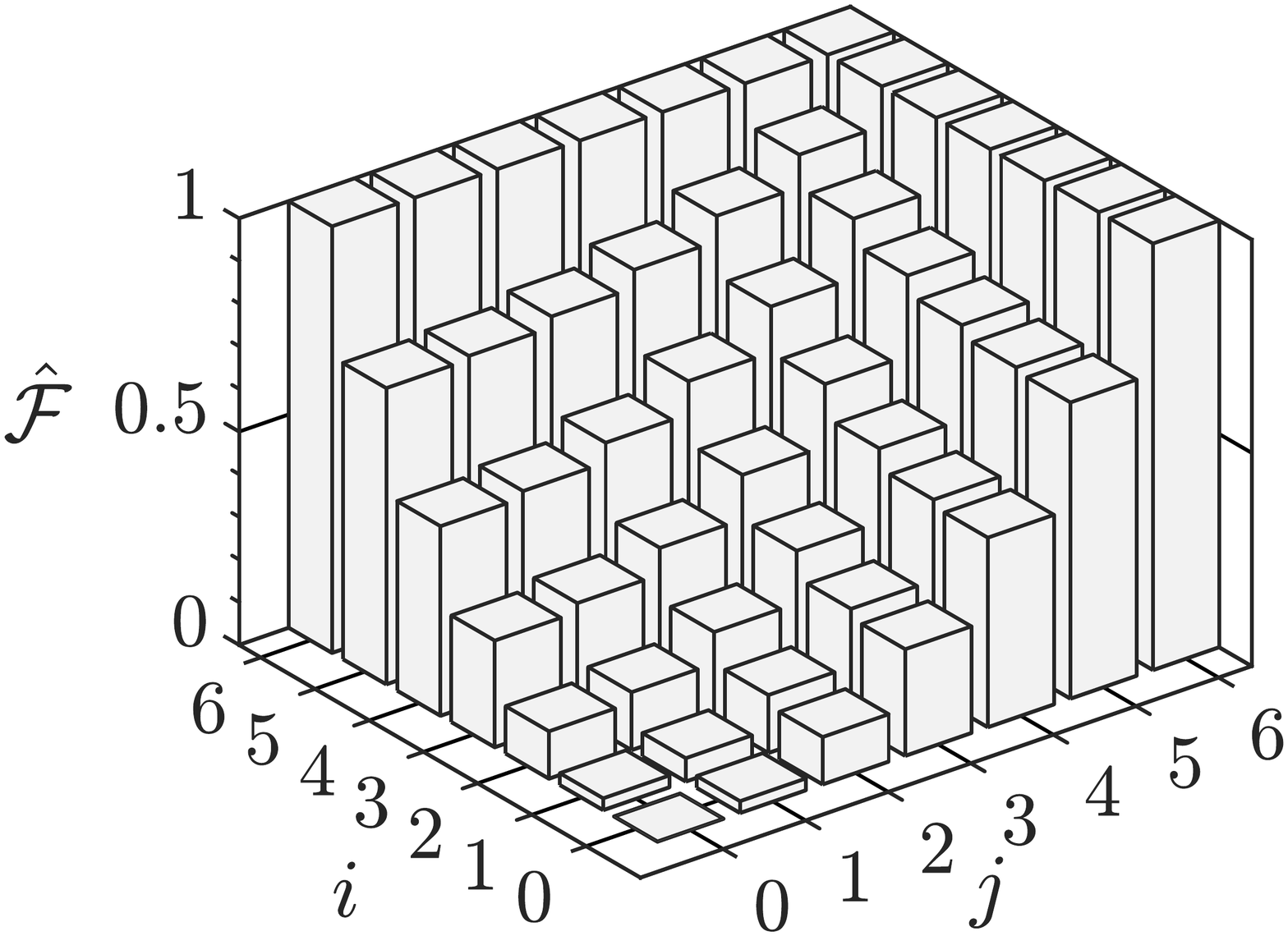}}
    \caption{Design of multidimensional filters from the tensor product of a one--dimensional
             filter. Two options are provided: a high--pass filter, where for example
             $\hat{\mathcal F}_{Nj}=\hat{\mathcal F}^{\mathrm{1D}}_{j}$, and a non high--pass
             version, where~$\hat{\mathcal F}_{Nj}=1$. Recall that modes with high filter
             coefficient values are more dissipated. For this plot,~$N=6$ and
             $\hat{\mathcal F}_{i}^{\mathrm{1D}}=(i/N)^{2}$}
\end{figure}

The methodology presented here is valid for any positive filtering
functions~$\hat{\mathcal{F}}_{ijk}$. However, we only present numerical results for the tensor
product filter ${\hat{\mathcal{F}}_{ijk}=\hat{\mathcal{F}}_{i}\hat{\mathcal{F}}_{j}
\hat{\mathcal{F}}_{k}}$. The matrix~$\mathcal{H}$ is that defined for the one--dimensional
filtering. From one--dimensional filters,~$\hat{F}^{1D}_{i}$, we can construct the three dimension
version in a high--pass version, Fig.~\ref{fig:DG:filtering:3D:hihgpass},
\begin{equation}
    \hat{\mathcal F}_{ijk} = \hat{\mathcal F}_{i}^{1D}\hat{\mathcal F}_{j}^{1D}
        \hat{\mathcal F}_{k}^{1D},
\end{equation}
or in a non high--pass version, Fig.~\ref{fig:DG:filtering:3D:lowpass},
\begin{equation}
    \hat{\mathcal F}_{ijk} = 1-\left(1-\hat{\mathcal F}_{i}^{1D}\right)
        \left(1-\hat{\mathcal F}_{j}^{1D}\right)\left(1-\hat{\mathcal F}_{k}^{1D}\right).
\end{equation}

Finally, we extend Theorem~\ref{theorem:1} to three dimensions,
\begin{prop}\label{prop:6}
The inner product of the filtered state vector~$\bar{\stvec{Q}}$ and the state vector~$\stvec{Q}$
is positive. We compute the inner product in the modal coefficients,
\begin{equation}
    \left\langle\stvec{Q},\bar{\stvec{Q}}\right\rangle_{\rE,N} =
        \sum_{i,j,k=0}^{N}\Vert L_i\Vert_N^2\Vert L_j\Vert_N^2\Vert L_k\Vert_N^2
        \hat{\stvec{Q}}_{ijk}\hat{\bar{\stvec{Q}}}_{ijk} =
        \sum_{i,j,k=0}^{N}\Vert L_i\Vert_N^2\Vert L_j\Vert_N^2\Vert L_k\Vert_N^2
        \hat{\mathcal F}_{ijk}\hat{\stvec{Q}}_{ijk}^2 \geqslant 0,
\end{equation}
that we can represent as,
\begin{equation}
    \left\langle {\stvec{Q}},\bar{\stvec{Q}}\right\rangle_{\rE,N} =
        \left\langle \stvec{Q},\mathcal F{\star}\stvec{Q}\right\rangle_{\rE,N} \geqslant 0.
\end{equation}
\end{prop}

\section{Entropy stable filtered dissipation}\label{sec:entropy-stable-filtered}

We now apply the findings in \Sec\ref{sec:DG:SVV} to construct entropy stable filtered dissipative
fluxes for the Navier--Stokes equations. Depending on the form of the
fluxes,~$\ssvec{F}_{a}=\bmat{B}_{a}(\stvec{Q})\ssvec{G}$, we proceed differently:
\begin{enumerate}
\item $\bmat{B}_{a}$ is a non--constant coefficient,~$\alpha(\stvec{q}) \geqslant 0$, times a
constant positive definite
matrix,~$\bmat{C}$,~${\bmat{B}_{a}=\alpha\left(\stvec{q}\right)\bmat{C}}$: the filtered flux is,
\begin{equation}
    \ssvec{F}_{a}^{\mathcal F} = \sqrt{\frac{\alpha}{\mathcal J}}\bmat{C}\mathcal F{\star}
        \left(\sqrt{\mathcal J\alpha}\ssvec{G}\right).
\label{eq:DG:SVV2}
\end{equation}
For example, the physical dissipation computed with the kinetic energy variables as in
\Sec\ref{sec:NSE:kin-visc}.

\item $\bmat{B}_{a}$ is a general positive definite matrix: any positive definite matrix has a
Cholesky decomposition (see~\cite{Fisher2013}),~$\bmat{B}_{a}=\bmat{L}^{T}\bmat{D}\bmat{L}$. Thus,
we define the filtered flux as,
\begin{equation}
    \ssvec{F}_a^{\mathcal F} = \frac{1}{\sqrt{\mathcal J}}\bmat{L}^T\sqrt{\bmat{D}}\mathcal F \star
        \left(\sqrt{\mathcal J\bmat{D}}\bmat{L}\ssvec{G}\right),
\label{eq:DG:SVV3}
\end{equation}
where the square root of a diagonal matrix represents a diagonal matrix whose entries are the square
roots of the original matrix. Examples of this form are the physical dissipation in thermodynamic
entropy variables, \Sec\ref{sec:NSE:ther-visc}, and the Guermond--Popov fluxes,
\Sec\ref{sec:NSE:GP}.

\end{enumerate}

The specific form of these fluxes is later justified by stability analysis. We conclude this section
by filtering the artificial viscosity flux in the volume and surface quadratures of the DG
scheme~\eqref{eq:dg:pol:NSE-split-form},
\begin{equation}
\begin{split}
    \left\langle \mathcal J\stvec{Q}_{t},\stvecg{\phi}\right\rangle_{\rE,N} &+
        \int_{\partial e,N}\stvecg{\phi}^{T}\left(\ssvec{F}_{e}^{\star} -
        \ssvec{F}_{e}\right)\cdot\svec{n}\diff{S} +
        \left\langle \mathbb{D}\left(\ssvec{F}_{e}^{\#}\right),\stvecg{\phi}\right\rangle_{\rE,N} \\
    =& \int_{\partial e,N}\stvecg{\phi}^{T}\left(\ssvec{F}_{v}^{\star} +
        \ssvec{F}_{a}^{\mathcal F,\star}\right)\cdot\svec{n}\diff{S} - \left\langle \cssvec{F}_{v} +
        \cssvec{F}_{a}^{\mathcal F},\svec{\nabla}_{\xi}\stvecg{\phi}\right\rangle_{\rE,N}.
\end{split}
\label{eq:dg:filtering:NSE}
\end{equation}

\subsection{Stability analysis}\label{sec:stability}

We study the stability of the DGSEM with filtered artificial viscosity. Since inviscid and viscous
terms have already been studied in other works~\cite{gassner2016split,Gassner2018}, we restrict
ourselves to outline the main ideas to focus on the novel form of the filtered dissipation.

First, we apply the discrete Gauss law to the gradient equation~\eqref{eq:dg:pol:NSE-split-form:G},
and we replace its test function by the viscous and filtered artificial viscosity fluxes
$\ssvecg{\varphi}=\ssvec{F}_{v}+\ssvec{F}_{a}^{\mathcal F}$,
\begin{equation}
    \left\langle \mathcal J\ssvec{G},\ssvec{F}_{v}+\ssvec{F}_{a}^{\mathcal F}\right\rangle_{\rE,N} =
    \int_{\partial e,N}\left(\stvec{W}^{\star}-\stvec{W}\right)^{T}\left(\ssvec{F}_{v} +
        \ssvec{F}_{a}^{\mathcal F}\right)\cdot\svec{n}\diff{S} +
        \left\langle\svec{\nabla}_{\xi}\stvec{W},\cssvec{F}_{v} +
        \cssvec{F}_{a}^{\mathcal F}\right\rangle_{\rE,N}.
\label{eq:stability:Grad}
\end{equation}
Second, we replace the test function~$\stvecg{\phi}=\stvec{W}$ in~\eqref{eq:dg:filtering:NSE}
\begin{equation}
\begin{split}
    \left\langle \mathcal J\stvec{Q}_{t},\stvec{W}\right\rangle_{\rE,N} &+
        \int_{\partial e,N}\stvec{W}^{T}\left(\ssvec{F}_{e}^{\star} -
        \ssvec{F}_{e}\right)\cdot\svec{n}\diff{S} +
        \left\langle \mathbb{D}\left(\ssvec{F}_{e}^{\#}\right),\stvec{W}\right\rangle_{\rE,N} \\
    =& \int_{\partial e,N}\stvec{W}^{T}\left(\ssvec{F}_{v}^{\star} +
        \ssvec{F}_{a}^{\mathcal F,\star}\right)\cdot\svec{n}\diff{S} - \left\langle \cssvec{F}_{v} +
        \cssvec{F}_{a}^{\mathcal F},\svec{\nabla}_{\xi}\stvec{W}\right\rangle_{\rE,N},
\end{split}
\label{eq:stability:NSE-repl-test}
\end{equation}
and we replace the viscous and artificial viscosity volume term from~\eqref{eq:stability:Grad}
into~\eqref{eq:stability:NSE-repl-test} to get a single equation,
\begin{equation}
\begin{split}
    &\left\langle \mathcal J\stvec{Q}_{t},\stvec{W}\right\rangle_{\rE,N} +
        \int_{\partial e,N}\stvec{W}^{T}\left(\ssvec{F}_{e}^{\star} -
        \ssvec{F}_{e}\right)\cdot\svec{n}\diff{S} +
        \left\langle \mathbb{D}\left(\ssvec{F}_{e}^{\#}\right),\stvec{W}\right\rangle_{\rE,N} \\
    =& \int_{\partial e,N}\left(\stvec{W}^{T}\left(\ssvec{F}_{v}^{\star} +
        \ssvec{F}_{a}^{\mathcal F,\star}\right) + \left(\stvec{W}^{\star} -
        \stvec{W}\right)^{T}\left(\ssvec{F}_{v} +
        \ssvec{F}_{a}^{\mathcal F}\right)\right)\cdot\svec{n}\diff{S} -
        \left\langle \mathcal J\ssvec{G},\ssvec{F}_{v} +
        \ssvec{F}_{a}^{\mathcal F}\right\rangle_{\rE,N}.
\end{split}
\label{eq:stability:NSE-single-eq}
\end{equation}

We first study the volume terms. From~\cite{Gassner2018,Gassner2014},
\begin{equation}
\begin{aligned}
    \left\langle \mathcal J\stvec{Q}_{t},\stvec{W}\right\rangle_{\rE,N} &=
        \left\langle\mathcal J\mathcal E_{t}\right\rangle_{\rE,N}, \\
    \left\langle\mathbb{D}(\ssvec{F}_{e}^{\#}),\stvec{W}\right\rangle_{\rE,N} &=
        \int_{\partial e,N}\svec{F}^{\mathcal E}_{e}\cdot\svec{n}\diff{S}+\Theta\mathcal W_p, \\
    \left\langle\mathcal J\ssvec{G},\ssvec{F}_{v}\right\rangle_{\rE,N} &=
        \left\langle\mathcal J\ssvec{G},\bmat{B}_{v}\ssvec{G}\right\rangle_{\rE,N} =
        \langle D_{v}\rangle^{\rE,N} \geqslant 0.
\end{aligned}
\end{equation}
Therefore, we transform~\eqref{eq:stability:NSE-single-eq} into
\begin{equation}
\begin{split}
    \left\langle\mathcal J\mathcal E_{t}\right\rangle_{\rE,N} &+
        \int_{\partial e,N}\stvec{W}^{T}\left(\ssvec{F}_{e}^{\star} -
        \ssvec{F}_{e}\right)\cdot\svec{n}\diff{S} +
        \int_{\partial e,N}\svec{F}^{\mathcal E}_{e}\cdot\svec{n}\diff{S}+
        \Theta \mathcal W_p + \langle D_{v}\rangle^{\rE,N} +
        \left\langle\mathcal J \ssvec{G},\ssvec{F}_{a}^{\mathcal F}\right\rangle_{\rE,N} \\
    =& \int_{\partial e,N}\left(\stvec{W}^{T}\left(\ssvec{F}_{v}^{\star} +
        \ssvec{F}_{a}^{\mathcal F,\star}\right) +
        \left(\stvec{W}^{\star}-\stvec{W}\right)^{T}\left(\ssvec{F}_{v} +
        \ssvec{F}_{a}^{\mathcal F}\right)\right)\cdot\svec{n}\diff{S}.
\end{split}
\end{equation}

We now study the stability of the volume term associated to the filtered artificial viscosity. We
study the two options provided in \Sec\ref{sec:entropy-stable-filtered}.
\begin{enumerate}
\item The non--constant coefficient times a constant positive definite matrix~\eqref{eq:DG:SVV2}
\begin{equation}
    \langle D_{a}\rangle^{\rE,N} = \left\langle\mathcal J\ssvec{G},
        \sqrt{\frac{\alpha}{\mathcal J}}\bmat{C}\mathcal F{\star}
        \left(\sqrt{\mathcal J\alpha}\ssvec{G}\right)\right\rangle_{\rE,N} =
        \left\langle \sqrt{\mathcal J\alpha}\ssvec{G},\bmat{C}\mathcal F{\star}
        \left(\sqrt{\mathcal J\alpha}\ssvec{G}\right)\right\rangle_{\rE,N} \geqslant 0.
  \end{equation}

\item The general positive definite matrix \eqref{eq:DG:SVV3}
\begin{equation}
    \langle D_{a}\rangle^{\rE,N} = \left\langle\mathcal J\ssvec{G},
        \frac{1}{\sqrt{\mathcal J}}\bmat{L}^{T}\sqrt{\bmat{D}}\mathcal F{\star}
        \left(\sqrt{\mathcal J\bmat{D}}\bmat{L}\ssvec{G}\right)\right\rangle_{\rE,N} =
    \left\langle\sqrt{\mathcal J\bmat{D}}\bmat{L}\ssvec{G},\mathcal F{\star}
        \left(\sqrt{\mathcal J\bmat{D}}\bmat{L}\ssvec{G}\right)\right\rangle_{\rE,N} \geqslant 0.
  \end{equation}
\end{enumerate}

Where we used Property~\ref{prop:6} to conclude that the two possibilities introduce volume
dissipative terms. A similar conclusion was reached by Lundquist and Nordstr\"om
in~\cite{Lundquist2020}, where they defined the filter kernel
as~$\mathcal{J}^{-1/2}\mathcal{F}\mathcal{J}^{1/2}$ to ensure the so called \emph{time--stability}.
Therefore, we can write,
\begin{equation}
\begin{split}
    \left\langle\mathcal J\mathcal E_{t}\right\rangle_{\rE,N} &+
        \int_{\partial e,N}\stvec{W}^{T}\left(\ssvec{F}_{e}^{\star} -
        \ssvec{F}_{e}\right)\cdot\svec{n}\diff{S} +
        \int_{\partial e,N}\svec{F}^{\mathcal E}_{e}\cdot\svec{n}\diff{S} +
        \Theta \mathcal W_p + \langle D_{v}\rangle^{\rE,N} +
        \langle D_{a}\rangle^{\rE,N} \\
    =& \int_{\partial e,N}\left(\stvec{W}^{T}\left(\ssvec{F}_{v}^{\star} +
        \ssvec{F}_{a}^{\mathcal F,\star}\right)+\left(\stvec{W}^{\star}-\stvec{W}\right)^{T}
        \left(\ssvec{F}_{v}+\ssvec{F}_{a}^{\mathcal F}\right)\right)\cdot\svec{n}\diff{S}.
\end{split}
\label{eq:stability:stability-1-element}
\end{equation}

With the volume terms addressed, we study the boundary terms. To do this, we
sum~\eqref{eq:stability:stability-1-element} over all the mesh elements,
\begin{equation}
    \langle\mathcal E_t\rangle^{N} + \mathrm{IBT} + \mathrm{PBT} +
        \sum_{e} \langle D_{v}\rangle^{\rE,N} +
        \sum_{e} \langle D_{a}\rangle^{\rE,N} + \Theta\sum_{e}\mathcal{W}_p = 0,
\label{eq:stability:stability-all-elements}
\end{equation}
where $\langle\mathcal E_{t}\rangle^{N}=\sum_{e}\left\langle \mathcal J\mathcal
E_t\right\rangle_{\rE,N}$ is the total entropy derivative, $\mathrm{IBT}$ are the interior boundary
terms,
\begin{equation}
\begin{split}
    \mathrm{IBT} = \mathrm{IBT}_{e} + \mathrm{IBT}_{v} + \mathrm{IBT}_{a} =
        &-\sum_{\interiorfaces}\int_{f,N}\left(\jump{\svec{F}_{e}^{\mathcal E}} +
        \jump{\stvec{W}}^{T}\ssvec{F}_{e}^{\star}-\jump{\stvec{W}^{T}\ssvec{F}_{e}}\right)
        \cdot\svec{n}_{L} \diff{S} \\
    &+ \sum_{\interiorfaces}\int_{f,N}\left(\jump{\stvec{W}}^{T}\ssvec{F}_{v}^{\star} +
        \stvec{W}^{\star,T}\jump{\ssvec{F}_{v}}-\jump{\stvec{W}^{T}\ssvec{F}_{v}}\right)
        \cdot\svec{n}_{L} \diff{S} \\
    &+ \sum_{\interiorfaces}\int_{f,N}\left(\jump{\stvec{W}}^{T}\ssvec{F}_{a}^{\mathcal F, \star} +
        \stvec{W}^{\star, T}\jump{\ssvec{F}^{\mathcal F}_{a}} -
        \jump{\stvec{W}^{T}\ssvec{F}^{\mathcal F}_{a}}\right) \cdot\svec{n}_{L} \diff{S},
\end{split}
\end{equation}
and PBT are the physical boundary terms,
\begin{equation}
\begin{split}
    \mathrm{PBT} = \mathrm{PBT}_{e}+\mathrm{PBT}_{v}+\mathrm{PBT}_{a} =
        &\sum_{\boundaryfaces}\int_{f,N}\left(\svec{F}_{e}^{\mathcal E} +
        \stvec{W}^{T}\ssvec{F}_{e}^{\star}-\stvec{W}^{T}\ssvec{F}_{e}\right)
        \cdot\svec{n}\diff{S} \\
    &-\sum_{\boundaryfaces}\int_{f,N}\left({\stvec{W}}^{T}\ssvec{F}_{v}^{\star} +
        \stvec{W}^{\star,T}{\ssvec{F}_{v}}-{\stvec{W}^{T}\ssvec{F}_{v}}\right)
        \cdot\svec{n}\diff{S} \\
    &-\sum_{\boundaryfaces}\int_{f,N}\left({\stvec{W}}^{T}\ssvec{F}_{a}^{\mathcal F, \star} +
        \stvec{W}^{\star,T}{\ssvec{F}^{\mathcal F}_{a}} -
        {\stvec{W}^{T}\ssvec{F}^{\mathcal F}_{a}}\right) \cdot\svec{n} \diff{S}.
  \end{split}
\end{equation}
The interior boundary terms for the inviscid and viscous terms have been studied
in~\cite{gassner2016split,Gassner2018}. If one uses the two--point entropy conserving flux as the
numerical flux,~$\mathrm{IBT}_{e}=0$, whereas for the dissipative
flux,~$\mathrm{IBT}_{e}\geqslant 0$. For viscous fluxes, the BR1 scheme gives~$\mathrm{IBT}_{v}=0$.
We have also used the BR1 scheme for the artificial dissipative fluxes, and therefore we also
get~$\mathrm{IBT}_{a}=0$,
\begin{equation}
    \jump{\stvec{W}}^{T}\ssvec{F}_{a}^{\mathcal F, \star} +
        \stvec{W}^{\star,T}\jump{\ssvec{F}^{\mathcal F}_{a}} -
        \jump{\stvec{W}^{T}\ssvec{F}^{\mathcal F}_{a}} =
    \jump{\stvec{W}}^{T}\aver{\ssvec{F}^{\mathcal F}_{a}} +
        \aver{\stvec{W}}^{T}\jump{\ssvec{F}^{\mathcal F}_{a}} -
        \jump{\stvec{W}^{T}\ssvec{F}^{\mathcal F}_{a}} = 0.
\end{equation}

Similarly, the physical boundary terms for the inviscid and viscous terms were
studied in~\cite{hindenlang2020stability} for the wall boundary condition. This analysis concluded
that~$\mathrm{PBT}_{e}=0$ if one simply imposes the boundary data,
and~$\mathrm{PBT}_{e}\geqslant 0$ if one uses the dissipative numerical flux with a ghost state. In
the case of the viscous fluxes,~$\mathrm{PBT}_{v}=0$, and for the artificial viscosity fluxes, the
physical boundary terms are also zero. For the no--slip wall,
\begin{equation}
    \stvec{W}^{T}\ssvec{F}_{a}^{\mathcal F,\star} +
        \left(\stvec{W}^{\star}-\stvec{W}\right)^{T}\ssvec{F}^{\mathcal F}_{a} =
        \stvec{W}^{T}\left(0, F_{a,2}, F_{a,3}, F_{a,4}, 0\right) +
        \left(0, -W_{2}, -W_{3}, -W_{4}, 0\right)^{T}\ssvec{F}^{\mathcal F}_{a} = 0,
\end{equation}
and for the free--slip wall,
\begin{equation}
    \stvec{W}^{T}\ssvec{F}_{a}^{\mathcal F, \star} +
        \left(\stvec{W}^{\star}-\stvec{W}\right)^{T}\ssvec{F}^{\mathcal F}_{a} =
        \left(\stvec{W}-\stvec{W}\right)\ssvec{F}^{\mathcal F}_{a} = 0.
\end{equation}

We conclude the stability analysis as we confirm that the thermodynamic entropy is monotonic and
decreases due to the numerical dissipation of the inviscid numerical flux at the inter--element and
physical faces, and the artificial and physical dissipation at the interior of the elements,
\begin{equation}
    \langle\mathcal E_t\rangle^{N} + \Theta\sum_{e}\mathcal{W}_p =
        -\mathrm{IBT}_{e} - \mathrm{PBT}_{e} -
        \sum_{e} \left\langle D_{v}\right\rangle^{\rE,N} -
        \sum_{e} \left\langle D_{a}\right\rangle^{\rE,N} \leqslant 0.
\label{eq:stability:stability-final}
\end{equation}
Regarding the kinetic energy, the work introduced by the pressure gradient cannot be bounded.
However, although the kinetic energy is not a mathematical entropy in the strict sense, schemes
based on kinetic energy preservation~\cite{Jameson2008} represent a popular choice for subsonic
flows, given that the resulting schemes are still robust and have a lower computational cost than
those based on the thermodynamic entropy (which are the preferred option to solve supersonic flows,
where oscillations in the energy equation become important).

In any case, we are able to adjust the artificial dissipation~$D_{a}^{\rE,N}$ with the
selection of the artificial viscosities, and the filter kernel~$\mathcal F$.

\subsection{Linear von Neumann analysis}\label{sec:vonNeumann}

Before continuing with the numerical experiments, as in other
works~\cite{Moura2016,Chavez2018,Manzanero2020,solan2021application,kou2021eigensolution}, we
perform a von Neumann analysis to characterise the SVV in the one--dimensional advection--diffusion
equation with constant coefficients. This will allow us to better understand the dissipative
behaviour of the SVV filtering technique at different scales, knowledge that we can use to
understand the results we find in the posterior Navier--Stokes experiments in \Sec\ref{sec:results}.
Starting from the continuous formulation,
\begin{equation}
    u_t + au_x = \mu u_{xx},
\label{eq:vonneumann:advdiff}
\end{equation}
where~$a>0$ is the advection velocity and~$\mu>0$, the viscosity, we introduce a wave-like
solution~${u(x,t) = e^{i(kx-\omega t)} = u_0(x)e^{-i\omega t}}$, to determine the dispersion
relation,~$\omega(k)$, which shows the behavior of a monochromatic wave when traversing the domain
of interest. In this case,
\begin{equation}
    -i\omega + iak = -\mu k^2, \quad \omega = ak - i\mu k^2,
\label{eq:vonneumann:wk}
\end{equation}
implying that~$u(x,t)=e^{-\mu k^2t}e^{ik(x-at)}$ and thus, the wave velocity is~$a$ and the
dissipation is~$\mu k^2$. The discretization of~\eqref{eq:vonneumann:advdiff} induces errors into
the solution that lead to a modified dispersion relation. If~$\Re(\omega)\neq ak$, the numerical
scheme is introducing a dispersion error that makes the waves to move at velocities different from
the exact one and, if this velocity depends on~$k$, the different components of non--monochromatic
waves will drift, deforming it as it traverses the domain. The imaginary part of~$\omega$ is
directly related to the dissipation of the scheme and thus, a negative value means that the scheme
is stable.

\begin{figure}
    \centering
    \subfigure[Dispersion]{%
        \includegraphics[width=0.49\textwidth]{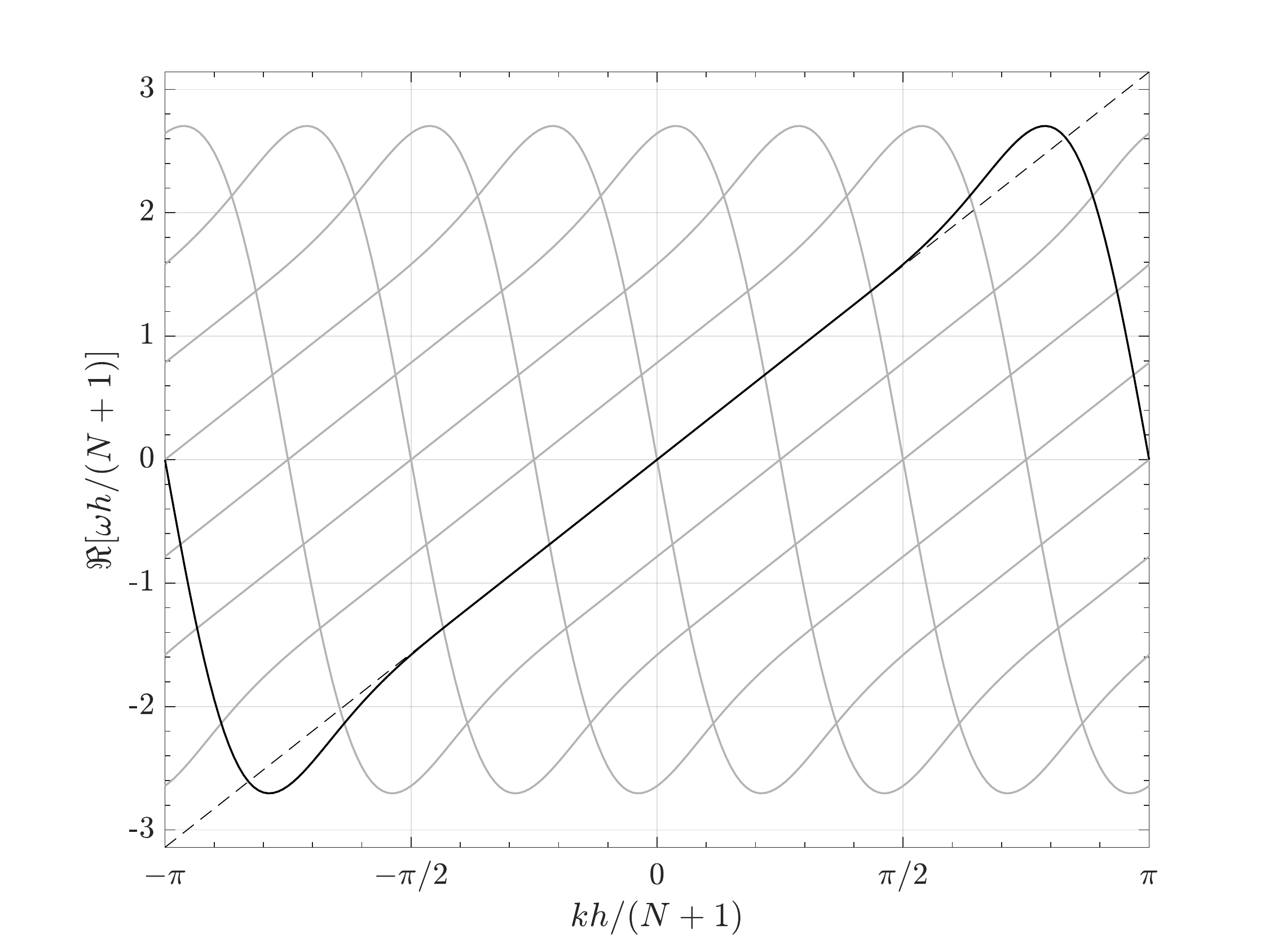}}
    \hfill
    \subfigure[Diffusion]{%
        \includegraphics[width=0.49\textwidth]{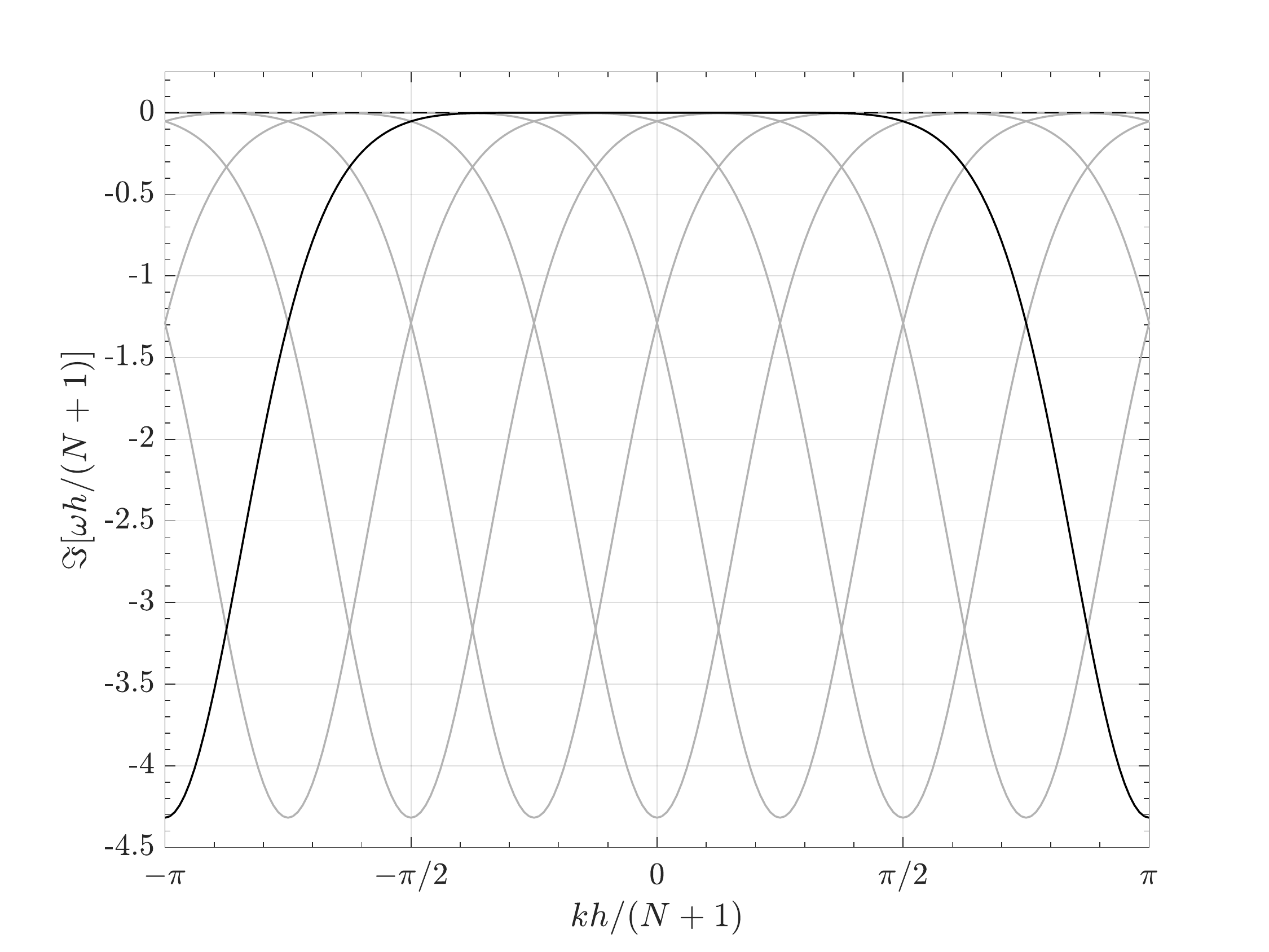}}
    \caption{Dispersion--diffusion curves of the DGSEM discretization
             of~\eqref{eq:vonneumann:advdiff} with~$N=7$ and no viscosity. Dashed lines represent
             the expected result from the continuous equation.}
    \label{fig:vonneumann:disp_diff}
\end{figure}

In this von Neumann analysis of the DGSEM with SVV dissipation we use a simple upwinding Riemann
solver for the advective term, and we keep the BR1 scheme for the viscous part,
\begin{equation}
    a U^{\star} = a\aver{U} - \frac{|a|}{2}\jump{U}, \quad
    \mu G^{\star} = \aver{\mu G}, \quad U^{\star} = \aver{U},
\end{equation}
where, as we have previously done, for instance, in~\eqref{eq:dg:geom:1st-order-ref}, we define the
gradient of the solution as~$G=U_x$. The dispersion--dissipation errors with no viscosity and~$N=7$
have been represented in Fig.~\ref{fig:vonneumann:disp_diff} to serve as a starting point.

\begin{figure}
    \centering
    \subfigure[Diffusion]{%
        \includegraphics[width=0.49\textwidth]{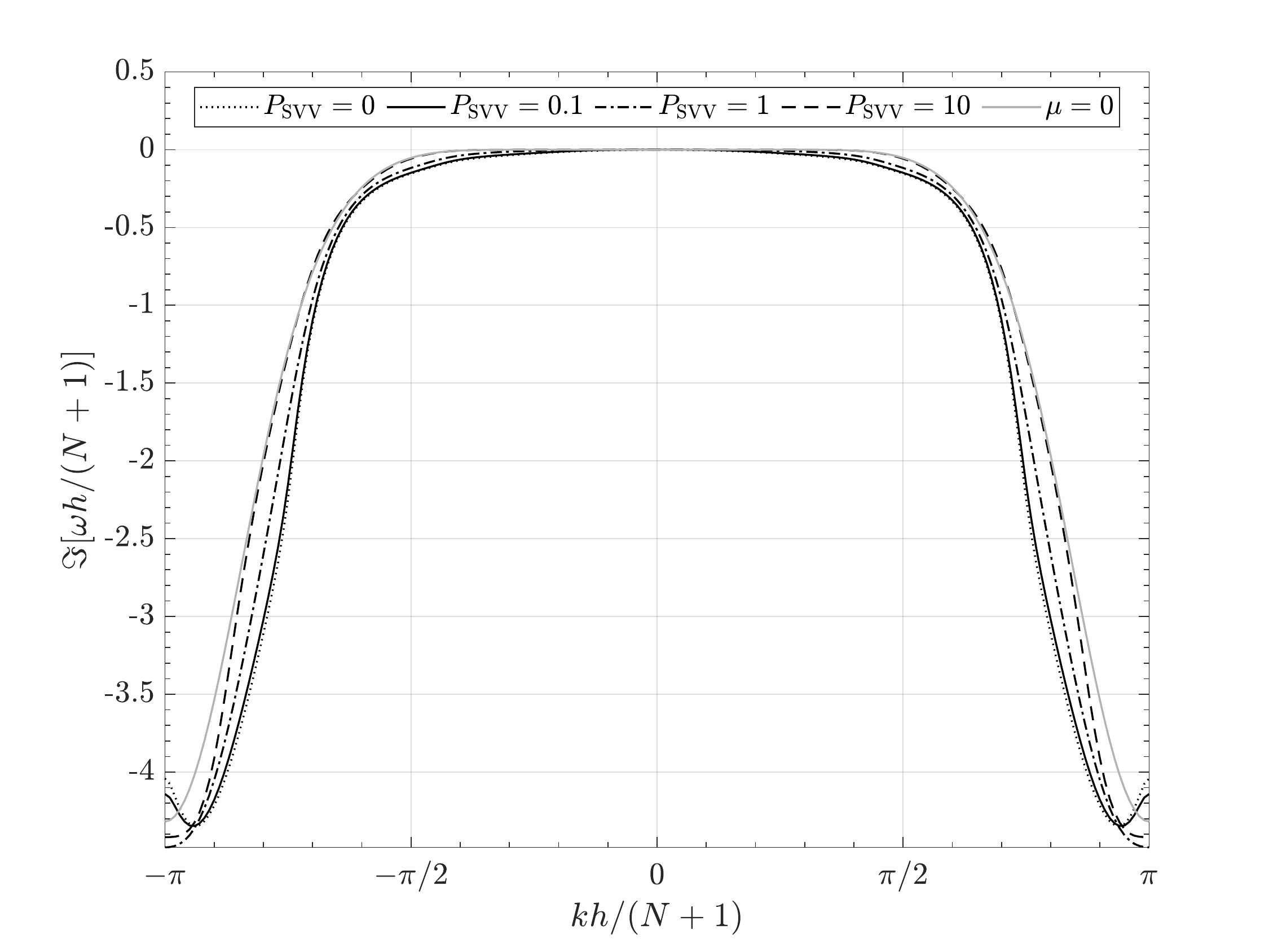}}
    \hfill
    \subfigure[Diffusion, detailed view]{%
        \includegraphics[width=0.49\textwidth]{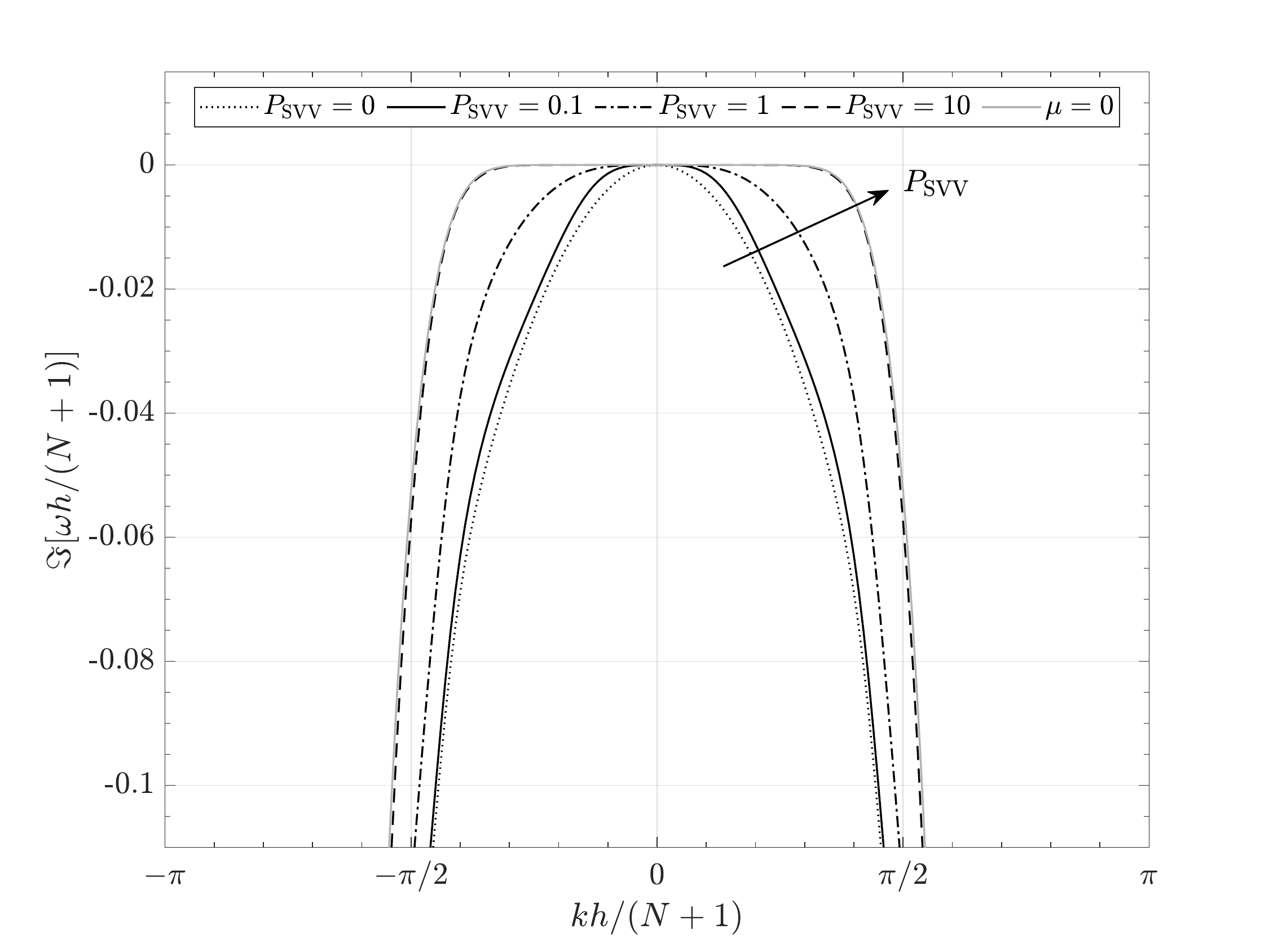}}
    \caption{Dispersion--diffusion curves with~$\mu_a=0.01$,~$N=7$ and different values
             of~$P_{\SVV}$.}
    \label{fig:vonneumann:svv_diff}
\end{figure}

The SVV is introduced in~\eqref{eq:vonneumann:advdiff} by redefining the viscous term as,
\begin{equation}
    \mu G = \mu_a \mathcal{H} G, \quad \mathcal{H} = B\hat{\mathcal F}F, \quad
        \hat{\mathcal F}_i = \diag\left[\left(\frac{i}{N}\right)^{P_{\SVV}}\right],
\label{eq:vonneumann:filter}
\end{equation}
where the matrix~$\mathcal{H}$ applies the filter kernel defined by Moura et
al.~\cite{Moura2016},~$\hat{\mathcal F}$, to the gradients,~$G$, as we showed in
\Sec\ref{sec:DG:filtering:1D}. In this formulation, the artificial viscosity is controlled with two
parameters:~$\mu_a$ and~$P_{\SVV}$, and we can now compute the curve~$\omega(k)$ as shown in
Fig.~\ref{fig:vonneumann:svv_diff} for a constant value of~$\mu_a=0.01$. We want to note that this
choice of filter kernel allows us to recover the original viscous operator by selecting~$P_{\SVV}=0$
and thus,~$\mathcal{H}=\bmat{I}$, but the only way to completely eliminate the viscosity is by
setting~$\mu_a=0$.

The results are presented in Fig.~\ref{fig:vonneumann:svv_diff} and are in agreement with those
in~\cite{Manzanero2020}. As we expected, the effect of~$P_{\SVV}$ is more evident for low frequency
waves, where the SVV is responsible for almost all the dissipation in this region. For the rest of
wavenumbers, the addition of viscosity evidently increases dissipation and, even if it can be
controlled with the filter kernel to some point, it is comparatively smaller than the one already
added by the Riemann solver. These results are also in line with our purpose, since we can still
dissipate unwanted oscillations in regions with discontinuities while keeping smoother regions
almost unaffected.

\section{Numerical experiments}\label{sec:results}

In this section we present numerical experiments where the benefits of the filtered artificial
viscosity are compared against the non--filtered DGSEM and other reference solutions from the
literature. We apply the filter kernel, taken from~\cite{Moura2016}, in all the cases,
\begin{equation}
    \hat{\mathcal F}_i^{\mathrm{1D}}=\left(\frac{i}{N}\right)^{P_{\SVV}}.
\end{equation}

\subsection{Taylor--Green vortex}

We repeat the experiment of~\cite{Manzanero2020} with the entropy stable version of the SVV. The
Taylor--Green vortex case is initialized with a periodic flow field representing a set of vortices
that evolve in time, transferring their energy to lower scales, where it is finally dissipated.
Thus, it is a good test for analyzing the dissipative features of our numerical scheme, including
the SVV. In our simulations, the computational domain is a periodic~$[0,2\pi]^3$ cube divided
into~$8^3$ elements with approximation degree~$N=8$, and the initial condition is
\begin{equation}
\begin{split}
    \rho &=  1, \\
    u    &=  \sin(x) \cos(y) \cos(z), \\
    v    &= -\cos(x) \sin(y) \cos(z), \\
    w    &=  0, \\
    p    &=  100 + \frac{1}{16} \left[\cos(2x) \cos(2z) + 2\cos(2y) +
             2\cos(2x) + \cos(2y) \cos(2z)\right]
\end{split}
\end{equation}

We aim to prove that the addition of the filtered artificial viscosity in the DGSEM does not
destabilize the energy preserving properties of the scheme, so we discretize the inviscid term with
the kinetic energy preserving split--form from Pirozzoli~\cite{Pirozzoli2010} and use a low
dissipation Roe, Riemann solver~\cite{Toro2009} to add dissipation in the region of higher
wavenumbers.

We apply a high--pass filtered Navier--Stokes viscosity to additionally increase dissipation in
lower/medium scales, controlling it through~$P_{\SVV}$. Regarding the viscosity, we compute its
value at runtime from a LES (Smagorinsky) formulation,
\begin{equation}
    \mu_a = C_{\mathrm{S}}^2 \Delta^2 |\tens{S}|, \quad
    \Delta^3 = \frac{\mathrm{Cell\,volume}}{(N+1)^3}, \quad
    S_{ij} = \frac{1}{2} \left(\frac{\partial u_i}{\partial x_j} +
        \frac{\partial u_j}{\partial x_i}\right),
\end{equation}
with~$C_{\mathrm{S}}=0.2$. In this way, we eliminate one parameter and, at the same time, we also
exploit the benefits of this LES approach to resolve turbulence. Similarly as we did in von
Neumann analysis of the one--dimensional advection--diffusion, we want to find the effect of the
artificial viscosity in different scales. To do so, we acknowledge that the theoretical kinetic
energy spectrum of the TGV obeys the law~$E(k) \propto k^{-5/3}$~\cite{Pope2001}, thus decaying for
higher wavenumbers. Due to the previously mentioned lack of numerical dissipation of the DGSEM and
the use of Pirozzoli's split--form, an accumulation of kinetic energy is expected for high
wavenumbers, where the Riemann solver and the SVV replace the dissipation that takes place in
well resolved cases.

\begin{figure}
    \centering
    \subfigure[Kinetic energy dissipation\label{fig:tgv:kindissipation}]{%
        \includegraphics[width=0.49\textwidth]{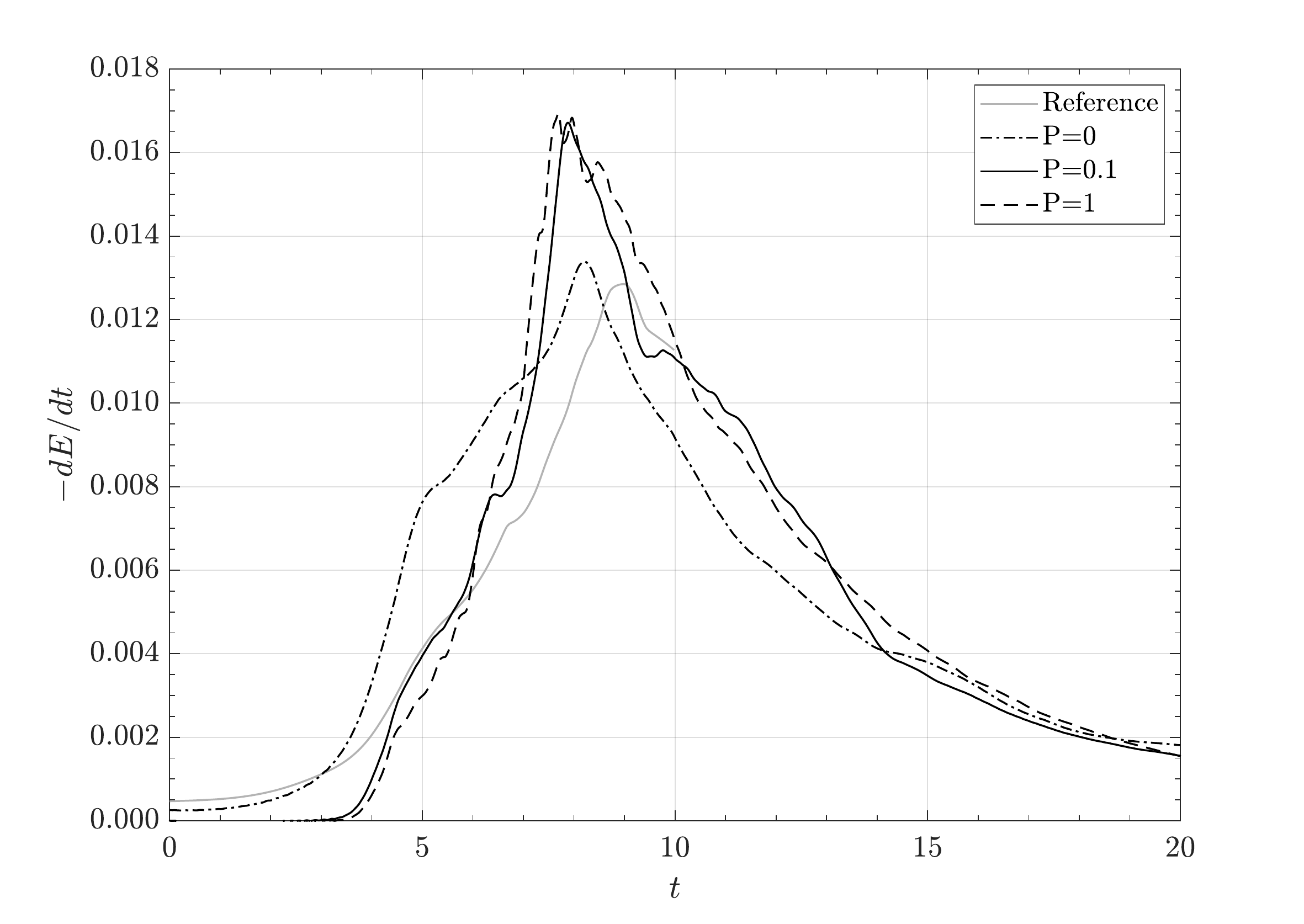}}
    \hfill
    \subfigure[Kinetic energy spectrum]{%
        \includegraphics[width=0.49\textwidth]{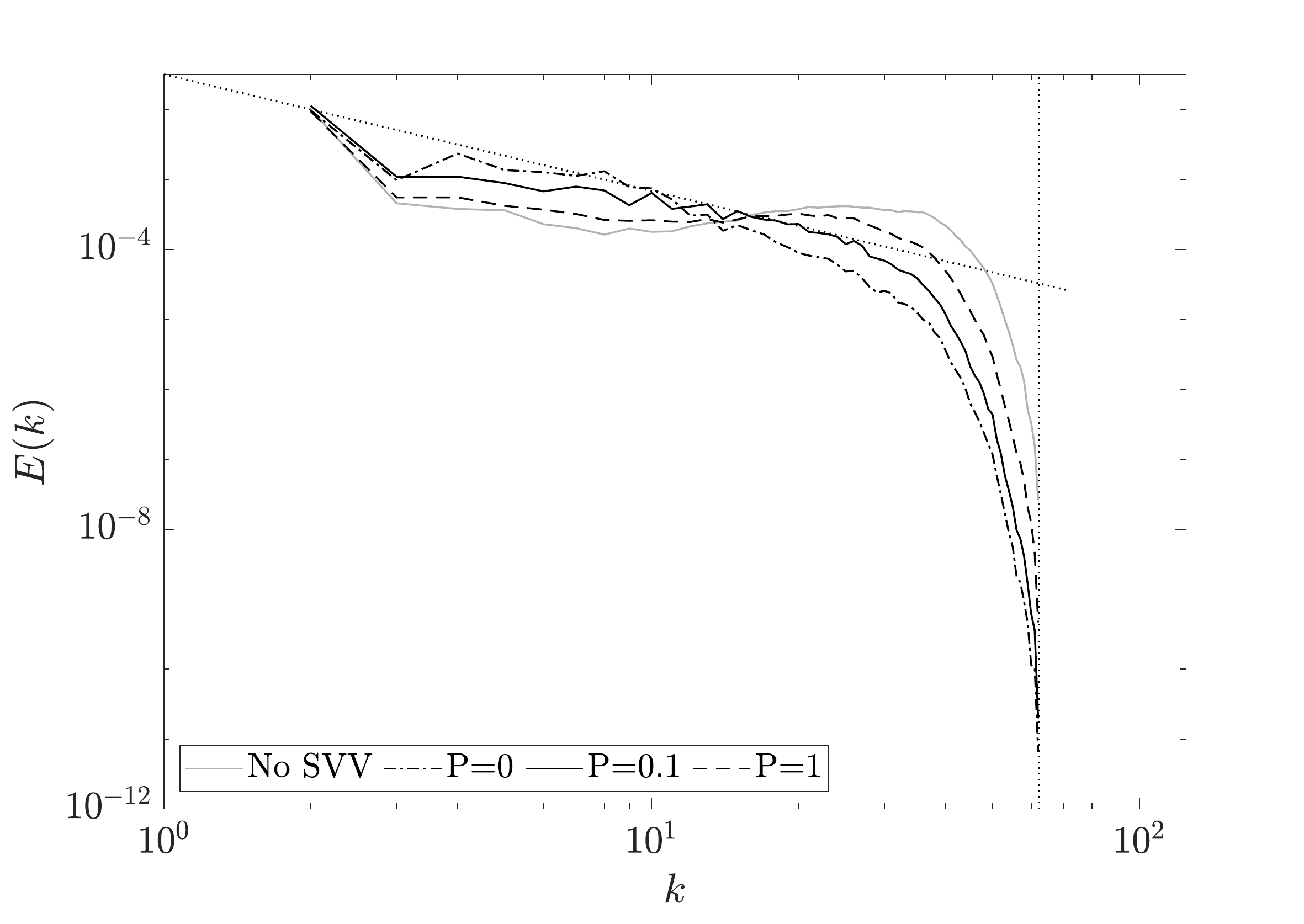}}
    \caption{Kinetic energy dissipation and spectrum of the inviscid Taylor--Green case at~$t=20$
        with a~$8^3$ elements cartesian grid and~${N=8}$. The LES viscosity is excessively
        dissipative and the SVV can control it, having an optimum value~$P_{\SVV} \approx 0.1$.
        Reference data obtained from~\cite{Rees2011} with a resolution of~$768^3$.}
    \label{fig:tgv:spectra}
\end{figure}

A representation of the evolution of the kinetic energy,~$E$, and its spectrum for different values
of~$P_{\SVV}$ at~$t=20$ (see Fig.~\ref{fig:tgv:spectra}) shows that the LES model ($P_{\SVV}=0$)
introduces excessive dissipation. This represents the perfect use case for the methodology
presented, as we aim to control this excess by adjusting the filter kernel coefficient. In agreement
with the results from von Neumann analysis, the SVV does not change significantly the behaviour of
the higher frequency modes, concentrating its effect in the low to medium wavenumbers. As we
expected, the introduction of viscosity in these regions helps to dissipate energy in all the scales
and not only in the smallest ones, avoiding its accumulation in high--frequency modes and flattening
the energy spectrum. For an optimum that we find to be~$P_{\SVV} \approx 0.1$, it closely resembles
that of the~$k^{-5/3}$ theoretical law in most parts of the range of wavenumbers resolved by the
scheme. This can also be seen in Fig.~\ref{fig:tgv:kindissipation}, where we compare the energy
dissipation of our scheme with a reference solution from Rees et al.~\cite{Rees2011}
at~$\mathrm{Re}=1600$. During the first five seconds, the flow is essentially laminar and our scheme
adds less dissipation than the reference. After that, higher-frequency modes gain importance and the
dissipation rate increases to a similar extent of that of the reference to prevent energy
accumulation in the smaller scales.

\subsection{Shock capturing: Shu--Osher problem}

The one--dimensional Shu--Osher problem describes the evolution of a shock wave that swallows a
density fluctuation as it advances. Since the solution combines a strong shock with smooth
oscillations, this benchmark is useful to assess how the dissipation is able to both, control the
shock and vanish in non--shocky regions at the same time.

The computational domain is $x \in [-4.5,4.5]$, and the initial condition is
\begin{equation}
    \left(\rho,u,p\right) =\begin{cases}
        (3.857143, 2.629369, 10.3333) & \quad \mathrm{if}~x \leqslant -4, \\
        (1+0.2\sin 5x, 0, 1)          & \quad \mathrm{if}~x > -4.
    \end{cases}
\end{equation}

\begin{figure}
    \centering
    \subfigure[Full view]{%
        \includegraphics[width=0.49\textwidth]{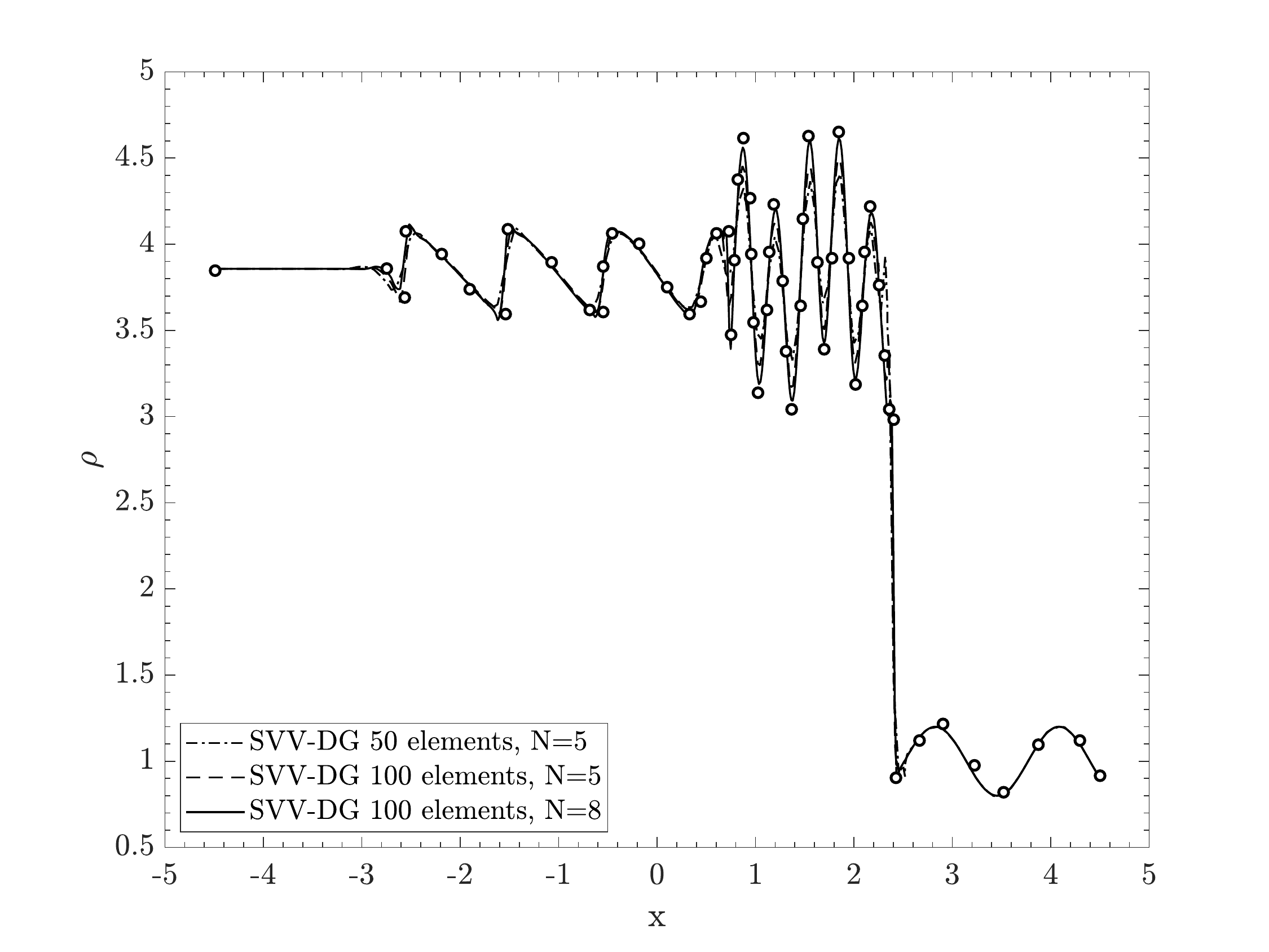}}
    \subfigure[Detailed view]{%
        \includegraphics[width=0.49\textwidth]{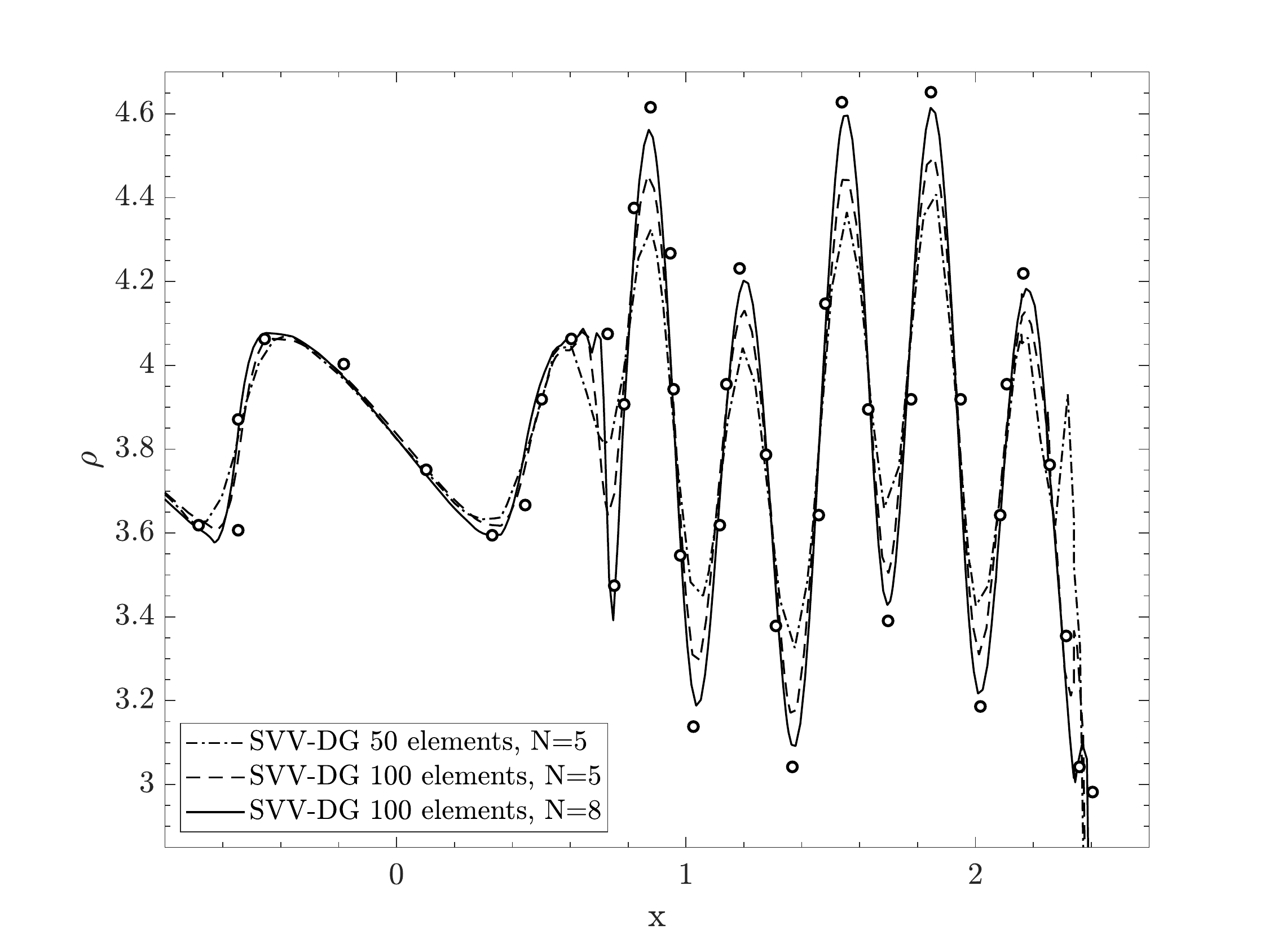}}
    \caption{Density solution of the one--dimensional Shu--Osher problem. The solid lines represent
             the DG solution with the SVV dissipation developed in this work, and the dots represent
             the reference solution of~\cite{shu1989efficient} using an ENO-RF-S-3 scheme with 1600
             degrees of freedom.}
    \label{fig:shu-osher:vs}
\end{figure}

To ensure that the solution is physically meaningful, here we use an entropy stable formulation,
discretizing the inviscid term with Chandrashekar's split--form~\cite{Chandrashekar2013} and the
matrix dissipation Riemann solver~\eqref{eq:dg:f-star:mat-diss}, and adding SVV filtered
Guermond--Popov fluxes with thermodynamic entropy variables.  We have found that the SVV filtering
is not enough to control the oscillations near the strong shock, and we have implemented a simple
sensor,~$s_1$, based on the density gradient to set~$P_{\SVV}=0$ in the discontinuity ($s_1>10$)
and~$P_{\SVV}=2$ elsewhere,
\begin{equation}
    s_1 = \sqrt{\sum_{i=0}^N w_i\left(\frac{\partial \rho}{\partial x}\right)_i^2}.
\label{eq:shu-osher:sensor}
\end{equation}
The results are shown in Fig.~\ref{fig:shu-osher:vs}, where we have tested this approach with two
different meshes and aproximation orders. The coarsest mesh contains 50, 5th order elements, while
the finest one has 100 elements with approximation orders~$N=5$ and~$N=8$, leading to 300, 600 and
900 solution nodes for each one of the three test cases, respectively. As expected, an increasing
number of nodes returns a better approximation to the reference solution
(see~\cite{shu1989efficient}); however, all the cases show a similar behaviour, reproducing all the
features of the flow. The main shock is well resolved due to the effect of the sensor. The density
regularization of the Guermond--Popov fluxes removes further fluctuations in the lower region of the
shock, and the introduction of the SVV filtering is responsible for the smoothing of the solution in
the downstream of the strong shock. This is very successful at eliminating most of the spurious
oscillations while it leaves the lower order modes, corresponding to the flow features, almost
unaltered.

\subsection{Shock capturing: Mach 3 forward facing step}

Finally, we solve the flow over a forward facing step at Mach 3 with a very coarse mesh to put to
test the capabilites of the SVV shock capturing scheme. The mesh contains 3653 quadrilateral
elements with approximation order~$N=7$ and, to make the case more challenging, there is no
intentional alignment between the element interfaces and the flow structure. Since we are interested
in properly capturing the strong shocks that appear before the step, we keep using the entropy
stable scheme and use the two--point flux from Chandrashekar~\cite{Chandrashekar2013} with matrix
dissipation~\eqref{eq:dg:f-star:mat-diss} for the inter--element fluxes, and the Guermond--Popov
filtered viscous fluxes for the artificial viscosity term.

\begin{figure}
    \centering
    \includegraphics[width=\textwidth]{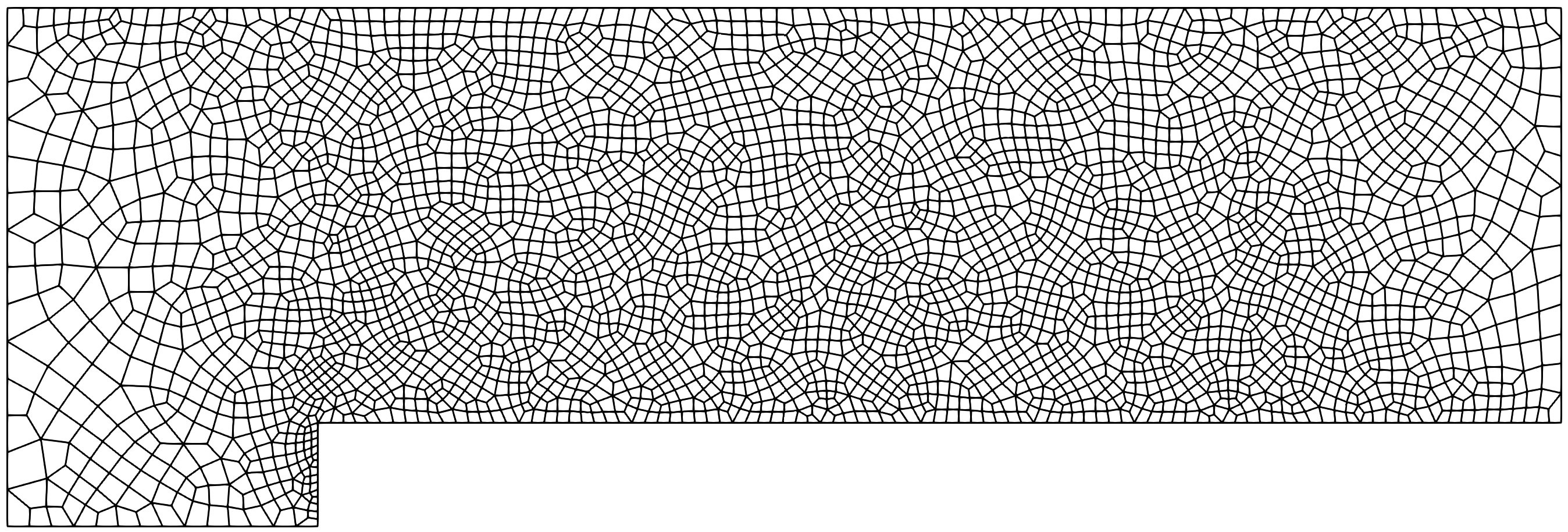}
    \caption{Mesh for the forward facing step case with 3653 elements.}
\end{figure}

As we did in the previous section, we detect the discontinuities with the sensor~$s_2$, defined as
an extension of~\eqref{eq:shu-osher:sensor} to two dimensions,
\begin{equation}
    s_2 = \sqrt{\sum_{i,j=0}^N w_i w_j\left(\svec{\nabla}\rho\right)_{ij}^2},
\end{equation}
and set~$\mu_1 = \mu_2 = \alpha_1 = 0.0005$ and~$\alpha_2 = 0$ according to
\begin{itemize}
    \item $s_1 < 1$: $P_{\SVV}=4$,~$\mu_{\alpha}=\mu_1$ and~$\alpha_{\alpha}=\alpha_1$,
    \item $s_1 \geqslant 1$: $P_{\SVV}=0$,~$\mu_{\alpha}=\mu_2$ and~$\alpha_{\alpha}=\alpha_2$.
\end{itemize}

\begin{figure}
    \centering
    \includegraphics[width=\textwidth]{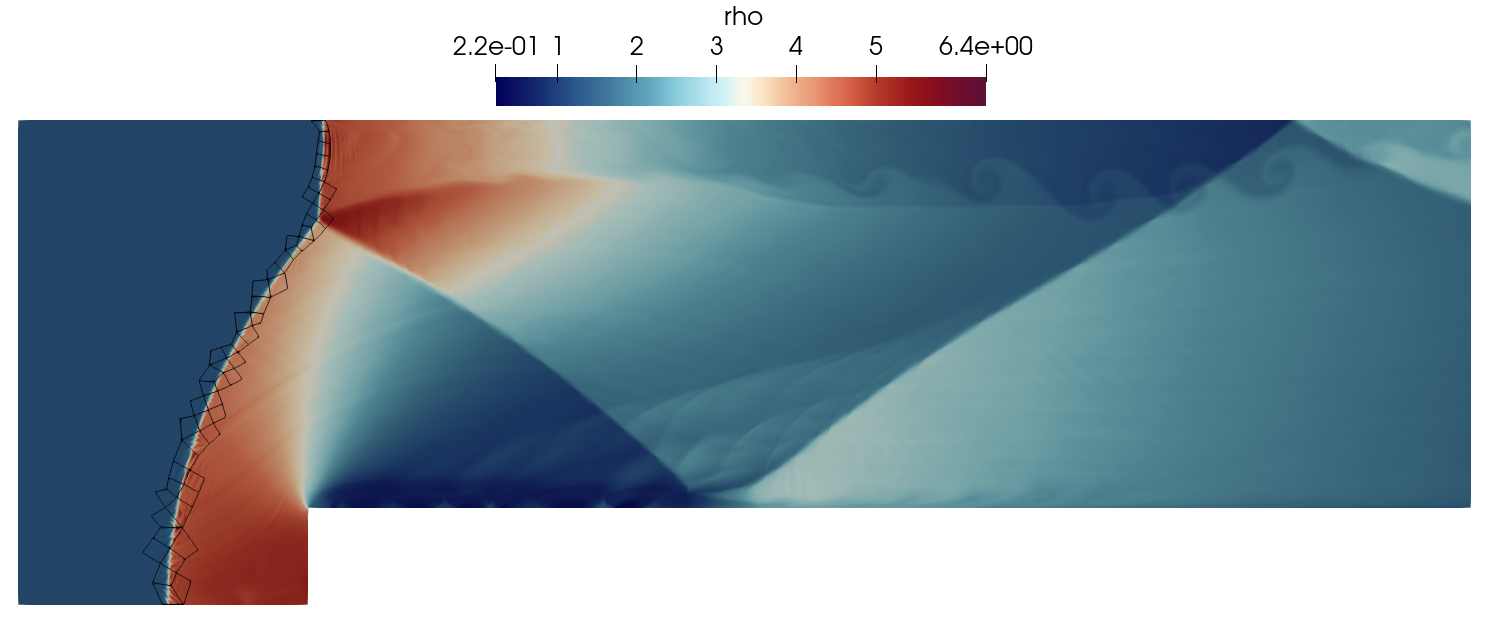}
    \caption{Density contour of the two--dimensional forward facing step at $t=10$.}
    \label{fig:FFS:density}
\end{figure}

We use a short run (until~$t=1$) with higher viscosity
($\mu_1=\mu_2=\alpha_1=0.001$ and~$\alpha_2=0$) to start the case with the entropy stable SVV, and
represent the solution at~$t=10$ in Fig.~\ref{fig:FFS:density}. The thickness of the main shock is
well resolved within the size of one element and, at the same time, the smoother features of the
flow, such as the turbulent wakes that start at the mixing layer of the top or the corner of the
step, are also captured.

\section{Conclusions}

In this work we have proposed an entropy stable SVV filtered artificial viscosity for the DGSEM that
is able to tackle some of the usual drawbacks of the discontinuous Galerkin methods. The revision of
the classical SVV filtering technique has allowed us to prove its entropy stability after small
modifications. As it is usual in SVV methods, the intensity of the dissipation is controlled through
a modulated viscosity coefficient, whereas its wavenumber profile is specified by the filter kernel.
Moreover, we also include an entropy stable discretization of the Guermond--Popov
fluxes~\cite{Guermond2014} using the Bassi--Rebay~1 scheme. These fluxes have been filtered with the
SVV technique and applied to solve usual shock capturig benchmark problems. Having accounted for
this properties, we have tested its capabilities regarding two different mathematical entropy
functions: the kinetic energy (used in subsonic flows) and the thermodynamical entropy (used in
supersonic flows).

For the kinetic energy, we have used a turbulent case and we have shown that a specific choice of
the parameters reproduces the expected energy decay for the scales resolved by the DGSEM.
Additionally, we have also found that this approach can be easily coupled with other schemes to
provide better results. Specifically, as in~\cite{Manzanero2020}, we implement an SVV-LES method
that automates the computation of the viscosity and improves the results for a wide range of values
of the SVV exponent, reducing the impact of the choice of~$P_{\SVV}$ on the solution. However, in
contrast to the method presented in~\cite{Manzanero2020}, the approach described herein is entropy
stable.

For the thermodynamic entropy, we have simulated the well known Shu--Osher one--dimensional case and
a supersonic two--dimensional forward facing step, where the SVV approach based on the
Guermond--Popov fluxes~\cite{Guermond2014} has proven to capture strong shocks fully contained
within a single element and, at the same time, resolve detailed features of the flow in smoother
regions.

\section{Acknowledgments}

Andrés Mateo has received funding from the Universidad Polit\'ecnica de Madrid under the Programa
Propio PhD programme. Eusebio Valero acknowledge the funding received from the European Commission
through the Global Fellowship Grant FLOWCID (Grant Agreement-101019137)\revtwo{, and from the
Ministerio de Ciencia, Innovaci\'on y Universidades of Spain under the project SIMOPAIR
(Ref: RTI2018-480 097075-B-I00).}

\appendix

\section{Effect of the filter kernel parameter}
\label{sec:filter_kernel_effect}

The entropy--stable SVV approach developed in this work adds a new parameter,~$P_\SVV$, controlling
the shape of the filter kernel,~$\mathcal{F}$. However, there is no specific strategy to set its
value, and a set of cases must be run to find it. Our simulations indicate that values of~$O(1)$ are
appropriate for supersonic flows, while turbulent cases need a more flat dissipation profile,
with~$P_\SVV \approx O(10^{-2})$ even without the help of the LES approach.

\begin{figure}
    \centering
    \subfigure[$P_\SVV = 0$]{%
        \includegraphics[width=0.49\textwidth]{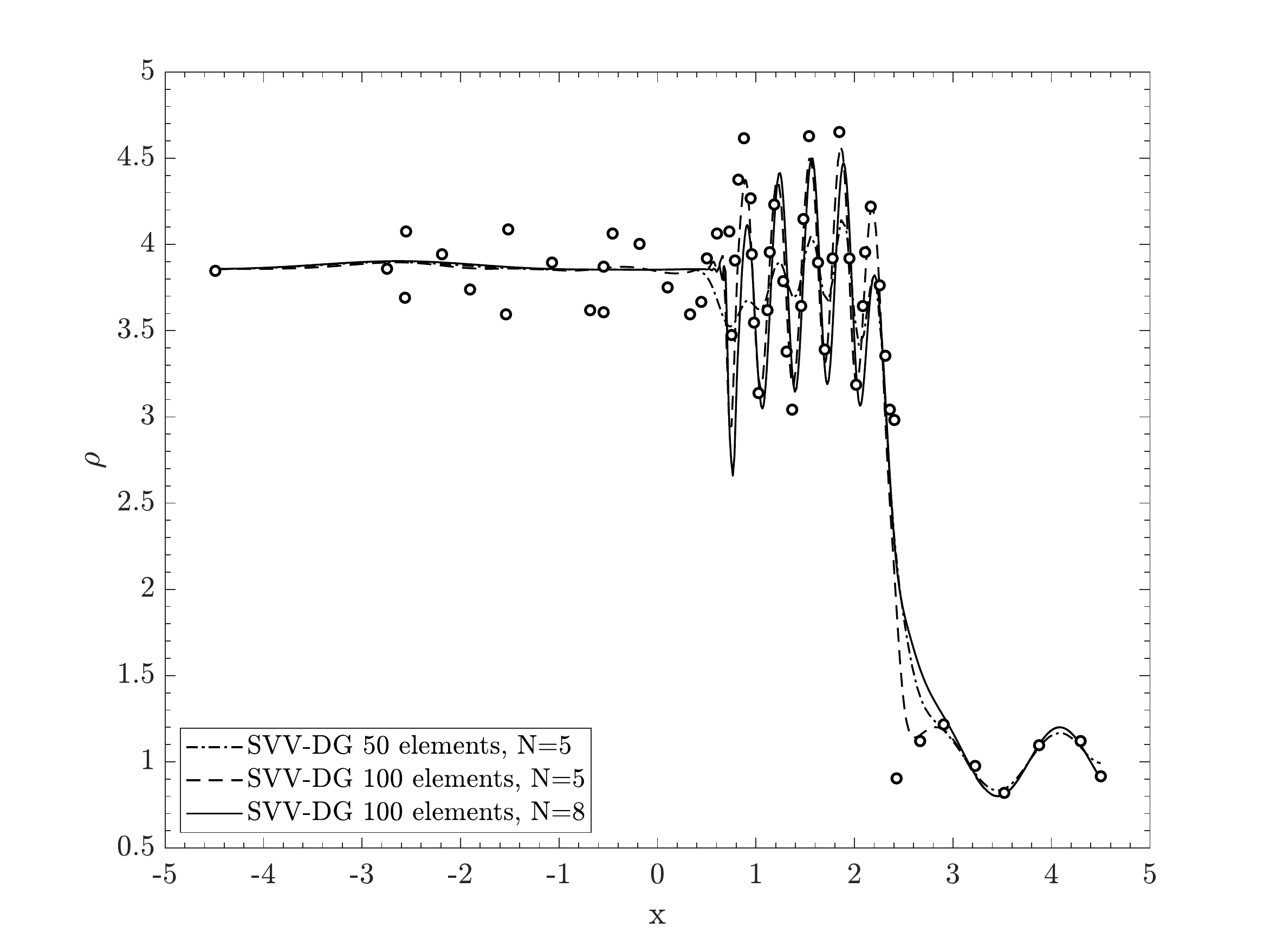}}
    \subfigure[$P_\SVV = 1$]{%
        \includegraphics[width=0.49\textwidth]{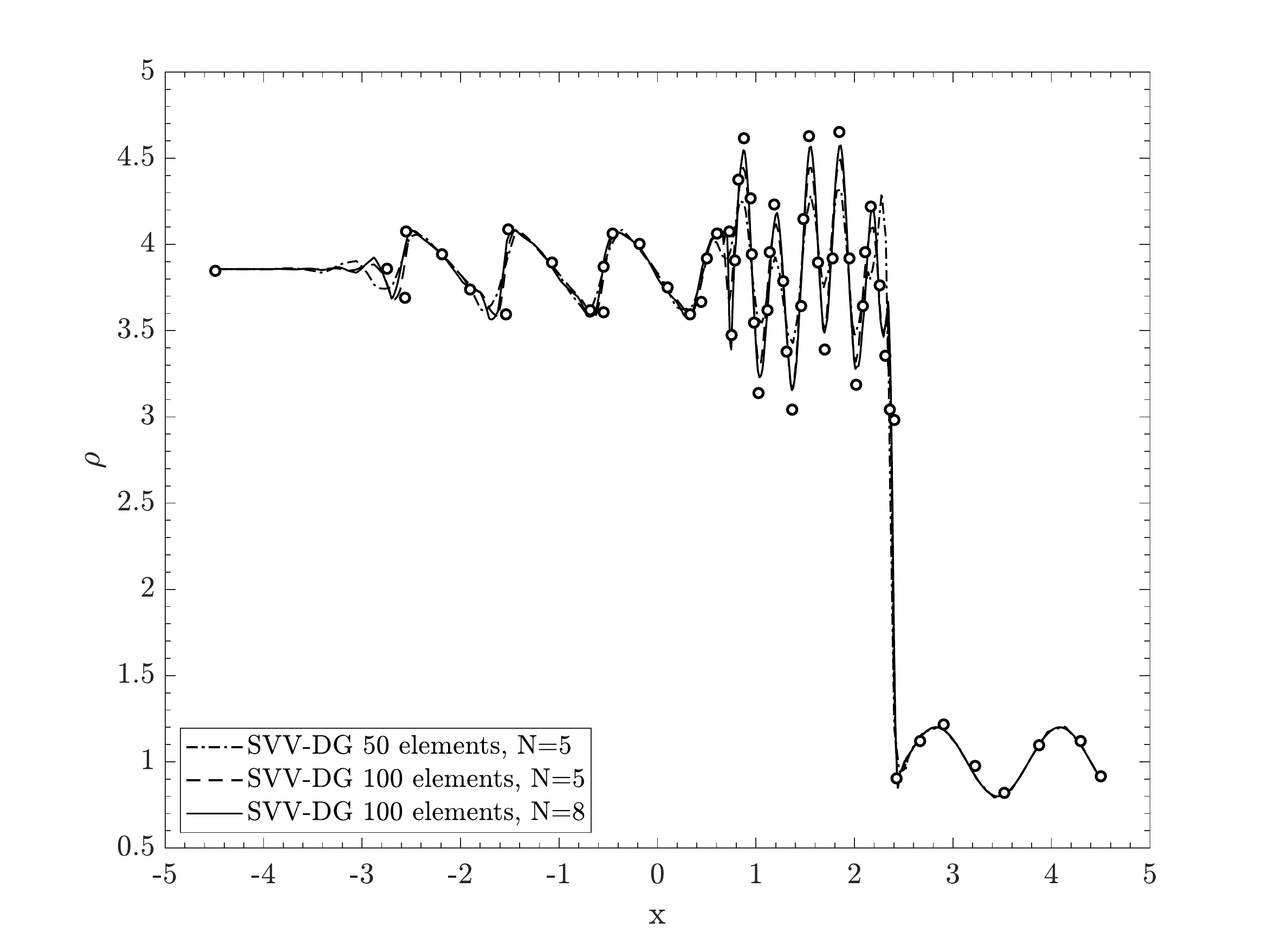}}
    \subfigure[$P_\SVV = 4$]{%
        \includegraphics[width=0.49\textwidth]{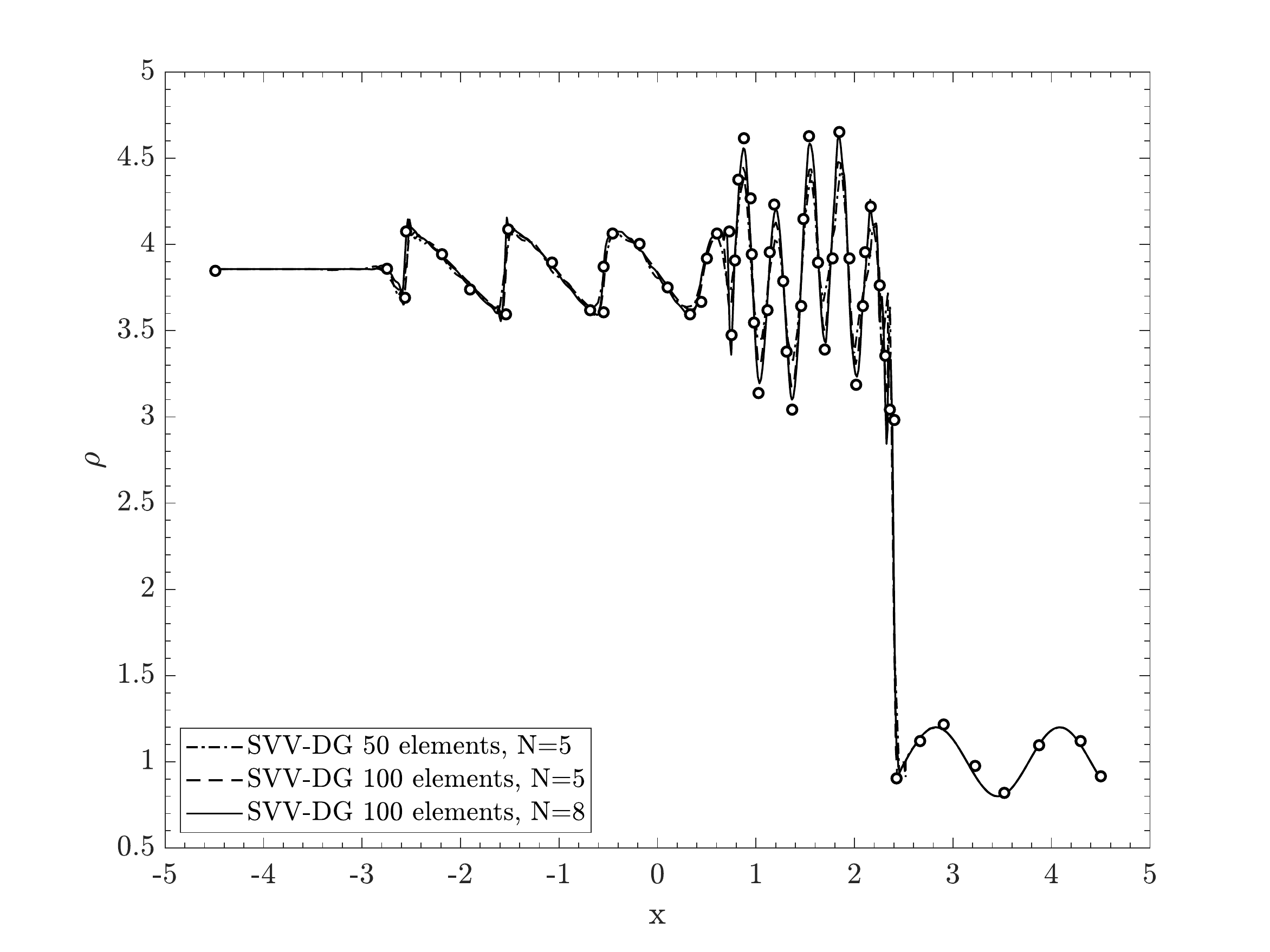}}
    \subfigure[$P_\SVV = 10$]{%
        \includegraphics[width=0.49\textwidth]{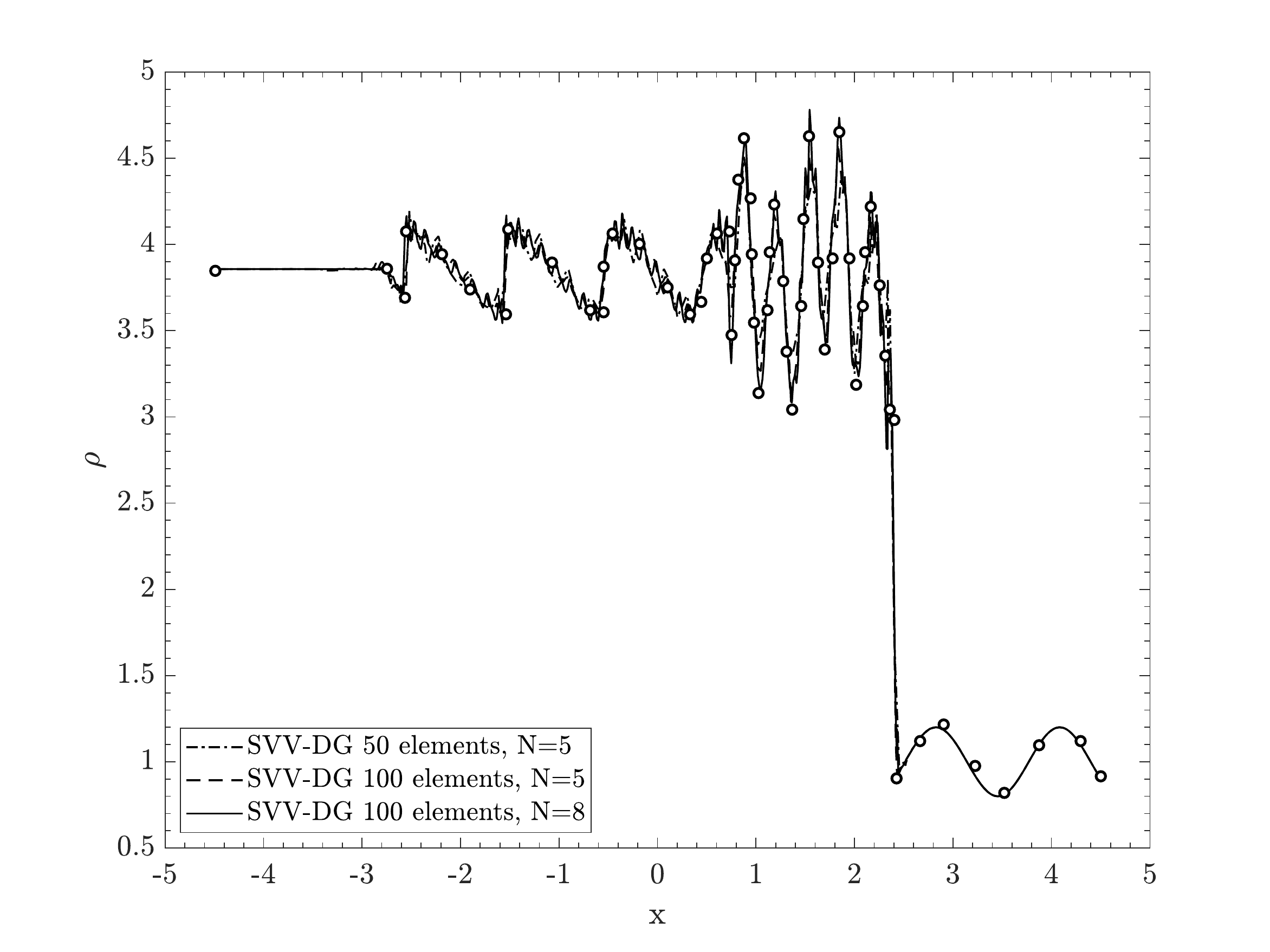}}
    \caption{Density solution of the one--dimensional Shu--Osher problem with different filter
             kernel shapes. The solid lines represent the DG solution with the SVV dissipation
             developed in this work, and the dots represent the reference solution
             of~\cite{shu1989efficient} using an ENO-RF-S-3 scheme with 1600 degrees of freedom.}
    \label{fig:shu-osher:appendix}
\end{figure}

We show in this Appendix (Figure~\ref{fig:shu-osher:appendix}) some one--dimensional test cases
performed with different values of~$P_\SVV$. Although we have chosen the Shu--osher problem in this
section, the effect of the shape parameter is similar in other cases and dimensions.

It is easily seen from Figure~\ref{fig:shu-osher:appendix} that there is only a range of values
for which the solution remains close to the reference, quickly degrading when~$P_\SVV$ is out it.
This loss of accuracy is observed from the first time iterations, being more evident if
discontinuities are present. This makes the trial and error process of finding a ``good'' value
for~$P_\SVV$ somewhat qualitative, but also not especially expensive.

\section{Navier--Stokes viscous flux with kinetic energy variables}
\label{sec:NS_kinetic_matrices}

Kinetic energy stability constitutes a particular case of stability but, since the last entropy
variable~\eqref{eq:NSE:kinetic-energy-variables} is zero, the energy equation is never involved in
the stability analysis. Because of that, we have not included the dissipative terms of the energy
equation into matrix~$\bmat{B}_v^{\mathcal K}$, but we add them separately to the viscous fluxes,
\begin{equation}
    \ssvec{f}_{v} = \mu \bmat{C}_{v}^{\mathcal K}\svec{\nabla}\stvec{w}^{\mathcal K} +
        \left(\tens{\tau}\cdot\svec{u}+\svec{q}\right)\stvec{e}_{5},
\label{eq:app:visc-diss-kin-final}
\end{equation}
with~$\smat{C}_{v,ij}^{\mathcal K} = \smat{C}_{v,ji}^{\mathcal K,T}$ and,
\begin{equation}
    \begin{split}
        &\smat{C}_{v,11}^{\mathcal K} = \left(\begin{array}{ccccc}
            0 & 0 & 0 & 0 & 0 \\
            0 & \frac{4}{3} & 0  & 0 & 0 \\
            0 & 0 & 1 & 0 & 0 \\
            0 & 0 & 0 & 1 & 0 \\
            0 & 0 & 0 & 0 & 0 \\
        \end{array}\right),~~
        \smat{C}_{v,22}^{\mathcal K} = \left(\begin{array}{ccccc}
            0 & 0 & 0 & 0 & 0 \\
            0 & 1 & 0 & 0 & 0 \\
            0 & 0 & \frac{4}{3}  & 0 & 0 \\
            0 & 0 & 0 & 1 & 0 \\
            0 & 0 & 0 & 0 & 0 \\
        \end{array}\right),~~
        \smat{C}_{v,33}^{\mathcal K} = \left(\begin{array}{ccccc}
            0 & 0 & 0 & 0 & 0 \\
            0 & 1 & 0 & 0 & 0 \\
            0 & 0 & 1  & 0 & 0 \\
            0 & 0 & 0 & \frac{4}{3} & 0 \\
            0 & 0 & 0 & 0 & 0 \\
        \end{array}\right),\\
        &\smat{C}_{v,12}^{\mathcal K}=\left(\begin{array}{ccccc} 0 & 0 & 0 & 0 & 0 \\
            0 & 0 & -\frac{2}{3} & 0 & 0 \\
            0 & 1 & 0 & 0 & 0 \\
            0 & 0 & 0 & 0 & 0 \\
            0 & 0 & 0 & 0 & 0 \\
        \end{array}\right),~~\smat{C}_{v,13}^{\mathcal K}=\left(\begin{array}{ccccc}
            0 & 0 & 0 & 0 & 0 \\
            0 & 0 & 0 & -\frac{2}{3} & 0 \\
            0 & 0 & 0  & 0 & 0 \\
            0 & 1 & 0 & 0 & 0 \\
            0 & 0 & 0 & 0 & 0 \\
        \end{array}\right),~~
        \smat{C}_{v,23}^{\mathcal K} = \left(\begin{array}{ccccc}
            0 & 0 & 0 & 0 & 0 \\
            0 & 0 & 0 & 0 & 0 \\
            0 & 0 & 0  & -\frac{2}{3} & 0 \\
            0 & 0 & 1 & 0 & 0 \\
            0 & 0 & 0 & 0 & 0 \\
        \end{array}\right),
    \end{split}
    \label{eq:app:visc-diss-kin}
\end{equation}
such that~$\bmat{C}_v^{\mathcal K}$ is positive definite.

\section{Navier--Stokes viscous flux with thermodynamic entropy variables}
\label{sec:NS_entropy_matrices}

The viscous flux with thermodynamic entropy variables is expressed in terms of a non--linear block
matrix,
\begin{equation}
    \ssvec{f}_v = \frac{\mu p}{\rho} \tilde{\bmat{B}}^{\mathcal S}_v
        \svec{\nabla}\stvec{w}^{\mathcal S},
\end{equation}
where the matrix~$\tilde{\bmat{B}}_v^{\mathcal S}$ is defined as,
\begin{equation}
\begin{split}
    &\smat{\tilde{B}}_{v,11}^{\mathcal S} = \left(\begin{array}{ccccc}
        0 & 0 & 0 & 0 & 0\\
        0 & \frac{4}{3} & 0 & 0 & \frac{4}{3}u \\
        0 & 0 & 1 & 0 & v \\
        0 & 0 & 0 & 1 & w \\
        0 & \frac{4}{3}u & v & w & \frac{1}{3}u^2+v_{\mathrm{tot}}^{2}+\frac{\theta p}{\rho}
    \end{array}\right),~~
    \smat{\tilde{B}}_{v,22}^{\mathcal S} = \left(\begin{array}{ccccc}
        0 & 0 & 0 & 0 & 0\\
        0 & 1 & 0 & 0 & u \\
        0 & 0 & \frac{4}{3} & 0 & \frac{4}{3}v \\
        0 & 0 & 0 & 1 & w \\
        0 & u & \frac{4}{3}v & w & \frac{1}{3}v^2+v_{\mathrm{tot}}^{2}+\frac{\theta p}{\rho}
    \end{array}\right),\\
    &\smat{\tilde{B}}_{v,33}^{\mathcal S} = \left(\begin{array}{ccccc}
        0 & 0 & 0 & 0 & 0\\
        0 & 1 & 0 & 0 & u \\
        0 & 0 & 1 & 0 & v \\
        0 & 0 & 0 & \frac{4}{3} & \frac{4}{3}w \\
        0 & u & v & \frac{4}{3}w & \frac{1}{3}w^2+v_{\mathrm{tot}}^{2}+\frac{\theta p}{\rho}
    \end{array}\right),~~
    \smat{\tilde{B}}_{v,12}^{\mathcal S} = \left(\begin{array}{ccccc}
        0 & 0 & 0 & 0 & 0\\
        0 & 0 & -\frac{2}{3} & 0 & -\frac{2}{3}v \\
        0 & 1 & 0 & 0 & u \\
        0 & 0 & 0 & 0 & 0\\
        0 & v & -\frac{2}{3}u & 0 & \frac{1}{3}uv
    \end{array}\right),\\
    &\smat{\tilde{B}}_{v,13}^{\mathcal S} = \left(\begin{array}{ccccc}
        0 & 0 & 0 & 0 & 0\\
        0 & 0 & 0 & -\frac{2}{3} & -\frac{2}{3}w \\
        0 & 0 & 0 & 0 & 0 \\
        0 & 1 & 0 & 0& u \\
        0 & w & 0 & -\frac{2}{3}u & \frac{1}{3}uw
    \end{array}\right),~~
    \smat{\tilde{B}}_{v,23}^{\mathcal S} = \left(\begin{array}{ccccc}
        0 & 0 & 0 & 0 & 0\\
        0 & 0 & 0 & 0 & 0 \\
        0 & 0 & 0 & -\frac{2}{3} & -\frac{2}{3}w \\
        0 & 0 & 1 & 0 & v\\
        0 & 0 & w & -\frac{2}{3}v & \frac{1}{3}vw
    \end{array}\right),
\end{split}
\label{eq:app:visc-diss-S}
\end{equation}

and can be decomposed as~$\tilde{\bmat{B}}_v^{\mathcal S} =
\bmat{L}_{v}^{\mathcal S,T}\bmat{D}_{v}^{\mathcal S}\bmat{L}_{v}^{\mathcal S}$ with,
\begin{equation}
\begin{split}
    &\bmat{D}_{v}^{\mathcal S} = \mathrm{diag}\left( 0, \frac{4}{3}, 1, 1, \frac{\theta p}{\rho},
        0, 0, 1, 1, \frac{\theta p}{\rho}, 0, 0, 0, 0, \frac{\theta p}{\rho}\right), \\
    &\smat{L}_{v,11}^{S} = \left(\begin{array}{ccccc}
        0 & 0 & 0 & 0 & 0 \\
        0 & 1 & 0  & 0 & u \\
        0 & 0 & 1 & 0 & v \\
        0 & 0 & 0 & 1 & w \\
        0 & 0 & 0 & 0 & 1
    \end{array}\right),~~
    \smat{L}_{v,22}^{S} = \left(\begin{array}{ccccc}
        0 & 0 & 0 & 0 & 0 \\
        0 & 0 & 0  & 0 & 0 \\
        0 & 0 & 1 & 0 & v \\
        0 & 0 & 0 & 1 & w \\
        0 & 0 & 0 & 0 & 1
    \end{array}\right),~~
    \smat{L}_{v,33}^{S} = \left(\begin{array}{ccccc}
        0 & 0 & 0 & 0 & 0 \\
        0 & 0 & 0  & 0 & 0 \\
        0 & 0 & 0 & 0 & 0 \\
        0 & 0 & 0 & 0 & 0 \\
        0 & 0 & 0 & 0 & 1
    \end{array}\right),\\
    &\smat{L}_{v,12}^{\mathcal S} = \left(\begin{array}{ccccc}
        0 & 0 & 0 & 0 & 0 \\
        0 & 0 & -\frac{1}{2} & 0 & \frac{v}{2} \\
        0 & 1 & 0 & 0 & u \\
        0 & 0 & 0 & 0 & 0 \\
        0 & 0 & 0 & 0 & 0
    \end{array}\right),~~
    \smat{L}_{v,13}^{\mathcal S} = \left(\begin{array}{ccccc}
        0 & 0 & 0 & 0 & 0 \\
        0 & 0 &  0 & -\frac{1}{2} & \frac{w}{2} \\
        0 & 0 & 0 & 0 & 0 \\
        0 & 1 & 0 & 0 & u \\
        0 & 0 & 0 & 0 & 0
    \end{array}\right),~~
    \smat{L}_{v,23}^{\mathcal S} = \left(\begin{array}{ccccc}
        0 & 0 & 0 & 0 & 0 \\
        0 & 0 & 0 & 0 & 0 \\
        0 & 0 &  0 & -1 & -w \\
        0 & 0 & 1 & 0 & v \\
        0 & 0 & 0 & 0 & 0
    \end{array}\right)\\
    &\smat{L}_{21}^{\mathcal S} = 0,~~
    \smat{L}_{31}^{\mathcal S} = 0,~~
    \smat{L}_{32}^{\mathcal S} = 0.
\end{split}
\end{equation}

\section{Guermond--Popov viscous flux with thermodynamic entropy variables}
\label{sec:GP_entropy_matrices}

The artificial viscosity of Guermond and Popov~\cite{Guermond2014} that we showed
in~\eqref{eq:NSE:GP-fluxes} can be rewritten in matrix form when expressed in terms of the entropy
variables,
\begin{equation}
    \ssvec{f}_{\mathrm{GP}}^{\mathcal S} =
        \bmat{B}_{\mathrm{GP}}\svec{\nabla}\stvec{w}^{\mathcal S},~~
    \bmat{B}_{\mathrm{GP}} = \alpha_a\rho\bmat{B}_{\mathrm{GP}}^{\alpha} +
        \mu_a p\bmat{B}_{\mathrm{GP}}^{\mu},
\end{equation}
where both block matrices,~$\bmat{B}_{\mathrm{GP}}^{\alpha}$ and~$\bmat{B}_{\mathrm{GP}}^{\mu}$ are,
\begin{subequations}
\begin{equation}
    \smat{B}_{\mathrm{GP},ii}^{\alpha} = \frac{1}{\rho^2}\stvec{q}\stvec{q}^{T} +
        \Lambda^{2}\stvec{e}_{5}\stvec{e}_{5}^{T} = \left(\begin{array}{ccccc}
            1 & u & v & w & e \\
            u & u^2 & uv & uw & ue \\
            v & uv & v^2 & vw & ve \\
            w & uw & vw & w^2 & ve \\
            e & ue & ve & we & e^2 + \Lambda^{2}
        \end{array}\right),\quad
    \smat{B}_{\mathrm{GP},ij}^{\alpha} = 0~(i\neq j),\quad \Lambda = \frac{p/\rho}{\sqrt{\gamma-1}},
\end{equation}
\begin{equation}
\begin{split}
    &\smat{B}_{\mathrm{GP},11}^{\mu} = \left(\begin{array}{ccccc}
        0 & 0 & 0 & 0 & 0 \\
        0 & 1 & 0 & 0 & u \\
        0 & 0 & \frac{1}{2} & 0 & \frac{v}{2} \\
        0 & 0 & 0 & \frac{1}{2} & \frac{w}{2} \\
        0 & u & \frac{v}{2} & \frac{w}{2} & \frac{1}{2}\left(u^2+v_{\mathrm{tot}}^{2}\right)
    \end{array}\right),~~
    \smat{B}_{\mathrm{GP},22}^{\mu} = \left(\begin{array}{ccccc}
        0 & 0 & 0 & 0 & 0 \\
        0 & \frac{1}{2} & 0 & 0 & \frac{u}{2} \\
        0 & 0 & 1 & 0 & v \\
        0 & 0 & 0 & \frac{1}{2} & \frac{w}{2} \\
        0 & \frac{u}{2} & v & \frac{w}{2} & \frac{1}{2}\left(v^2+v_{\mathrm{tot}}^{2}\right)
    \end{array}\right),\\
    &\smat{B}_{\mathrm{GP},33}^{\mu} = \left(\begin{array}{ccccc}
        0 & 0 & 0 & 0 & 0 \\
        0 &\frac{1}{2} & 0 & 0 & \frac{u}{2} \\
        0 & 0 & \frac{1}{2} & 0 & \frac{v}{2} \\
        0 & 0 & 0 & 1 & w \\
        0 & \frac{u}{2} & \frac{v}{2} & w &
            \frac{1}{2}\left(w^2+v_{\mathrm{tot}}^{2}\right)
    \end{array}\right),~~
    \smat{B}_{\mathrm{GP},12}^{\mu} = \left(\begin{array}{ccccc}
        0 & 0 & 0 & 0 & 0 \\
        0 & 0 & 0 & 0 & 0 \\
        0 & \frac{1}{2} & 0 & 0 & \frac{u}{2} \\
        0 & 0 & 0 & 0 & 0 \\
        0 & \frac{v}{2} & 0 & 0 & \frac{uv}{2}
    \end{array}\right),\\
    &\smat{B}_{\mathrm{GP},13}^{\mu} = \left(\begin{array}{ccccc}
        0 & 0 & 0 & 0 & 0 \\
        0 & 0 & 0 & 0 & 0 \\
        0 & 0 & 0 & 0 & 0 \\
        0 & \frac{1}{2} & 0 & 0  & \frac{u}{2} \\
        0 & \frac{w}{2} & 0 & 0 & \frac{uw}{2}
    \end{array}\right),~~
    \smat{B}_{\mathrm{GP},23}^{\mu} = \left(\begin{array}{ccccc}
        0 & 0 & 0 & 0 & 0 \\
        0 & 0 & 0 & 0 & 0 \\
        0 & 0 & 0 & 0 & 0 \\
        0 & 0 & \frac{1}{2} & 0 & \frac{v}{2} \\
        0 & 0 & \frac{w}{2} & 0 & \frac{vw}{2}
    \end{array}\right).
\end{split}
\end{equation}
\end{subequations}
\revtwo{We want to highlight here that, in this case, the density gradient is easily computed as,
\begin{equation}
    \svec{\nabla}\rho = \rho \svec{\nabla}w_{1}^{\mathcal S} +
        \rho u\svec{\nabla}w_{2}^{\mathcal S} + \rho v\svec{\nabla}w_{3}^{\mathcal S} +
        \rho w\svec{\nabla}w_{4}^{\mathcal S} + \rho e\svec{\nabla}w_{5}^{\mathcal S} =
        \stvec{q}^{T}\svec{\nabla}\stvec{w}^{\mathcal S},
\end{equation}
an expression that we have found useful when implementing some of the algorithms derived in this
document as, for instance, the computation of the shock sensor in~\eqref{eq:shu-osher:sensor}.} A
Cholesky decomposition can also be performed to rewrite~$\bmat{B}_{\mathrm{GP}} =
\bmat{L}_{\mathrm{GP}}^T\bmat{D}_{\mathrm{GP}}\bmat{L}_{\mathrm{GP}}$ where,
\begin{equation}
\begin{split}
    &\bmat{D}_{\GP} = \mathrm{diag}\left(\alpha_{a}\rho, \mu_{a}p, \frac{1}{2}\mu_{a}p,
        \frac{1}{2}\mu_{a}p, \alpha_{a} \rho,\alpha_{a} \rho, 0, \mu_{a}p, \frac{1}{2}\mu_{a}p,
        \alpha_{a} \rho, \alpha_{a} \rho, 0, 0, \mu_{a}p, \alpha_{a} \rho\right), \\
    &\smat{L}_{11}^{\GP}=\left(\begin{array}{ccccc}
        1 & u & v & w & e \\
        0 & 1 & 0 & 0 & u \\
        0 & 0 & 1 & 0 & v \\
        0 & 0 & 0 & 1 & w \\
        0 & 0 & 0 & 0 & \Lambda
    \end{array}  \right),\quad
    \smat{L}_{22}^{\GP}=\left(\begin{array}{ccccc}
        1 & u & v & w & e \\
        0 & 0 & 0 & 0 & 0 \\
        0 & 0 & 1 & 0 & v \\
        0 & 0 & 0 & 1 & w \\
        0 & 0 & 0 & 0 & \Lambda
    \end{array}\right),\quad
    \smat{L}_{33}^{\GP}=\left(\begin{array}{ccccc}
        1 & u & v & w & e \\
        0 & 0 & 0 & 0 & 0 \\
        0 & 0 & 0 & 0 & 0 \\
        0 & 0 & 0 & 1 & w \\
        0 & 0 & 0 & 0 & \Lambda
    \end{array}\right), \\
    &\smat{L}_{12}^{\GP}=\left(\begin{array}{ccccc}
        0 & 0 & 0 & 0 & 0 \\
        0 & 0 & 0 & 0 & 0 \\
        0 & 1 & 0 & 0 & u \\
        0 & 0 & 0 & 0 & 0 \\
        0 & 0 & 0 & 0 & 0
    \end{array}\right),\quad
    \smat{L}_{13}^{\GP}=\left(\begin{array}{ccccc}
        0 & 0 & 0 & 0 & 0 \\
        0 & 0 & 0 & 0 & 0 \\
        0 & 0 & 0 & 0 & 0 \\
        0 & 1 & 0 & 0 & u \\
        0 & 0 & 0 & 0 & 0
    \end{array}\right),\quad
    \smat{L}_{23}^{\GP}=\left(\begin{array}{ccccc}
        0 & 0 & 0 & 0 & 0 \\
        0 & 0 & 0 & 0 & 0 \\
        0 & 0 & 0 & 0 & 0 \\
        0 & 0 & 1 & 0 & v \\
        0 & 0 & 0 & 0 & 0
    \end{array}\right), \\
    &\smat{L}_{21}^{\GP} = \smat{L}_{31}^{\GP} = \smat{L}_{32}^{\GP} = 0.
\end{split}
\label{eq:app:GP-Lmat}
\end{equation}

\bibliographystyle{elsarticle-num}
\bibliography{library}

\begin{thebibliography}{10}
\expandafter\ifx\csname url\endcsname\relax
  \def\url#1{\texttt{#1}}\fi
\expandafter\ifx\csname urlprefix\endcsname\relax\def\urlprefix{URL }\fi
\expandafter\ifx\csname href\endcsname\relax
  \def\href#1#2{#2} \def\path#1{#1}\fi

\bibitem{Reed1977}
W.~H. Reed, T.~R. Hill, Triangular mesh methods for the neutron transport
  equation, 1977, p.~23.

\bibitem{Cockburn2000}
B.~Cockburn, G.~E. Karniadakis, C.-W. Shu, The development of discontinuous
  {Galerkin} methods (2000).

\bibitem{Woopen2014}
M.~Woopen, A.~Balan, G.~May, J.~Sch{\"u}tz, A comparison of hybridized and
  standard {DG} methods for target-based hp-adaptive simulation of compressible
  flow, Computers \& Fluids 98 (2014) 3--16.

\bibitem{Kompenhans2016}
M.~Kompenhans, G.~Rubio, E.~Ferrer, E.~Valero, Adaptation strategies for high
  order discontinuous {Galerkin} methods based on tau-estimation, Journal of
  Computational Physics 306 (2016) 216--236.

\bibitem{Friedrich2018}
L.~Friedrich, A.~R. Winters, D.~C. D.~R. Fern{\'a}ndez, G.~J. Gassner,
  M.~Parsani, M.~H. Carpenter, An entropy stable h / p non-conforming
  discontinuous {Galerkin} method with the summation-by-parts property, Journal
  of Scientific Computing 77 (2018) 689--725.

\bibitem{Ntoukas2021}
G.~Ntoukas, J.~Manzanero, G.~Rubio, E.~Valero, E.~Ferrer, A free--energy stable
  p--adaptive nodal discontinuous {Galerkin} for the {Cahn--Hilliard} equation,
  Journal of Computational Physics 442 (2021) 110409.

\bibitem{Toro2009}
E.~F. Toro, Riemann Solvers and Numerical Methods for Fluid Dynamics, Springer
  Berlin Heidelberg, 2009.

\bibitem{Harten1983}
A.~Harten, On the symmetric form of systems of conservation laws with entropy,
  Journal of Computational Physics 49 (1983) 151--164.

\bibitem{Tadmor1986}
E.~Tadmor, A minimum entropy principle in the gas dynamics equations, Applied
  Numerical Mathematics 2 (1986) 211--219.

\bibitem{Tadmor2003}
E.~Tadmor, Entropy stability theory for difference approximations of nonlinear
  conservation laws and related time-dependent problems, Acta Numerica 12
  (2003) 451--512.

\bibitem{Kopriva2009}
D.~A. Kopriva, Implementing Spectral Methods for Partial Differential
  Equations, Springer Netherlands, 2009.

\bibitem{Strand1994}
B.~Strand, Summation by parts for finite difference approximations for d/dx,
  Journal of Computational Physics 110 (1994) 47--67.

\bibitem{Abgrall2020}
R.~Abgrall, J.~Nordström, P.~Öffner, S.~Tokareva, Analysis of the {SBP-SAT}
  stabilization for finite element methods part i: Linear problems, Journal of
  Scientific Computing 85 (2020).

\bibitem{Kravchenko1997}
A.~Kravchenko, P.~Moin, On the effect of numerical errors in large eddy
  simulations of turbulent flows, Journal of Computational Physics 131 (1997)
  310--322.

\bibitem{Hennemann2021}
S.~Hennemann, A.~M. Rueda-Ramírez, F.~J. Hindenlang, G.~J. Gassner, A provably
  entropy stable subcell shock capturing approach for high order split form dg
  for the compressible euler equations, Journal of Computational Physics 426
  (2021) 109935.

\bibitem{Persson2006}
P.-O. Persson, J.~Peraire, Sub-cell shock capturing for discontinuous
  {Galerkin} methods, Collection of Technical Papers - 44th AIAA Aerospace
  Sciences Meeting 2 (2006) 1408--1420.

\bibitem{Klockner2011}
A.~Kl{\"o}ckner, T.~Warburton, J.~S. Hesthaven, Viscous shock capturing in a
  time-explicit discontinuous {Galerkin} method, Mathematical Modelling of
  Natural Phenomena 6 (2011) 57--83.

\bibitem{Gottlieb2001}
D.~Gottlieb, J.~Hesthaven, Spectral methods for hyperbolic problems, Journal of
  Computational and Applied Mathematics 128 (2001) 83--131.

\bibitem{Guermond2014}
J.~L. Guermond, B.~Popov, Viscous regularization of the {Euler} equations and
  entropy principles, SIAM Journal on Applied Mathematics 74 (2014) 284--305.

\bibitem{Karamanos2000}
G.-S. Karamanos, G.~E. Karniadakis, A spectral vanishing viscosity method for
  large-eddy simulations, Journal of Computational Physics 163 (2000) 22--50.

\bibitem{Kirby2006}
R.~M. Kirby, S.~J. Sherwin, Stabilisation of spectral/hp element methods
  through spectral vanishing viscosity: Application to fluid mechanics
  modelling, Computer Methods in Applied Mechanics and Engineering 195 (2006)
  3128--3144.

\bibitem{Moura2016}
R.~Moura, S.~Sherwin, J.~Peir{\'o}, Eigensolution analysis of spectral/hp
  continuous {Galerkin} approximations to advection--diffusion problems:
  Insights into spectral vanishing viscosity, Journal of Computational Physics
  307 (2016) 401--422.

\bibitem{Lodares2021}
D.~Lodares, J.~Manzanero, E.~Valero, An entropy–stable discontinuous
  {G}alerkin approximation of the {S}palart--{A}llmaras turbulence model for
  the compressible {R}eynolds {A}veraged {N}avier--{S}tokes equations, Under
  review in Journal of Computational Physics (2021).

\bibitem{Manzanero2020}
J.~Manzanero, E.~Ferrer, G.~Rubio, E.~Valero, Design of a {Smagorinsky}
  spectral vanishing viscosity turbulence model for discontinuous {Galerkin}
  methods, Computers \& Fluids 200 (2020) 104440.

\bibitem{Kirby2002}
R.~M. Kirby, G.~E. Karniadakis, Coarse resolution turbulence simulations with
  spectral vanishing viscosity-large-eddy simulations ({SVV-LES}), Journal of
  Fluids Engineering 124 (2002) 886--891.

\bibitem{Pasquetti2008}
R.~Pasquetti, E.~S{\'e}verac, E.~Serre, P.~Bontoux, M.~Sch{\"a}fer, From
  stratified wakes to rotor-stator flows by an {SVV-LES} method, Comput. Fluid
  Dyn 22 (2008) 261--273.

\bibitem{Gassner2018}
G.~J. Gassner, A.~R. Winters, F.~J. Hindenlang, D.~A. Kopriva, The {BR1} scheme
  is stable for the compressible {Navier-Stokes} equations, Journal of
  Scientific Computing 77 (2018) 154--200.

\bibitem{Friedrichs1971}
K.~O. Friedrichs, P.~D. Lax, Systems of conservation equations with a convex
  extension, Proceedings of the National Academy of Sciences 68 (1971)
  1686--1688.

\bibitem{Merriam1989}
M.~L. Merriam, An entropy-based approach to nonlinear stability (1989).

\bibitem{Fisher2013}
T.~C. Fisher, M.~H. Carpenter, High-order entropy stable finite difference
  schemes for nonlinear conservation laws: Finite domains, Journal of
  Computational Physics 252 (2013) 518--557.

\bibitem{Jameson2008}
A.~Jameson, Formulation of kinetic energy preserving conservative schemes for
  gas dynamics and direct numerical simulation of one-dimensional viscous
  compressible flow in a shock tube using entropy and kinetic energy preserving
  schemes, Journal of Scientific Computing 34 (2008) 188--208.

\bibitem{Gassner2014}
G.~J. Gassner, A kinetic energy preserving nodal discontinuous galerkin
  spectral element method, International Journal for Numerical Methods in
  Fluids 76 (2014) 28--50.

\bibitem{gassner2016split}
G.~J. Gassner, A.~R. Winters, D.~A. Kopriva, Split form nodal discontinuous
  galerkin schemes with summation-by-parts property for the compressible euler
  equations, Journal of Computational Physics 327 (2016) 39--66.

\bibitem{Kopriva2006}
D.~A. Kopriva, Metric identities and the discontinuous spectral element method
  on curvilinear meshes, Journal of Scientific Computing 26 (2006) 301--327.

\bibitem{Winters2021}
A.~R. Winters, D.~A. Kopriva, G.~J. Gassner, F.~Hindenlang, Construction of
  modern robust nodal discontinuous galerkin spectral element methods for the
  compressible navier–stokes equations (2021).

\bibitem{Kopriva2017}
D.~A. Kopriva, A polynomial spectral calculus for analysis of dg spectral
  element methods, arXiv preprint arXiv:1704.00709 (2017).

\bibitem{gassner2013skew}
G.~J. Gassner, A skew-symmetric discontinuous galerkin spectral element
  discretization and its relation to sbp-sat finite difference methods, SIAM
  Journal on Scientific Computing 35~(3) (2013) A1233--A1253.

\bibitem{Carpenter1999}
M.~H. Carpenter, J.~Nordström, D.~Gottlieb, A stable and conservative
  interface treatment of arbitrary spatial accuracy, Journal of Computational
  Physics 148 (1999) 341--365.

\bibitem{Nordstrom2013}
J.~Nordström, T.~Lundquist, Summation-by-parts in time, Journal of
  Computational Physics 251 (2013) 487--499.

\bibitem{DelRey2014}
D.~C. D.~R. Fernández, J.~E. Hicken, D.~W. Zingg, Review of summation-by-parts
  operators with simultaneous approximation terms for the numerical solution of
  partial differential equations, Computers \& Fluids 95 (2014) 171--196.

\bibitem{Pirozzoli2010}
S.~Pirozzoli, Generalized conservative approximations of split convective
  derivative operators, Journal of Computational Physics 229 (2010) 7180--7190.

\bibitem{Chandrashekar2013}
P.~Chandrashekar, Kinetic energy preserving and entropy stable finite folume
  schemes for compressible {Euler} and {Navier-Stokes} equations,
  Communications in Computational Physics 14 (2013) 1252--1286.

\bibitem{BR1}
{\relax F. Bassi and S. Rebay}, A high-order accurate discontinuous finite
  element method for the numerical solution of the compressible
  {N}avier-{S}tokes equations, Journal of Computational Physics 131~(2) (1997)
  267 -- 279.

\bibitem{Carlson2011}
J.-R. Carlson, Inflow/outflow boundary conditions with application to {FUN3D}
  (2011).

\bibitem{Mengaldo2014}
G.~Mengaldo, D.~D. Grazia, F.~Witherden, A.~Farrington, P.~Vincent, S.~Sherwin,
  J.~Peiro, A guide to the implementation of boundary conditions in compact
  high-order methods for compressible aerodynamics, American Institute of
  Aeronautics and Astronautics, 2014.

\bibitem{Tadmor1989}
E.~Tadmor, Convergence of spectral methods for nonlinear conservation laws,
  SIAM Journal on Numerical Analysis 26 (1989) 30--44.

\bibitem{Maday1993}
Y.~Maday, S.~M.~O. Kaber, E.~Tadmor, Legendre pseudospectral viscosity method
  for nonlinear conservation laws, SIAM Journal on Numerical Analysis 30 (1993)
  321--342.

\bibitem{Kaber1996}
S.~M.~O. Kaber, A legendre pseudospectral viscosity method, Journal of
  Computational Physics 128 (1996) 165--180.

\bibitem{Nordstrom2021}
J.~Nordström, A.~R. Winters, Stable filtering procedures for nodal
  discontinuous galerkin methods, Journal of Scientific Computing 87 (2021)
  1--9.

\bibitem{Lundquist2020}
T.~Lundquist, J.~Nordström, Stable and accurate filtering procedures, Journal
  of Scientific Computing 82 (2020) 16.

\bibitem{hindenlang2020stability}
F.~J. Hindenlang, G.~J. Gassner, D.~A. Kopriva, Stability of wall boundary
  condition procedures for discontinuous galerkin spectral element
  approximations of the compressible euler equations, Spectral and High Order
  Methods for Partial Differential Equations (2020) 3.

\bibitem{Chavez2018}
M.~Chávez-Modena, E.~Ferrer, G.~Rubio, Improving the stability of
  multiple-relaxation lattice boltzmann methods with central moments, Computers
  \& Fluids 172 (2018) 397--409.

\bibitem{solan2021application}
P.~Sol{\'a}n-Fustero, A.~Navas-Montilla, E.~Ferrer, J.~Manzanero,
  P.~Garc{\'\i}a-Navarro, Application of approximate dispersion-diffusion
  analyses to under--resolved {B}urgers turbulence using high resolution {WENO}
  and {UWC} schemes, Journal of Computational Physics 435 (2021) 110246.

\bibitem{kou2021eigensolution}
J.~Kou, A.~Hurtado-de Mendoza, S.~Joshi, S.~L. Clainche, E.~Ferrer,
  Eigensolution analysis of immersed boundary method based on volume
  penalization: applications to high-order schemes, arXiv preprint
  arXiv:2107.10155 (2021).

\bibitem{Pope2001}
S.~B. Pope, Turbulent flows, Measurement Science and Technology 12 (2001)
  2020--2021.

\bibitem{Rees2011}
W.~M. van Rees, A.~Leonard, D.~I. Pullin, P.~Koumoutsakos, A comparison of
  vortex and pseudo-spectral methods for the simulation of periodic vortical
  flows at high reynolds numbers, Journal of Computational Physics 230 (2011)
  2794--2805.

\bibitem{shu1989efficient}
C.-W. Shu, S.~Osher, Efficient implementation of essentially non-oscillatory
  shock-capturing schemes, ii, in: Upwind and High-Resolution Schemes,
  Springer, 1989, pp. 328--374.

\end{thebibliography}

\end{document}